\documentclass[onefignum,onetabnum,reqno]{siamart220329}

\usepackage{a4}
\usepackage{amsmath}
\usepackage{amsfonts}
\usepackage{graphicx}
\usepackage{longtable}
\usepackage{subcaption}
\usepackage{algorithmic}
\usepackage{algorithm}
\usepackage{multirow}
\usepackage{tikz}
\usepackage{overpic}
\usepackage{wrapfig}
\usetikzlibrary{arrows,decorations.pathmorphing,backgrounds,positioning,fit,petri,ext.paths.ortho}
\usetikzlibrary{shapes.geometric}
\floatstyle{plain}
\restylefloat{algorithm}
\usepackage{url}
\usepackage{hyperref}
\graphicspath{{figures/}{surffigures/}}

\usepackage{geometry}
\makeatother

\newcommand{\Real}{\mathbb R}

\newcommand{\set}[1]{\left\{ \, #1 \, \right\}}
\newcommand{\seq}[1]{\left<#1\right>}

\newcommand{\restr}[2]{\kern-\nulldelimiterspace % automatically resize the bar with \right
  #1 % the function
\ensuremath{{\mathbin\upharpoonright}%
\raise-.5ex\hbox{$#2$}}}

\def\X {{\mathbf{X}}}

\def\Ai {{\mathbf{A_i}}}
\def\Aie {{\mathbf{A_{ie}}}}
\def\M {{\mathbf{M}}}

\def\I {{\mathbf{I}}}

\def\u {{\mathbf{u}}}
\def\v {{\mathbf{v}}}
\def\w {{\mathbf{w}}}

\def\c {{\mathbf{c}}}
\def\wvec {{\vec{w}}}
\def\cvec {{\vec{c}}}

\theoremstyle{remark} %style : remark (italic title, roman body) 

\theoremstyle{definition} %style : definition (fat title, roman body) 

\theoremstyle{plain} %style : plain (fat title, italic body) 

\usepackage[textsize=scriptsize]{todonotes}
\usepackage[xcolor=dvipdf]{changes}

\tikzstyle{startstop} = [rectangle, rounded corners, minimum width=3cm, minimum height=0.8cm,text width=5cm,text centered, draw=black]
\tikzstyle{startstop1} = [rectangle, rounded corners, minimum width=3cm, minimum height=0.8cm,text centered, draw=black]
\tikzstyle{process} = [rectangle, minimum width=3cm, minimum height=0.8cm, text centered, text width=5cm, draw=black]
\tikzstyle{process1} = [rectangle, minimum width=3cm, minimum height=0.8cm, text centered, text width=5.5cm, draw=black]
\tikzstyle{arrow} = [thick,->,>=stealth]
\headers{Novel bidomain partitioned strategies}{Gopika P.B.,  Peter B. \& Nagaiah C.}

\title{Novel bidomain partitioned strategies for the simulation
 of ventricular fibrillation dynamics} % on realistic 3D geometry}
\author{Gopika P B\footnotemark[3] \thanks{School of Mathematics, IISER Thiruvananthapuram, Kerala, India (\email{gopikapb23@iisertvm.ac.in}, \email{nagaiah.chamakuri@iisertvm.ac.in})} \and Peter Bastian\thanks{Interdisciplinary Center for Scientific Computing, Heidelberg University, Heidelberg, Germany (\email{peter.bastian@iwr.uni-heidelberg.de})} \and Nagaiah Chamakuri\footnotemark[1] \thanks{Center for High-Performance Computing, IISER Thiruvananthapuram, Kerala, India}
}
\date{\today}
\usepackage{amsopn}

\begin{document}
\maketitle

\begin{abstract}
The numerical tools to simulate the bidomain model in cardiac electrophysiology are constantly developing due to
the great clinical interest and scientific advances in mathematical models and computational power.
The bidomain model consists of an elliptic partial differential equation (PDE)
and a non-linear parabolic PDE of reaction-diffusion type,
where the reaction term is described by a set of ordinary differential equations (ODEs).
We propose and analyze a suite of numerical strategies for the efficient and accurate simulation of cardiac electrophysiology, with a particular focus on ventricular fibrillation in realistic geometries. Specifically, we develop and compare a fully coupled strategy, a traditional decoupled strategy, and a novel partitioned strategy. The centerpiece of this work is a bidomain partitioned strategy enhanced with spectral deferred correction, designed to balance numerical stability and computational efficiency.
To address the substantial memory requirements posed by biophysically detailed ionic models, we adopt a compile-time memory-efficient sparse matrix technique. This enables the efficient solution of the coupled nonlinear parabolic PDE and the associated large systems of ODEs that govern ionic gating and concentration dynamics. We perform comprehensive numerical experiments using the Luo–Rudy and Ten Tusscher cell models in both two- and three-dimensional geometries. In addition, we demonstrate the applicability of our approach to bidomain-bath coupling scenarios. The results confirm that the proposed partitioned strategy achieves high accuracy and efficiency compared to standard decoupled strategies.
Our findings support the use of advanced partitioned strategies for large-scale simulations in cardiac electrophysiology and highlight their potential for future investigations into cardiac arrhythmias and other pathological conditions.
%\commentnc{will be modified at the end}
\end{abstract}

% REQUIRED

% REQUIRED

% \commentgp{add your comments like this}

%%%%%%%%%%%%%%%%%%%%%%%%%%%%%%%%%%%%%%%%%%%%%%%%%%%%%
%%%%%%%%%%%%%%%%%%%%%%%%%%%%%%%%%%%%%%%%%%%%%%%%%%%%%
\section{Introduction}
%%%%%%%%%%%%%%%%%%%%%%%%%%%%%%%%%%%%%%%%%%%%%%%%%%%%%
%%%%%%%%%%%%%%%%%%%%%%%%%%%%%%%%%%%%%%%%%%%%%%%%%%%%%

Over the past few decades, computer simulations of non-invasive cardiac electrical outputs, such as electrocardiograms (ECGs) and body surface potential maps, have played a crucial role in advancing our understanding of cardiac diseases like atrial fibrillation \cite{Zhou_PLOSOne16}, ventricular tachycardia \cite{Lopez-Bishop-FPhy19}, and cardiac ischemia \cite{Philip-Gernot-CBM20}. The success of these simulations is largely attributed to their ability to accurately model both the cardiac and torso geometries, along with incorporating detailed, physics-based models that capture the complex biological processes involved. Efforts to close the gap between real-world (in vivo or in vitro) observations and computational (in silico) simulations \cite{Vadakkumpadan-IEEETMI12, Bishop-ajpheart10} have led to the creation of highly detailed and anatomically complex heart-torso models. However, achieving high anatomical accuracy, along with sophisticated mathematical models of cardiac physiology, demands a fine spatial resolution in the discretization of these geometries \cite{Niederer-PTRS11}. This increased level of detail, while necessary for precision, significantly raises computational demands and costs.

The bidomain model is a key tool in simulating the electrical activity of the heart. It consists of a set of strongly coupled partial differential equations that describe the diffusion of intracellular and extracellular potential fields across cardiac tissue. Additionally, ionic currents through the cell membrane and the associated cellular reactions are represented by systems of ordinary differential equations (ODEs), which are based on either phenomenological or electrophysiological cell models \cite{tung78:_bidomain}. 
Solving these models numerically requires substantial computational time, often limiting simulations to just a few heartbeats for the ventricles and atria. This is a significant challenge for studying long-term cardiac function.
Given these limitations, much effort has been directed toward improving the speed of cardiac electrophysiology simulations. The complexity of the problem arises from the multiple temporal and spatial scales involved in wave excitation across the membrane. As a result, solving the bidomain equations numerically is highly challenging
 \cite{Vigmond08:_solvers, Southern_TBE09, Pathmanathan_PBMB10, ColliFranzone_book14} and references therein.
Furthermore, longer observation times are necessary to study wavefront propagation without boundary effects and to analyze events like arrhythmias, e.g., reentry.

%\commentnc{Why implicit schemes are preferred?}
A widely adopted strategy for simulating complex cardiac electrophysiology models is operator splitting (OS). This approach involves decoupling the fully coupled, strongly nonlinear system of partial differential equations (PDEs) into simpler subproblems that can be solved more efficiently, often with the help of specialized numerical techniques (see \cite{Vigmond08:_solvers}). OS methods have been successfully employed for several decades and continue to be a significant area of active research due to their effectiveness in handling computationally challenging problems. In the context of the bidomain model, numerous methods have been proposed to solve its equations. The finite element method (FEM) has proven particularly effective for the spatial discretization of these equations. When the bidomain model is discretized spatially using FEM, the result is a system of differential-algebraic equations (DAEs). These can be integrated in time either using coupled or decoupled strategies.
In coupled approaches, all state variables in the DAE system are updated simultaneously—typically using a fully implicit time-stepping scheme \cite{Murillo_NLA04, ColliFranzone_ESIAM-M2AN13, Nagaiah_JMB13}. This allows for greater numerical stability and accuracy but comes at a high computational cost due to the need to solve a nonlinear system at each time step. On the other hand, decoupled approaches such as implicit-explicit (IMEX) time-stepping \cite{Boulakia_AMRe08, ColliFranzone_MBS05, Ethier_Bourgault_SINUM08, Pathmanathan_PBMB10, Plank_AMG_TBE07, Qu_TBE99, Southern_TBE09, Whiteley_TBE06} or operator splitting methods \cite{Trangenstein_JCP04, Sundnes_MBS05, MacLachlan_Sundnes_CMBBE07, Southern_TBE09, Plank_AMG_TBE07} treat the reaction and diffusion terms separately. This separation simplifies the numerical process, requiring only the solution of linear systems for the PDE variables at each time step. Although coupled methods tend to provide higher accuracy, they are significantly more computationally intensive and demand greater memory resources. This is especially limiting in realistic cardiac simulations, where detailed physiological cell models can result in millions of degrees of freedom. Consequently, coupled solvers with fully implicit time-stepping have so far only been applied to simplified phenomenological models with relatively few state variables.

We can see that the semi-discretized systems from spatial discretization of the bidomain model, which is under consideration, have forms similar to
those of index-1 differential algebraic equations (DAEs)\cite{brenan1995numerical}.
For differential-algebraic equations, which combine algebraic constraints
with ODEs, splitting schemes separate the algebraic constraints from the differential equations. For example, when the ODEs and constraints arise from distinct but coupled physical phenomena, splitting schemes
can take full advantage of existing computer codes tuned for each subproblem \cite{strang1968construction,dae1_cdcs_2005}.
It is well known that, in the numerical solution of ODEs, the method of OS yields 
high-order accurate schemes based on separate, computationally convenient treatments of distinct physical effects. Such schemes are equally desirable but much less accurate for semi-explicit index-1 differential-algebraic 
equations (DAEs). In addition, in \cite{dae1_cdcs_2005}, it is shown that naive application of standard splitting schemes to DAEs suffers from order reduction in both first and second-order schemes, which are only first-order accurate for DAEs. 
In contrast, the spectral deferred corrections (SDC) \cite{dutt2000spectral, MR2216607} combined with OS strategies achieve the desired second-order accuracy and improve the accuracy for index-1 DAEs \cite{yoon2021spectral}.

This article examines various time-stepping methods: explicit, semi-implicit, operator splitting, coupled, and deferred correction. Due to varying naming conventions for operator splitting strategies, we use three terms in this article. The first is the \textsl{fully coupled (or monolithic) strategy}, where the complete bidomain model equations are solved in a single step at each time step of implicit schemes. This strategy resolves the instability issues of explicit and decoupled methods, offering greater robustness and accuracy. However, it requires more computational time due to the use of the Newton method and necessitates a strong preconditioner to solve the ill-conditioned elliptic system efficiently. The second is the \textsl{partitioned strategy}, introduced in this work, where we propose novel partitioning methods for the bidomain and bidomain-bath model equations to enhance accuracy. The third is the \textsl{decoupled strategy}, well-established in cardiac simulation research. We present results using several popular decoupled strategies from existing literature.

In this paper, we introduce a novel partitioned strategy for efficiently solving the bidomain model, complementing both fully coupled and decoupled strategies. This strategy simultaneously solves the parabolic equation alongside the set of ODEs. The complete implementation details are provided in Section~\ref{sec:PartitionedStrategy}.  The fully coupled approach to solving the elliptic PDE, parabolic PDE, and ODEs demands significant memory due to the large size of the physiological cell models. To address this challenge, we employ a memory-efficient compile-time sparse matrix technique (CTSM) \cite{NP_AMC19} in the FEM simulations. The CTSM technique optimizes the assembly of the global compressed row storage (CRS) matrix for the finite element discretization. By leveraging the loosely coupled nature of the cellular state and gating variables, we implement local block matrices for each CRS element as compile-time sparse matrices. Additionally, we utilize higher-order time-stepping schemes with adaptive step-size selection, substantially reducing computational costs and memory requirements. This makes it feasible to solve the bidomain equations in a coupled manner, even with complex physiological cell models, and enables accurate simulations of pathological phenomena such as cardiac reentry. In the latter half of this paper, we develop a new family of higher-order splitting schemes based on the spectral deferred correction (SDC) method for index-1 DAEs. These schemes follow the spectral deferred correction paradigm \cite{dutt2000spectral, MR2216607}, in which an error equation is numerically solved, offering a computationally efficient structure. The partitioned strategy, when coupled with the SDC method, achieves higher-order accuracy. Numerical results validate the anticipated order of accuracy and demonstrate the efficiency of the proposed approach.

The remaining part of the manuscript is organized as follows. In the next section, we present the bidomain model 
equations and list the cell models, Luo-Rudy \cite{Luo_CircRes91} and Ten Tusscher-Panfilov \cite{tenTusscher_AJP06}, used in our numerical experiments. 
In Section~\ref{sec:num_discre}, we illustrate the spatial discretization of the bidomain model, give a comparison between fully coupled and decoupled
solution strategies, and finally explain our central idea of the partitioned strategy by utilizing the core idea of compile-time local sparse matrices. Section~\ref{sec:num_res} contains detailed numerical results, where the accuracy and efficiency of the proposed strategies are analyzed. These include convergence analysis for the fully coupled, partitioned, and decoupled strategies, simulations of reentrant wave dynamics in a two-dimensional domain, and evaluations on both a realistic three-dimensional ventricular geometry and a bidomain-bath configuration that models tissue-bath interactions.
% we first contrast our approach with a prohibitive dense matrix implementation, and then compare it to a common decoupled solver based on the simulation of
% cardiac reentry with the two aforementioned cell models.

% 
% \end{document}

%%%%%%%%%%%%%%%%%%%%%%%%%%%%%%%%%%%%%%%%%%%%%%%%%%%%%
%%%%%%%%%%%%%%%%%%%%%%%%%%%%%%%%%%%%%%%%%%%%%%%%%%%%%
\section{Governing equations}
%%%%%%%%%%%%%%%%%%%%%%%%%%%%%%%%%%%%%%%%%%%%%%%%%%%%%
%%%%%%%%%%%%%%%%%%%%%%%%%%%%%%%%%%%%%%%%%%%%%%%%%%%%%

\subsubsection*{Bidomain Model}

%numerical simulation of bidomain model equations~\cite{henriquez93:_review,plonsey88:_bioelectric,tung78:_bidomain}
For brevity, we denote $\Omega \subset \Real^d$, $d=2,3$, 
a bounded connected domain with Lipschitz continuous boundary $\partial \Omega$.
The space-time domain denoted by $Q = \Omega \times (0,T]$ and its lateral
boundary is denoted by $\Sigma = \partial \Omega \times (0,T]$.
The dynamics of the intra- and extracellular potentials are described by a set of
reaction-diffusion equations, called the bidomain model, see \cite{Sundnes_book06, ColliFranzone_book14} for 
a detailed description, while the ionic activity is expressed in terms of ordinary
differential equations. 
In the bidomain model, the intra- and extra-cellular
electric potentials $u_i(x, t)$, $u_e (x, t)$ enter into the two 
partial differential equations where transmembrane voltage is denoted by $v=u_i-u_e$, while ionic gating variables
$\wvec(x, t)$ and ion concentrations $\cvec(x, t)$ arise in the ODE system.

The governing equations of the bidomain model are represented by the following equations.
\begin{eqnarray}
\label{eq:elliptic}
 0 &=& \nabla \cdot (\bar{\sigma_i}+\bar{\sigma_e}) \nabla u
            + \nabla \cdot \bar{\sigma_i} \nabla v   \quad \text{in}~ Q \,,\\
\label{eq:parabolic}
 \frac{\partial v}{\partial t} &=& \nabla \cdot \bar{\sigma_i} \nabla v
                    +\nabla \cdot \bar{\sigma_i} \nabla u- I_{ion}(v,\wvec,\cvec) + I_{stim}(x,t) \quad \text{in}~ Q \,,\\
      \label{eq:odes_con}
        \frac{\partial \cvec}{\partial t} &=& S(v,\wvec,\cvec)  \quad \text{in}~Q, \\
      \label{eq:odes}
        \frac{\partial \wvec}{\partial t} &=& G(v,\wvec)  \quad \text{in}~Q,
\end{eqnarray}
where $u \colon Q \to \Real$ is the extracellular potential,
$v \colon Q \to \Real$ is the transmembrane voltage,
$\wvec \colon Q \to \Real^n$ represents the ionic $n$ current variables,
$\cvec \colon Q \to \Real^m$ represents the ionic $m$ concentration variables,
$\bar{\sigma_i} \colon \Omega \to \Real^{d \times d} $ and
$\bar{\sigma_e} \colon \Omega \to \Real^{d \times d} $ are respectively
the intracellular and extracellular conductivity tensors. The term
$I_{stim}$ is the external current density stimulus to generate the excitation wavefronts. 
The $I_{ion}(v,\wvec,\cvec)$ is the current density flowing through the ionic channels
and the functions $G(v,\wvec)$ and $S(v,\wvec,\cvec)$  determines the evolution of the gating variables and concentrations,
which are determined by a specific ionic model, see \cite{CellML}
for more description of these models.
Eq.~(\ref{eq:elliptic}) is an elliptic type equation,
Eq.~(\ref{eq:parabolic}) is a parabolic type equation 
%\commentpk{is this correct in the strict mathematical sense ?}
and Eqs.~(\ref{eq:odes_con}-\ref{eq:odes})
is a set of ordinary differential equations (ODEs) which
need to be solved at each spatial point of the computational domain. 

The conductivity tensors, both intracellular and extracellular, exhibit anisotropy due to the inherent geometry of the myocardium \cite{ColliFranzone_book14}. In our approach, we simplify the computation by modeling the conductivity tensors as diagonal vectors with constant coefficients. However, it is important to note that the proposed method is applicable to cases where the conductivity tensors are anisotropic in nature.
We suppose that a conductive bath is absent, such that the intracellular and extracellular domains are
electrically isolated along the tissue boundaries. Here, we impose no flux boundary conditions.
The initial and boundary conditions for the bidomain model are given by
\begin{eqnarray}\label{eq:bc}
\label{eq:bcpde1}  \eta \cdot (\bar{\sigma_i} \nabla v + \bar{\sigma_i} \nabla u) &=& 0 \quad \textrm{on}~~ \Sigma \,,\\
\label{eq:bcpde2} \eta \cdot \bar{\sigma_e} \nabla u &=& 0 \quad  \textrm{on}~~ \Sigma_{N} \,,\\
u &=& 0 \quad \textrm{on} ~~\Sigma_{D} \,,\\
\label{eq:initcond} u(x,0) = u_0\,, \quad  v(x,0) = v_0\,, \quad  \wvec(x,0) = \wvec_0\quad \textrm{and}~~ \cvec(x,0) &=& \cvec_0 \quad  \textrm{on} ~~\Omega \,,
\end{eqnarray}
where  $\eta$ is the outward unit normal vector to the boundary of the domain, $\Sigma_{N}=\partial\Omega_{N}\times (0,T]$, $\Sigma_{D}=\partial\Omega_{D}\times (0,T]$, with $\partial\Omega_{N}, \partial\Omega_{D} \subset \Omega$ and $\partial\Omega_{N}\cup\partial\Omega_{D}=\partial\Omega$,  
$u_0: \Omega \to \Real^d$ denotes the initial extracellular potential, 
$v_0: \Omega \to \Real^d$ denotes the initial transmembrane potential, and $w_0:\Omega \to \Real^d$ is the initial gating variables at time $t=0$.
%\commentnc{missing Dirichlet}

\subsection*{Ionic models}
The ionic activity across the cell membrane is governed by a system of ordinary differential equations (ODEs), which must be solved at every node of the computational grid. In our study, we employ the following two widely accepted cellular models to describe the ionic current $I_{\text{ion}}$ through membrane channels, as well as the source terms $G(v, \wvec)$ and $S(v, \wvec, \cvec)$ that capture the dynamics of transmembrane voltage and associated state variables.

\subsubsection*{Luo-Rody model (LR91)}
The mammalian ventricular membrane model is described by the Luo-Rudy phase-1 model \cite{Luo_CircRes91}. The ionic current is comprised of the sum of 6 currents
\begin{equation*}
    I_{ion} = I_{Na}+I_{si}+I_{K1}+I_{Kp}+I_{b}
\end{equation*}
 where $I_b$ is the background current, $I_{Kp}$ is the plateau potassium current, $I_{K1}$ is time independent potassium current, $I_{si}$ is the slow inward calcium current, $I_{K}$ is the potassium current which is time-dependent and $I_{Na}$ is the fast sodium current. These currents, which are time-dependent, depend on the gating variables $m,h,j,d,f, X$, which are given by the ODEs of the form
 \begin{equation}
     \frac{dk}{dt} = \alpha_k(v) (1-k)-\beta_k (v) k,
 \end{equation}
for $k=m,h,j,d,f, X$. The rate constants $\alpha_k$ and $\beta_k$ depend on the transmembrane potential $v$. Apart from this, the calcium current $I_{si}$ depends also on the calcium concentration $c$, which is given by
\begin{equation}
    \frac{dc}{dt} = S(v,\wvec,c).
\end{equation}
In conclusion, the LR91 model is comprised of 8 state variables, which include membrane potential, 6 gating variables, and one ionic concentration. In our simulations, we adopt the parameter set as defined in the original study \cite{Luo_CircRes91}, ensuring consistency with the validated physiological model.

\subsubsection*{Ten Tusscher model (TP06)}
The human ventricle action potential for epicardial, endocardial, and
myocardial cells are given by the ten-Tusscher and Panfilov \cite{tenTusscher_AJP06} model 
which is based on data from human ventricular myocyte experiments. This cell model consists of 19 state variables, 
namely, one membrane voltage, 11 gating, and seven concentration variables.
The ionic membrane current $I_{ion}$ has the following 12 components: fast Na$^+$, L-type Ca$^{2+}$, rapid and slow components of the delayed rectifier K$^+$, inward rectifier K$^+$, transient outward K$^+$, plateau K$^+$, Na$^+$-Ca$^{2+}$ exchanger, Na$^+$-K$^+$ pump, sarcolemmal Ca$^{2+}$ pump, and background Na$^+$ and Ca$^{2+}$ currents.
In our simulations, we utilize the epicardial cell parameter set, and the model functions are implemented as specified in the original work \cite{tenTusscher_AJP06}, ensuring alignment with the established electrophysiological framework.

%%%%%%%%%%%%%%%%%%%%%%%%%%%%%%%%%%%%%%%%%%%%%%%%%%%%%
%%%%%%%%%%%%%%%%%%%%%%%%%%%%%%%%%%%%%%%%%%%%%%%%%%%%%
\section{Discretization and Numerical Strategies}
\label{sec:num_discre}
%%%%%%%%%%%%%%%%%%%%%%%%%%%%%%%%%%%%%%%%%%%%%%%%%%%%%
%%%%%%%%%%%%%%%%%%%%%%%%%%%%%%%%%%%%%%%%%%%%%%%%%%%%%
%In this section, we present different approaches for a solver of a strongly
%coupled PDE/ODE system.
%Here we provide a brief description of the spatial discretization by using FEM and 
%different time stepping strategies to solve the bidomain equations.

\subsection{Weak Formulation}
\label{sec:weak_form}
In the following, we denote the inner product and norm in $L^2(\Omega)$ by
$(\cdot,\cdot)$ and $|\cdot|$ respectively, and the
inner product and norm on $H^1(\Omega)$ by
$(\cdot,\cdot)_{H^1}$ and $||\cdot||$. Furthermore,
we set $$H^1_{D}(\Omega) = \set{u\in H^1(\Omega)~:~ u = 0~\text{on}~ \partial\Omega_{D} \subset \Omega}, $$
% $$L^2(\Omega)\slash_\Real = \set{[u]=u-\int_{\Omega} u ~:~ u
% \in L^2(\Omega)},$$
 and $H=L^2(\Omega)$, $V=H^1(\Omega)$,
$U=H^1_{D}(\Omega)$. 
% The purpose of the quotient space reflects the fact that the
% extracellular potential component $u$ of the solution to (\ref{eq:elliptic}) is only defined up to an additive constant. 
It is assumed throughout that for constant $0<m<M$
\begin{equation}
m |\xi|^2 \leq \xi^T \sigma_{i,e}\xi \leq M |\xi|^2 \,, \quad \textrm{for all } \xi \in \Real^2 \,.
\label{eq:uniformelliptic}
\end{equation}

Let $n_c$ and $n_w$ denote the number of concentration and gating variables, respectively. The weak formulation of the bidomain model reads as follows.
A quadruple $(u,v,\cvec,\wvec) \in L^2(0,T;H^1_{D}(\Omega)) \times \big( L^2(0,T;V)\cap C(0,T;H)$ $ \cap L^p(0,T;Q) \big) \times C(0,T;H^{n_c}) \times C(0,T;H^{n_w})$ with $v_t \in
L^2(0,T;V^*) \cap L^{p'}(0,T;Q)$, $\wvec_t \in L^2(0,T;H^{n_w})$ and $\cvec_t \in L^2(0,T;H^{n_c})$ is called
weak solution of (\ref{eq:elliptic}-\ref{eq:odes}), if $\left(v(0),\wvec(0), \cvec(0)\right) =
(v_0,\wvec_0,\cvec_0)$ and if
\begin{eqnarray}
\label{eq:vareqelliptic} \int_{\Omega} \sigma_i \nabla u(t)
\nabla \chi + \int_{\Omega} (\sigma_i+\sigma_e) \nabla v(t) \nabla
\chi &=&0 \,, \\[1.6ex]
\label{eq:vareqparabolic} \seq{v_t(t),\varphi}_{V^*,V} +
\int_{\Omega} \sigma_i \nabla (u(t)+v(t)) \nabla \varphi +
\int_{\Omega} I_{ion}(v(t), \wvec(t), \cvec(t)) \varphi &= &\seq{I_{stim}(t),\varphi}_{V^*,V} \,, \\[1.6ex]
 \label{eq:vareqode_con}
\seq{\cvec_t(t),\phi}_{V^*,V} + \int_{\Omega} S(u(t),\wvec(t),\cvec(t)) \phi &=& 0 \,, \\[1.6ex]
 \label{eq:vareqode}
\seq{\wvec_t(t),\psi}_{V^*,V} + \int_{\Omega} G(u(t),\wvec(t)) \psi &=& 0
\end{eqnarray}
is satisfied for all a.e.~$t \in (0,T)$ and $(\chi, \varphi,\phi, \psi) \in H^1_{D}(\Omega)
\times H^1(\Omega)\times H^{n_c} \times H^{n_w} $. Here, $V^*$ denotes the dual to $V$ with $H$ as pivot space and
$\frac{1}{p} + \frac{1}{p'} =1$. 

The existence and uniqueness of a weak solution to the bidomain model were demonstrated in \cite{Bourgault08} for some
simplified cell models, while existence was guaranteed in \cite{Boulakia_AMRe08} for the Mitchel-Schaeffer cell model.
The well-posedness of the bidomain model coupled with the Luo-Rudy cell model was shown in \cite{Veneroni_NARWA09}, but it is an open problem for other 
physiological cell models.

\subsection{Spatial Discretization}
\label{sec:space_discre}

In this paper, we use the conforming linear FEM for the spatial discretization. 
We choose $V_h \subset H^1(\Omega)$ as a finite-dimensional subspace spanned by
piecewise linear basis functions with respect to the spatial grid.
The  approximate solutions ${u}, {v}$, $\cvec$ and $\wvec$ 
%\commentpk{why boldface for scalar functions}
are expressed in the form
${u}(t) = \sum_{i=0}^{N} {u_{\,i}}(t) \omega_i$,
${v}(t) = \sum_{i=0}^{N} {v_{\,i}}(t) \omega_i$,
$\cvec(t) = \sum_{i=0}^{N} {\cvec_{\,i}}(t) \omega_i$ and
$\wvec(t) = \sum_{i=0}^{N} {\wvec_{\,i}}(t) \omega_i$, respectively, where
$\{\omega_i\}_{i=1}^N$ denote the basis functions.
This semi-discretization in space results in the differential-algebraic system:
\begin{align}
\label{eq:semi-elliptic} 0 &=\Aie \u + \Ai \v \\
\label{eq:semi-parabolic} \M \frac{d \v}{d t} &= -\Ai \v - \Ai \u - \I_{ion}(\v,\w,\c) + \I_{stim} \\
\label{eq:semi-ode_con} \M_{c} \frac{d \c}{d t} &= \mathbf{S}(\v,\w,\c),\\
\label{eq:semi-ode} \M_{w} \frac{d \w}{d t} &= \mathbf{G}(\v,\w),
\end{align}
together with initial conditions for $\v$, $\c$ and $\w$. Here $\Aie = \{\seq{ (\sigma_{i}+\sigma_{e})\nabla \omega_i,\nabla \omega_j}\}_{i,j=1}^N$ and
$\Ai = \{\seq{ \sigma_{i}\nabla \omega_i, \nabla \omega_j}\}_{i,j=1}^N$ are the stiffness matrices,
$\M = \{\seq{\omega_i,\omega_j}\}_{i,j=1}^N$ is the mass matrix, $\M_c$ and $\M_w$ are block diagonal matrices with $n_c$ and $n_w$ blocks of $\M$, respectively, and the vector
$\I_{stim}$ is defined by
$\I_{stim} = \{\langle I_{stim},\omega_j \rangle\}_{j=1}^N$.
The nonlinear reaction terms of the evolution equations are given by
\begin{align*}
 \I_{ion}(\v,\w,\c) &= \I_{ion}\left(\sum_{i=0}^{N} {v_{\,i}} \omega_i,\sum_{i=0}^{N} {\wvec_{\,i}} \omega_i,\sum_{i=0}^{N} {\cvec_{\,i}} \omega_i, \right)\,. \\
 \mathbf{S}(\v,\w,\c) &= \mathbf{S}\left(\sum_{i=0}^{N} {v_{\,i}} \omega_i,\sum_{i=0}^{N} {\wvec_{\,i}} \omega_i,\sum_{i=0}^{N} {\cvec_{\,i}} \omega_i \right)\,, \quad 
 \mathbf{G}(\v,\w) = \mathbf{G}\left(\sum_{i=0}^{N} {v_{\,i}} \omega_i,\sum_{j=0}^{N} {\wvec_{\,j}} \omega_j \right)\,.
\end{align*}

\subsection{Temporal discretization}
\label{sec:temp_discre}
Time discretization schemes for the numerical integration of (\ref{eq:semi-parabolic}-\ref{eq:semi-ode}) can be divided into coupled and decoupled methods.
For a detailed comparative study on the stability and accuracy of various implicit, semi-implicit, and explicit time discretizations of the bidomain model equations with simplified ionic models, we refer to \cite{Ethier_Bourgault_SINUM08}.
A comprehensive survey of methodologies for system decomposition and decoupling is presented in several review articles \cite{Vigmond08:_solvers, Clayton_PBMB11, Pathmanathan_PBMB10}.
The time dependent discretized equations \eqref{eq:semi-elliptic} and (\ref{eq:semi-parabolic}-\ref{eq:semi-ode}) can be expressed as a coupled system
\begin{align}
\label{eq:elliptic-ODE}
\Aie \u &= -\Ai \v, \qquad {\u}(t^0) = {\u_0}\\
\label{eq:genODE}
{\bf \tilde{M}}\frac{d {\bf X}}{d t} &= {\bf F}({\bf X}, \u),
\qquad {\bf X}(t^0) = {\bf X_0}
\end{align}
for the state variable ${\bf X} = (\v,\w,\c)$ and suitable initial solutions $\u_0$ and ${\bf X_0}$. 

\subsubsection{Fully coupled strategy}\label{sec:CoupledStrategy}
The coupled bidomain system can be compactly expressed in block matrix form as
\begin{equation}
 \begin{bmatrix}
   \M& 0&0&0\\ 0 &\M_{c}&0&0 \\0&0&\M_{w}&0\\0&0&0&0
 \end{bmatrix} \begin{bmatrix}
     \dot{\v} \\\dot\c\\\dot\w\\ \dot\u
 \end{bmatrix}
 + \begin{bmatrix}
    \Ai& 0&0& \Ai\\0&0&0&0\\0&0&0&0\\ \Ai&0&0 &\Aie 
 \end{bmatrix}
\begin{bmatrix} \v \\\c\\\w\\ \u
 \end{bmatrix} = \begin{bmatrix} -\I_{ion}(\v, \w, \c)+\I_{stim} \\{\bf S}(\v, \w, \c) \\{\bf G}(\v, \w)\\ 0
 \end{bmatrix}
 \label{eq:coupled_discrete_matrix}
\end{equation}
where $\dot{\v}=\frac{\partial \v}{\partial t}$, with similar notation for $\dot{\c}, \dot\w$ and $\dot{\u}$. For the time discretization of the system \eqref{eq:coupled_discrete_matrix}, 
we consider the second-order diagonally implicit Runge-Kutta method, referred to as the Alexander method \cite{alexander1977diagonally} with time step size $\tau^i$. This scheme is chosen for its favourable stability properties when applied to stiff problems and its ability to achieve better accuracy compared to the backward Euler method.
% implicit Backward
% Euler (BE) method with time step size $\tau^i$, due to its simplicity. 

After time discretization, this scheme requires the solution of a large nonlinear algebraic system at each time step, which is done by Newton's method. 
At each Newton iteration, a large, sparse, and nonsymmetric linear system needs to be solved, for which we employ the preconditioned Generalized Minimal Residual (GMRES) method.
As a preconditioner, a two-level Restricted Additive Schwarz (RAS) method, \cite{toselli2004domain},  with a multiplicative coarse space 
correction based on the partition-of-unity approach by Sarkis,
see \cite{sarkis2003partition}, is used. More specifically, each subdomain
contributes the constant and linear functions multiplied by the partition of unity for the $v$ and $u$ components. For the ion concentrations and the ionic gating variables, no coarse space contribution is necessary. This coarse space has the advantage that it is quite effective and can be built completely algebraically.

\subsubsection{Partitioned strategy}\label{sec:PartitionedStrategy}
The approach to solving the bidomain problem is split into two parts: the system of parabolic PDE and the ODEs \eqref{eq:genODE}, and the elliptic equation \eqref{eq:elliptic-ODE}. At each time step, the solution of \eqref{eq:genODE} determines the transmembrane potential, gating variables, and concentration variables, and \eqref{eq:elliptic-ODE} gives the extracellular potential.
A key advantage of this partitioned approach lies in its flexibility: each subproblem can be addressed independently using different numerical schemes, which enhances computational efficiency. In the following, we present various approaches aimed at solving the bidomain equations more effectively. We begin with a concise explanation of the method used to solve the coupled parabolic PDE and ODE system \eqref{eq:genODE}, followed by an overview of partitioning strategies.

In order to solve \eqref{eq:genODE}, we use the Rosenbrock method, which belongs to the
class of linearly implicit Runge-Kutta methods.  
For brevity, an $s$-stage {\it Rosenbrock} method of order $p$ with the embedding of order $\hat{p}
\neq p$ has the form
\begin{align}
\label{eq:rb_system}
(\frac{1}{\tau^i \gamma}{\bf \tilde{M}}-  {\bf J}){\bf k}_j &= {\bf F}
({t^i+\tau^i \alpha_j, {\bf X}^i + \sum_{l=1}^{j-1} a_{jl}{\bf k}_l}, \u^{i+1}) -
{\bf \tilde{M}}\sum_{l=1}^{j-1} \frac{c_{lj}}{\tau^i} {\bf k}_l \,, \quad j=1,\ldots,s\, , \\
{\bf X}^{i+1} &= {\bf X}^i + \sum_{l=1}^{s} m_l {\bf k}_l \, ,  \quad
\hat{\bf X}^{i+1} = {\bf X}^i + \sum_{l=1}^{s} \hat{m}_l {\bf k}_l \,. \label{eq:second_sol}
\end{align}
The coefficients $\gamma, \alpha_j, a_{jl}, c_{jl}, m_l$, and $\hat{m}_l$ are
chosen in such a way that certain consistency order conditions are
fulfilled to guarantee the convergence orders $p$ and $\hat{p}$.
 For the construction of the
Jacobian matrix ${\bf J}$, we used exact derivatives of the vector $\mathbf
{F}(\mathbf {x})$. Specifically, we used the stiff ODE solver {\sl ROS2} \cite{Lang01}, which has $s=2$ internal stages, is second-order
accurate ($p=2$), satisfies the $L$-stability property and has no order reduction in the PDE case. The second solution $\hat{\bf \X}^{i+1}$ is computed according to eq.~\eqref{eq:second_sol} to determine the future time step, see \cite{Lang01} for complete details.

In the following, we present a concise overview of partitioned strategies used to solve the bidomain model equations. After spatial discretization, the resulting system takes the form of differential-algebraic equations (DAEs). Partitioned schemes are particularly well-suited for such problems, as they effectively decouple the algebraic constraints from the differential components. This separation allows for tailored, computationally efficient treatments of the model's distinct physical phenomena.
In general, partitioned methods are capable of achieving high-order accuracy by independently solving subsystems; however, their application to semi-explicit index-1 DAEs presents challenges. Notably, standard partitioned schemes often suffer from order reduction in this context—both first- and second-order methods tend to degrade to first-order accuracy when applied naively to DAEs \cite{dae1_cdcs_2005}. To address this, we adopt a spectral deferred correction (SDC) approach integrated with partitioned strategies, as proposed in \cite{yoon2021spectral}. This combination not only preserves high-order accuracy but also mitigates the limitations typically encountered with index-1 DAEs. Specifically, we discretize the PDE-ODE subsystem using a Rosenbrock method, and embed this solver within a time-stepping framework that also incorporates the elliptic component of the bidomain model. Our partitioned scheme advances the parabolic-ODE system and the elliptic equation sequentially at each time step, enabling both accuracy and efficiency.

Let us denote ${\bf Z}_{\tau}$ the solver for the problem given by equation \eqref{eq:elliptic-ODE} and let $\Theta_{\tau}$ be the solver for the system of ODEs \eqref{eq:genODE}. Various splitting methods exist for solving the system of differential and algebraic equations (DAEs)\cite{dae1_cdcs_2005}. The one-pass algorithm $\Theta_{\tau}\circ{\bf Z}_{\tau}$, which is first-order accurate for solving the system of equations, proceeds as follows:
\begin{enumerate}
    \item First, solve for $\u^{i+1}$ with fixed $\v^{i}$, 
    $${\bf Z}_{\tau^{i}} : \quad \Aie \u^{i+1} = -\Ai {\v}^{i}.$$
    \item Next, solve the system of ODEs \eqref{eq:genODE} for $\X^{i+1}=(\v^{i+1},\w^{i+1},\c^{i+1})$,
    $$\Theta_{\tau^{i}}: \quad {\bf \tilde{M}}\frac{d {\bf X}}{d t} = {\bf F}({\bf X}, \u^{i+1}).$$
\end{enumerate}
To advance the solution in time, starting with the initial solution $\X_{0}$, the one-pass algorithm can be written as a sequence
\begin{equation*}
\ldots\left( \Theta_{\tau^{i}} \circ {\bf Z}_{\tau^{i}} \right)\ldots\left( \Theta_{\tau^{1}} \circ {\bf Z}_{\tau^{1}} \right) \circ \left( \Theta_{\tau^{0}} \circ {\bf Z}_{\tau^{0}} \right)(\X_{0}).
\end{equation*}
where $\tau^{i}$ denotes the time step size at the $i^{\text{th}}$ step.

To enhance the basic one-pass algorithm, we employ inner iterations at each time step. The algorithm for the resulting partitioned strategy to solve the bidomain equations for each time step is outlined in Algorithm~\ref{alg:partitioned}.
% \begin{algorithm}[]
% \centering
%     \begin{tikzpicture}[node distance=1.5cm]
%     \node (in1) [startstop] {Initial values $(\v^n, \w^n, \c^n, \u^n)$};
%     \node (inner) [startstop, below of = in1] {Initialize for inner iterations $k:=0$, $\X^k:=\X^n$, $\u^k:=\u^n$.};
%     \node (usol) [process,below of=inner] {Solve the elliptic equation \eqref{eq:coupled_discre_elli} for $\u^{k+1}$ with fixed $\v^k$};
%     \node (vsol) [process, below of = usol] {Solve the coupled system \eqref{eq:genODE} for $(\v^{k+1},\w^{k+1},\c^{k+1})$};
%     \node (converge) [startstop1, below of = vsol] {If $\|\X^k-\X^{k+1}\|<tol$ };
%     \node (end) [startstop1, below of = converge] {Update $\X^{n+1}:=\X^{k+1}$, $\u^{n+1}:=\u^{k+1}$};
%     \draw [arrow] (in1) -- (inner);  \draw [arrow] (inner) -- (usol);
%     \draw [arrow] (usol) -- (vsol); \draw [arrow] (vsol) -- (converge);
%     \draw [arrow] (converge) -- node[anchor=east] {yes}(end); 
%     \draw [arrow] (converge.east) -|- [distance=-6mm] node[anchor=south,pos=0.05] {no} node[anchor=east,pos=0.375]{$k:=k+1$} (usol);
%     \draw [arrow] (end.west) -|-[distance=-10mm] node[anchor=west,pos=0.3]{$n:=n+1$}  (in1);
% \end{tikzpicture}
%     \caption{Partitioned Strategy without SDC}
%     \label{alg:partitioned}
% \end{algorithm}

The standard partitioned strategy with inner iterations is expected to achieve second-order accuracy when applied to systems of ODEs. Since the present system (Eqs. (\eqref{eq:elliptic-ODE} and \eqref{eq:genODE}) comprises differential equations and an algebraic equation, partitioning of the system leads to a reduction of the order of accuracy, see for more details \cite{dae1_cdcs_2005} and references therein. In fact, they only achieve first-order accuracy. In order to overcome this issue, we employ a spectral deferred correction \cite{yoon2021spectral} for the solution of the system of ODEs \eqref{eq:genODE} and enhance the overall accuracy. 

In the following, we describe the partitioned strategy incorporating the SDC approach. At each time step, given the solution $\X^n$, we first solve the elliptic constraint equation
 \begin{equation}
 \Aie \u^{n+1} = -\Ai {\v}^{n},
  \label{eq:coupled_discre_elli}
 \end{equation}
 to determine $\u^{n+1}$. With $\u^{n+1}$ known, the variable $\X(t)$ evolves according to the system \eqref{eq:genODE}:
 \begin{equation}
    {\bf \tilde{M}}\frac{d {\bf X}}{d t} = {\bf F}({\bf X}, \u^{n+1})\,,
    \label{eq:discrete_ode}
\end{equation}
over the time interval $t \in [t^n, t^{n+1})$. The system described above is solved using the ROS2 method, as outlined in Eq. \eqref{eq:rb_system}. The solution is then refined through corrections computed via the SDC procedure, as detailed below. The Picard integral formulation of \eqref{eq:discrete_ode} reads:
\begin{equation}
    {\bf \tilde{M}}\X(t) = {\bf \tilde{M}}\X^n + \int_{t^n}^{t}{\bf F}({\bf X}(s), \u^{n+1})ds.
    \label{eq:picard_form}
\end{equation}
Suppose $\tilde{\X}^{n+1}$ is the predicted solution of system \eqref{eq:discrete_ode} at $t^{n+1}$. Then the residual error from the above formulation is
\begin{equation}
  {\bf R}(t^{n+1}) = {\bf \tilde{M}}\X^n + \int_{t^n}^{t^{n+1}}{\bf F}(\tilde{\X}^{n+1}, \u^{n+1})dt - {\bf \tilde{M}}\tilde{\X}^{n+1} \label{eq:residual}
\end{equation}
Let the solution be $\X(t^{n+1})=\tilde{\X}^{n+1}+\boldsymbol{\varepsilon}(t)$. Subtracting \eqref{eq:residual} from \eqref{eq:picard_form}, we get the error equation for the correction $\boldsymbol{\varepsilon}(t^{n+1})$, 
\begin{equation}
    {\bf \tilde{M}}{\bf \varepsilon}(t^{n+1}) = \int_{t^n}^{t^{n+1}}{\bf F}({\tilde{\X}^{n+1}}+{\bf \varepsilon}(t), \u^{n+1})dt - \int_{t^n}^{t^{n+1}}{\bf F}({ \tilde{\X}^{n+1}}, \u^{n+1})dt +{\bf R}(t^{n+1})
    \label{eq:sdc_ode}
\end{equation}
To compute the correction based on Eq. \eqref{eq:sdc_ode}, we employ the first-order implicit integration.  To approximate the integral term in ${\bf R}(t^{n+1})$, we use the trapezoidal rule with two quadrature nodes located at $t^n$ and $t^{n+1}$. The resulting equation is a nonlinear system in ${\bf \varepsilon}^{n+1}$:
\begin{align*}
         {\bf\tilde{M}}{\bf \varepsilon}^{n+1}&= \tau \left({\bf F}(\tilde{\X}^{n+1}+\varepsilon^{n+1}, \u^{n+1})-{\bf F}(\tilde{\X}^{n+1}, \u^{n+1})\right)+{\bf \tilde{M}}\X^n \\
     &\quad+ \frac{\tau}{2}\left({\bf F}(\tilde{\X}^{n+1}, \u^{n+1})+{\bf F}(\tilde{\X}^{n}, \u^{n})\right) - {\bf \tilde{M}}\tilde{\X}^{n+1}.
\end{align*}
Since only an approximate correction is required, we linearize ${\bf F} (\tilde{\X}^{n+1}+\varepsilon^{n+1}, \u^{n+1})$, thereby reducing the nonlinear system to the following linear system of equations.
 \begin{equation}
 \label{eq:sdc}
     \left({\bf \tilde{M}} -\tau {\bf J}\right)\boldsymbol{\varepsilon}^{n+1} = {\bf \tilde{M}}\X^n -{\bf \tilde{M}}\tilde{\X}^{n+1}+\frac{\tau}{2}\left({\bf F}(\X^n, \u^n)+{\bf F}(\tilde{\X}^{n+1},\u^{n+1})\right),
 \end{equation}
where ${\bf J}$ is the Jacobian matrix of partial derivatives of ${\bf F}(\X, \u)$ with respect to $\X$ at $\tilde{\X}^{n+1}$. The corrected solution is then given by $\X^{n+1}=\tilde{\X}^{n+1}+\boldsymbol{\varepsilon}^{n+1}$. 
\newline
 \begin{minipage}[]{0.48\textwidth}
\begin{algorithm}[H]
\centering
    \begin{tikzpicture}[node distance=1.5cm]
    \node (in1) [startstop] {Initial values $(\v^n, \w^n, \c^n, \u^n)$};
    \node (inner) [startstop, below of = in1] {Initialize for inner iterations $k:=0$, $\X^k:=\X^n$, $\u^k:=\u^n$.};
    \node (usol) [process,below of=inner] {Solve the elliptic equation \eqref{eq:coupled_discre_elli} for $\u^{k+1}$ with fixed $\v^k$};
    \node (vsol) [process, below of = usol] {Solve the coupled system \eqref{eq:genODE} for $(\v^{k+1},\w^{k+1},\c^{k+1})$};
    \node (converge) [startstop1, below of = vsol] {If $\|\X^k-\X^{k+1}\|<tol$ };
    \node (end) [startstop1, below of = converge] {Update $\X^{n+1}:=\X^{k+1}$, $\u^{n+1}:=\u^{k+1}$};
    \draw [arrow] (in1) -- (inner);  \draw [arrow] (inner) -- (usol);
    \draw [arrow] (usol) -- (vsol); \draw [arrow] (vsol) -- (converge);
    \draw [arrow] (converge) -- node[anchor=east] {yes}(end); 
    \draw [arrow] (converge.east) -|- [distance=-6mm] node[anchor=south,pos=0.05] {no} node[anchor=east,pos=0.375]{$k:=k+1$} (usol);
    \draw [arrow] (end.west) -|-[distance=-8mm] node[anchor=west,pos=0.3]{$n:=n+1$}  (in1);
\end{tikzpicture}
    \caption{Partitioned Strategy without SDC}
    \label{alg:partitioned}
\end{algorithm}
    \end{minipage}
      \hfill
     \begin{minipage}[]{0.48\textwidth}
\begin{algorithm}[H]
\centering
    \begin{tikzpicture}[node distance=1.5cm]
    \node (in1) [startstop] {Initial values $(\v^n, \w^n, \c^n, \u^n) $};
    \node (inner) [startstop, below of = in1] {Initialize for inner iterations $k:=0$, $\X^k:=\X^n$, $\u^k:=\u^n$.};
    \node (usol) [process,below of=inner] {Solve the elliptic equation \eqref{eq:coupled_discre_elli} for $\u^{k+1}$ with fixed $\v^k$};
    \node (vsol) [process, below of = usol] {Solve the coupled system \eqref{eq:discrete_ode} for $(\v^{k+1},\w^{k+1},\c^{k+1})$};
    \node (correction) [process1, below of = vsol] {Solve the system \eqref{eq:sdc} for $\boldsymbol{\varepsilon}^{k+1}$, and update $\X^{k+1}:=\X^{k+1}+\boldsymbol{\varepsilon}^{k+1}$};
    \node (converge) [startstop1, below of = correction] {If $\|\X^k-\X^{k+1}\|<tol$ };
    \node (end) [startstop1, below of = converge] {Update $\X^{n+1}:=\X^{k+1}$, $\u^{n+1}:=\u^{k+1}$};
    % \node (update) [process, right of=converge, xshift=5cm] {Update $k=k+1$};
    \draw [arrow] (in1) -- (inner);  \draw [arrow] (inner) -- (usol);
    \draw [arrow] (usol) -- (vsol); \draw [arrow] (vsol) -- (correction);
    \draw [arrow] (correction) -- (converge);
    \draw [arrow] (converge) -- node[anchor=east] {yes}(end); 
    % \draw [arrow] (converge) -- node[anchor=south] {no}(update);
    % \draw [arrow] (update) |- (usol);
    \draw [arrow] (converge.east) -|- [distance=-7mm] node[anchor=south,pos=0.05] {no} node[anchor=east,pos=0.33]{$k:=k+1$} (usol);
    \draw [arrow] (end.west) -|-[distance=-8mm] node[anchor=west,pos=0.3]{$n:=n+1$}  (in1);
\end{tikzpicture}
    \caption{Partitioned Strategy with SDC}
    \label{alg:partitioned_sdc}
\end{algorithm}
 \end{minipage}
\newline

Solving the coupled system Eq.~\eqref{eq:discrete_ode} using the implicit time integrations comprising a parabolic PDE and a set of gating and concentration ODEs poses significant memory challenges, particularly given the large-scale and detailed nature of the physiological cell models under consideration. To address this issue, we adopt a memory-efficient compile-time sparse matrix (CTSM) technique, as introduced in \cite{NP_AMC19}, within our FEM implementation. This approach exploits the inherent weak coupling among cellular state variables by representing local block matrices, generated at each finite element, as compile-time sparse matrices. This design significantly reduces memory overhead for matrix storage, accelerates the assembly of the global compressed row storage (CRS) matrix, and enhances the performance of the linear solvers. In addition, we employ higher-order time-stepping schemes with adaptive step size control, further improving computational efficiency. Together, these strategies substantially reduce both computational cost and memory usage, enabling the coupled parabolic PDE and ODE system to be solved efficiently. For more implementation details and performance benchmarks of the monodomain model, we refer the reader to \cite{NP_AMC19}. 
%This framework also allows for accurate and efficient simulation of complex pathological phenomena such as cardiac reentry.
In this work, we adopt the CTSM-based implementation to ensure scalability and performance. At each time step, the resulting systems of linear equations are solved using robust iterative methods. Specifically, for the internal stages of the ROS2 method \eqref{eq:rb_system}, we employ the BiCGSTAB algorithm \cite{vorst-BiCGSTAB} with block Jacobi preconditioning. The same BiCGSTAB solver with block Jacobi preconditioning is used to solve the linearized deferred correction system \eqref{eq:sdc}, ensuring consistency and efficiency across the stages of the time-stepping scheme. After full spatial discretization, the final algebraic system \eqref{eq:coupled_discre_elli} is solved using the BiCGSTAB method in combination with an algebraic multigrid (AMG) preconditioner, which is developed using a greedy heuristic
algorithm based on the strength of connection criterion \cite{MarkusBlatt_PhD}. 
The corresponding partitioned strategy with SDC is outlined in Algorithm \ref{alg:partitioned_sdc}.

\subsubsection{Decoupled strategy}\label{sec:DecoupledStrategy}
In this subsection, we review established decoupled strategies for solving the bidomain model equations \eqref{eq:semi-elliptic}–\eqref{eq:semi-ode}. These methods are widely adopted in the cardiac electrophysiology research community, and we use them as benchmarks to compare the performance of our proposed approach. Specifically, we consider a commonly used operator splitting technique that decouples the diffusion operator representing electrical conduction in the tissue from the reaction operator, which models ionic currents, as well as the dynamics of gating and concentration variables \cite{Austin_TBE06, Plank_AMG_TBE07, Vigmond08:_solvers}.
For this purpose, we employ the Strang splitting scheme, in which the diffusion term is discretized using the Crank–Nicolson method, and the reaction terms are integrated using a second-order explicit Runge–Kutta method. This approach, referred to as CN-RK2, has been shown to achieve a good balance between stability and accuracy \cite{roy2020analysis}. Additionally, we compare our results with various other semi-implicit and explicit methods for solving the bidomain equations, as described in \cite{Ethier_Bourgault_SINUM08}.

The semi-implicit time stepping schemes implemented include the first-order Implicit-Explicit (IMEX) method, the Implicit-Explicit first-order Gear (IMEX Gear) method, the Crank Nicolson-Forward Euler (CN-FE) method, the second-order Crank Nicolson-Adams Bashforth (CN-AB) method, and the second-order semi-implicit backward differentiation (SBDF) method. We also compare these schemes with a popular decoupled strategy for solving the bidomain equations \cite{colli2013comparison}, which we refer to as CN-IMEX for consistency in naming. 
% The numerical formulations for these methods are provided in Appendix~\ref{semi_implicit}.
Such operator splitting methods help reduce strong couplings between state variables and enhance computational efficiency by allowing the use of different numerical schemes for the diffusion and reaction terms. However, their primary drawback lies in the loss of numerical accuracy due to the temporal decoupling of state variables, which are not updated simultaneously. To restore this accuracy, smaller time steps are often required, leading to increased computational costs. We illustrate this trade-off with several examples in the numerical results section.

The Strang splitting method with the Crank-Nicolson and Runga Kutta method (CN-RK2)\cite{roy2020analysis} for the bidomain equations is as follows:
\begin{enumerate}
    \item  Initially, a half-time step of the second-order Runge-Kutta method is used on the reaction part of the parabolic equation, as well as for the gating and concentration variables.
     \begin{align}
     \frac{\M}{\Delta t/2}( \v^{*}-\v^{n}) &= - \I_{ion} \left(\v^n+\frac{\Delta t}{4}\I_{ion}(\v^n, \w^{n}, \c^{n}), \w^{n}+\frac{\Delta t}{4}\mathbf{G}^n, \c^{n}+\frac{\Delta t}{4}\mathbf{S}^n\right)\nonumber \label{eq:cnrk_odes}\\ 
     \frac{\w^{*}-\w^n}{\Delta t/2} &=\mathbf{G}\left(\v^n+\frac{\Delta t}{4}\I_{ion}(\v^n, \w^{n}, \c^{n}),\w^{n}+\frac{\Delta t}{4}\mathbf{G}(\v^n,\w^{n})\right) \\
     \frac{(\c^{*}-\c^n)}{\Delta t/2} &= \mathbf{S}\left(\v^n+\frac{\Delta t}{4}\I_{ion}(\v^n, \w^{n}, \c^{n}),\w^{n}+\frac{\Delta t}{4}\mathbf{G}^n+\c^{n}+\frac{\Delta t}{4}\mathbf{S}^n\right) \nonumber    
 \end{align}
 where $\mathbf{G}^{n}$ denotes for $\mathbf{G}(\v^n,\w^{n})$ and $\mathbf{S}^{n}$ denotes for $\mathbf{S}(\v^n,\w^{n},\c^n)$.
 \item  The diffusion component of the parabolic equation is then discretized using a step of the Crank–Nicolson method. Solve for the intermediate value $\v^{**}$ from the following equation.
 \begin{equation}
     \M\frac{\v^{**}-\v^{*}}{\Delta t} =- \frac{1}{2}\Ai(\v^{**}+\v^{*})-\Ai \u^n +\I_{stim}
     \label{eq:decoupled_discre_parabolic_cn}
 \end{equation}
 \item Subsequently, a half-time step method of Runge-Kutta is applied for the reaction part of the parabolic equation, and the same for that of the gating and concentration variables.
 \begin{align}
     \frac{\M}{\Delta t/2}( \v^{n+1}-\v^{**}) &= -\I_{ion} \left(\v^{**}+\frac{\Delta t}{4}\I_{ion}^{**}, \w^{**}+\frac{\Delta t}{4}\mathbf{G}^{**}, \c^{**}+\frac{\Delta t}{4}\mathbf{S}^{**}\right)\nonumber\\ 
     \frac{\w^{n+1}-\w^{**}}{\Delta t/2} &=\mathbf{G}\left(\v^{**}+\frac{\Delta t}{4}\I_{ion}^{**},\w^{**}+\frac{\Delta t}{4}\mathbf{G}^{**}\right) \\
     \frac{\c^{n+ 1}-\c^{**}}{\Delta t/2} &= \mathbf{S}\left(\v^{**}+\frac{\Delta t}{4}\I_{ion}^{**},\w^{**}+\frac{\Delta t}{4}\mathbf{G}^{**}+\c^{**}+\frac{\Delta t}{4}\mathbf{S}^{**}\right) \nonumber
 \end{align}
where $\I_{ion}^{**}$ denotes for $\I_{ion}(\v^{**}, \w^{**}, \c^{**})$, $\mathbf{G}^{**}$ denotes for $\mathbf{G}(\v^{**},\w^{**})$, and $\mathbf{S}^{**}$ denotes for $\mathbf{S}(\v^{**},\w^{**},\c^{**})$.
\item  Finally, we solve the elliptic equation for the extracellular potential:
 \begin{equation}
 \Aie \u^{n+1} = -\Ai \v^{n+1}.
 \label{eq:decoupled_discre_elli_2}
 \end{equation}
\end{enumerate}

For solving the systems of linear equations in Eqs.~\eqref{eq:cnrk_odes}, and \eqref{eq:decoupled_discre_parabolic_cn}, we used a BiCGSTAB method with block Jacobi preconditioning. 
Equation~\eqref{eq:decoupled_discre_elli_2} is solved using the BiCGSTAB method preconditioned with the AMG preconditioner, consistent with the approach used for Eq.~\eqref{eq:coupled_discre_elli}.
In parallel numerical simulations, the main essence is to distribute the
computational load evenly across all processors, which is mainly achieved by
domain decomposition algorithms.  
The internal {YASP} grid in DUNE \cite{Dune_Computing08} is used for parallel grid 
constructions in our numerical simulations.

\section{Numerics}\label{sec:num_res}
This section presents numerical results obtained using fully coupled, partitioned, and decoupled simulation strategies. The implementation of these numerical strategies is based on the open-source finite element software framework DUNE \cite{Dune_Computing08} and its module Dune-PDELab \cite{Dune_pdelab}. The computations were performed on the Padmanabha HPC cluster in IISER-TVM using up to two Intel Xeon Gold 6132 computing nodes, each consisting of 28 cores with 128 GB of main memory. The parameters for the physiological cell models TP06 \cite{tenTusscher_AJP06} have been set as specified in their original formulations. Some of the parameters for the LR91 \cite{Luo_CircRes91} cell model were changed from their original formulation to generate spiral wavefronts. These modified values are $G_{si}=0.02$, $G_{Na}=4.0$, and $G_{K}=0.423$. The intracellular and extracellular conductivity coefficients are set to $1.5\times10^{-3}$ in both the transversal and longitudinal directions in all the cases. In the following, we present convergence tests with respect to various spatial grid sizes and time step sizes for all implemented strategies. For all convergence tests involving the LR91 and TP06 models, the computational domain is taken as $[0,1] \times [0,1]$ cm$^2$, and the spatial grid is fixed at $100 \times 100$ quadrilateral elements unless otherwise specified. Based on the outcomes of these convergence tests, we provide a comprehensive analysis of the ventricular arrhythmia simulations.

\subsection{Convergence analysis of the fully coupled strategy}
\begin{figure}
\centering
\begin{subfigure}[t]{0.495\textwidth}
      \includegraphics[width=\linewidth]{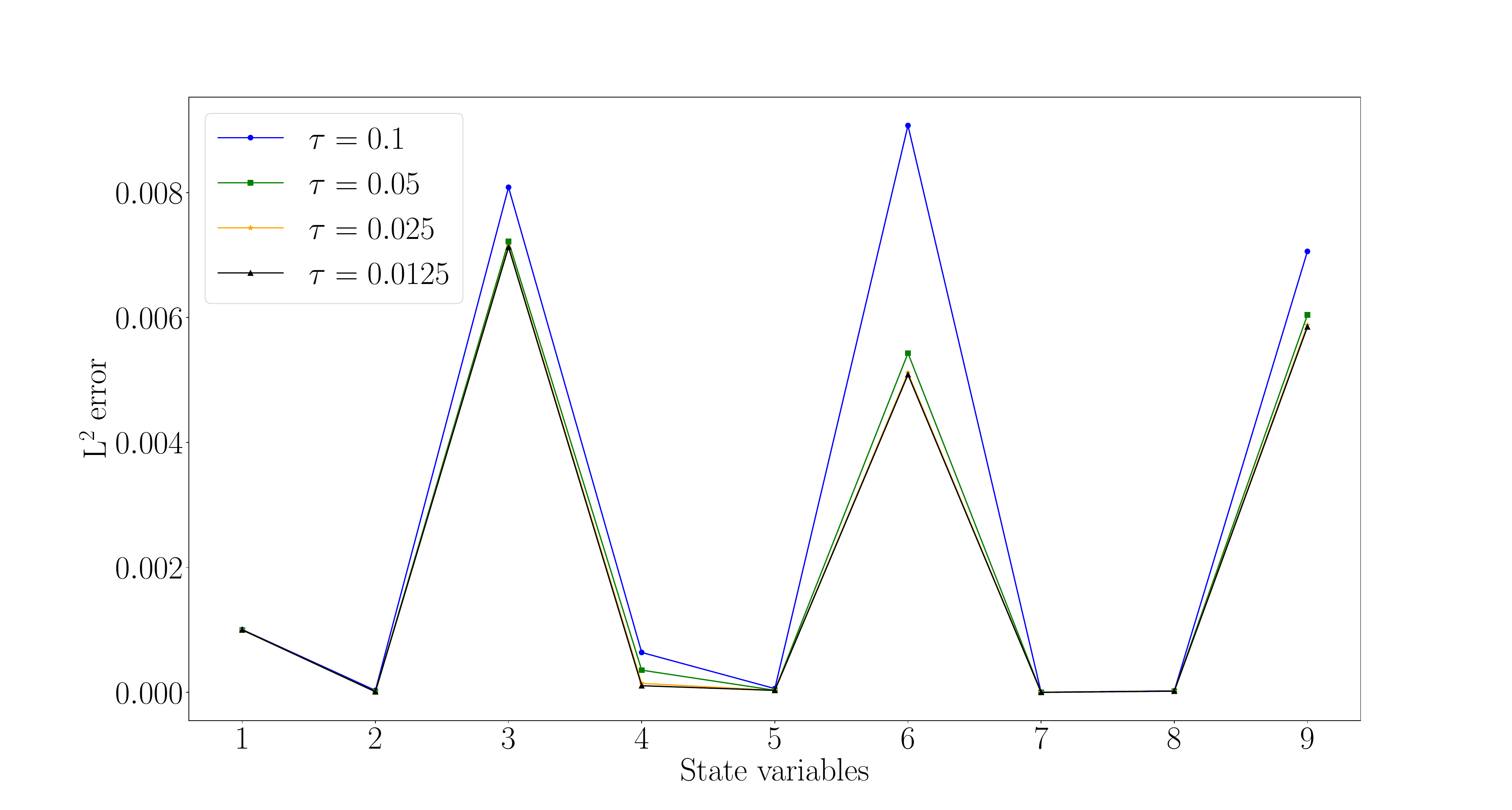}
      \caption{Relative L$^2$ error of the state variables at \mbox{$t=0.5\text{ ms}$} for different time steps.}
    \end{subfigure}
\hfill
\begin{subfigure}[t]{0.495\textwidth}
     \includegraphics[width=\linewidth]{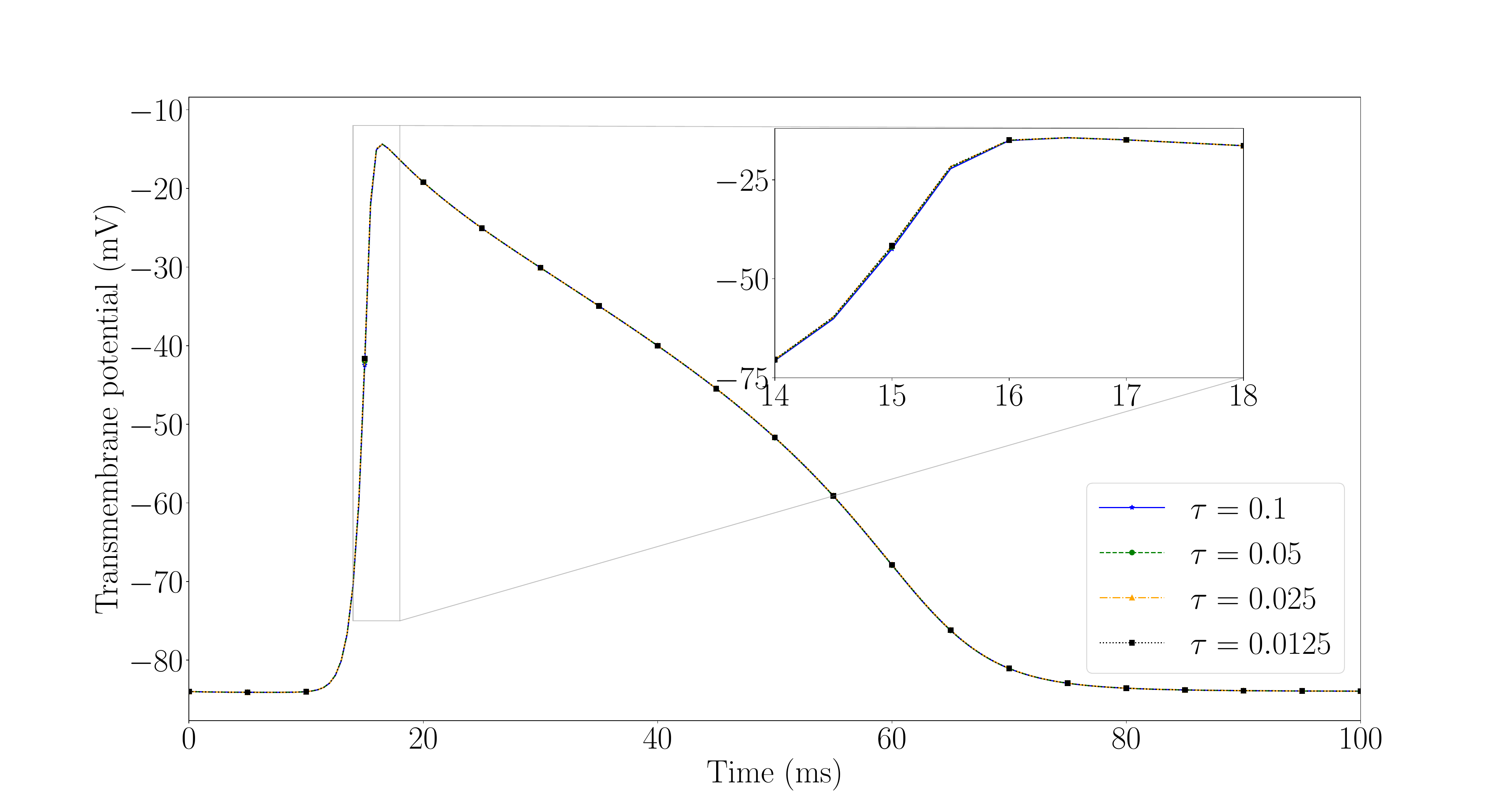}
      \caption{ Transmembrane potential $v$ over a period of time $100\text{ ms}$ at point $(0.5,0.5)$ on the computational grid.}
\end{subfigure} 
\caption{Convergence in time for the fully coupled strategy using the LR91 model.}
\label{fig:lr_time}
\end{figure}

\begin{figure}
\centering
\begin{subfigure}[t]{0.495\textwidth}
      \includegraphics[width=\linewidth]{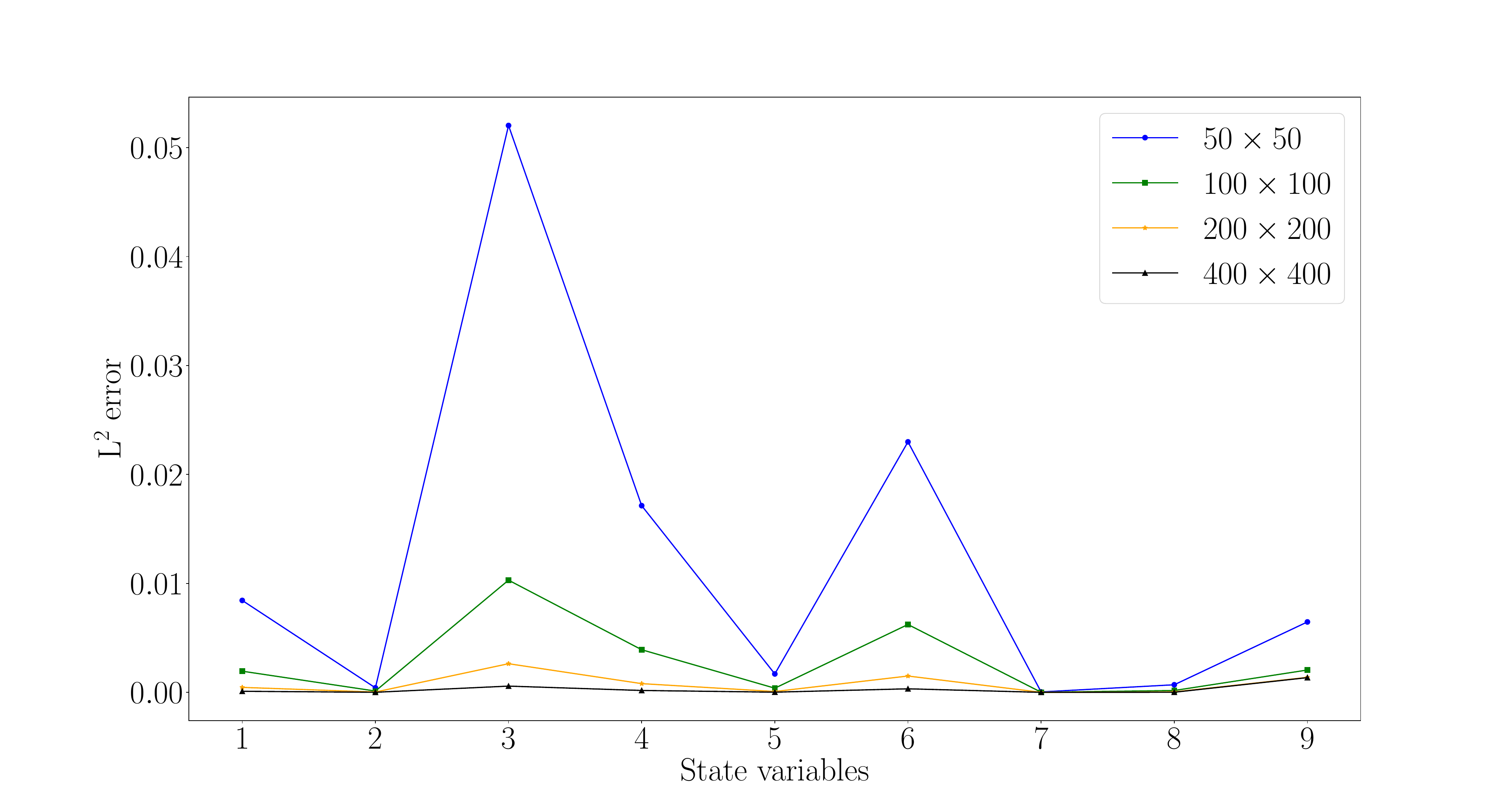}
      \caption{Relative L$^2$ error of the state variables at \mbox{$t=0.5$}\text{ ms} for different grid resolutions.}
    \end{subfigure}
\hfill
\begin{subfigure}[t]{0.495\textwidth}
     \includegraphics[width=\linewidth]{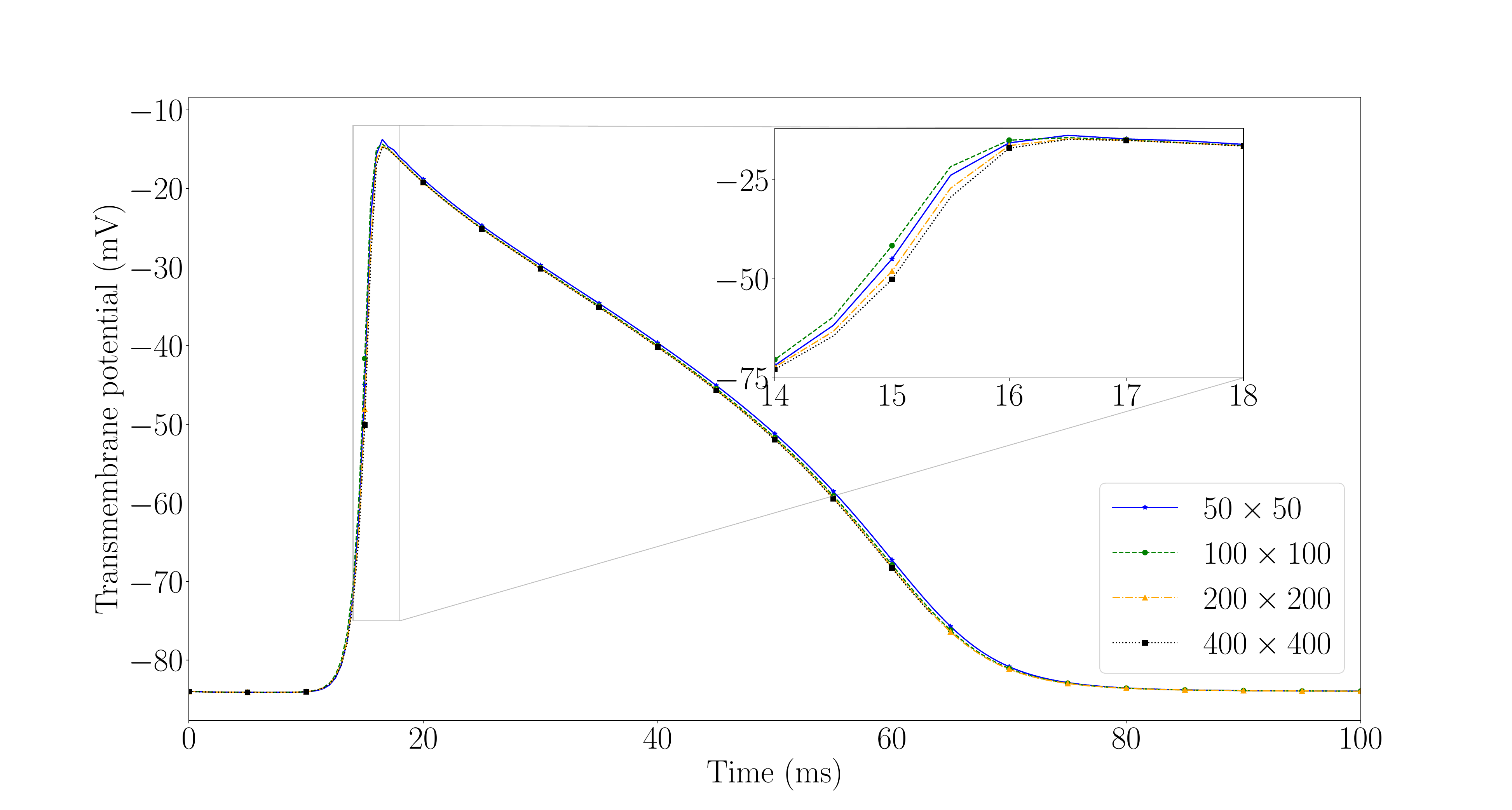}
      \caption{Transmembrane potential $v$ over a period of time $100\text{ ms}$ at the midpoint of the unit square domain.}
\end{subfigure}
\caption{Convergence in space for the fully coupled strategy using the LR91 model.}
\label{fig:lr_space}
\end{figure}

 Using the fully coupled strategy, the simulations are performed by varying the time-step size and spatial resolution to test the convergence of solutions in space and time. Figure~\ref{fig:lr_time} presents the results of the temporal convergence study using the LR91 model. Simulations are performed with the time step size successively halved from $\tau = 0.1$ to $\tau = 0.0125$. The relative $L^2$ error of the nine state variables at $t = 0.5\,\text{ms}$ is shown in Figure~\ref{fig:lr_time}(a), while the transmembrane potential $v$ over a duration of $100\,\text{ms}$ at the midpoint of the domain is plotted in Figure~\ref{fig:lr_time}(b). Since the exact solution of the bidomain model is not available, a reference solution is computed on a finer resolution using a grid of $800 \times 800$ cells and a time step of $\tau = 0.001$. To assess spatial convergence, simulations are conducted with a fixed time step size of $\tau = 0.0125$, which is sufficiently small to ensure temporal accuracy. The grid resolution is refined from $50 \times 50$ to $400 \times 400$ cells. The spatial convergence results are depicted in Figure~\ref{fig:lr_space}. It is observed that the fully coupled strategy achieves convergence for $\tau = 0.05$ and a grid resolution of $100 \times 100$ cells on the unit square domain.


\subsection{Convergence analysis of the partitioned strategies}
In this section, a comparison of partitioned strategies without and with SDC for different time steps is presented. 
%The domain $\Omega=[0,1]^2$ with $100\times 100$ cells is chosen for the computation to analyse the convergence of the partitioned strategies.
The reference solution is computed using a fully coupled strategy with a fixed time step size of $\tau = 0.001$, serving as a benchmark for evaluating the accuracy of the partitioned strategy. The time step size $\tau$ is varied from $0.1$ to $0.003125$ to examine the convergence trends. The relative L$^2$ error for the solution of transmembrane potential $v$, for each time step size, along with the observed numerical order of accuracy for the partitioned strategies without and with SDC, is reported in Table \ref{tab:partitioned_order}. The comparison is made at $0.5\,\text{ms}$, and a tolerance of $10^{-3}$ is prescribed for the convergence of the inner iterations. Consistently comparable relative errors were observed across all eight remaining state variables, reinforcing the robustness of the results. As theoretically expected, the standard partitioned strategy demonstrates only first-order accuracy, see Table \ref{tab:partitioned_order}. In contrast, the SDC-enhanced partitioned strategy clearly achieves second-order accuracy, validating its effectiveness in improving the convergence rate.
\begin{table}
    \centering
    \begin{tabular}{|l|ll|ll|}
   \hline
       Time Step & \multicolumn{2}{|c|}{Partitioned without SDC}& \multicolumn{2}{|c|}{Partitioned with SDC} \\
       $\tau$& $\text{L}^2$ error ($v$) &Order&$\text{L}^2$ error ($v$) &Order\\
       \hline
        $0.1$ &0.00188564 &1.41431339 &0.00090457&1.531963964\\
       $0.05$ &0.00070747 &0.7912878925 &0.000312806&1.875682613\\
       $0.025$ &0.000408796 & 0.8396921041&0.000085239 &1.919346218\\
       $0.0125$ &0.00022842&0.9117779554 &0.000022535&2.113995058\\
       $0.00625$ &0.000121412&0.9527042513&0.00000520573 &1.992054262\\
       $0.003125$ &0.0000627291 &0.9765556426&0.00000130862&1.991413224\\
       \hline
    \end{tabular}
    \caption{Numerical order of the partitioned strategy with and without SDC at $t = 0.5\text{ ms}$, with $tol = 10^{-6}$.}
    \label{tab:partitioned_order}
\end{table}

\begin{figure}
\centering
  \includegraphics[width=0.495\textwidth]{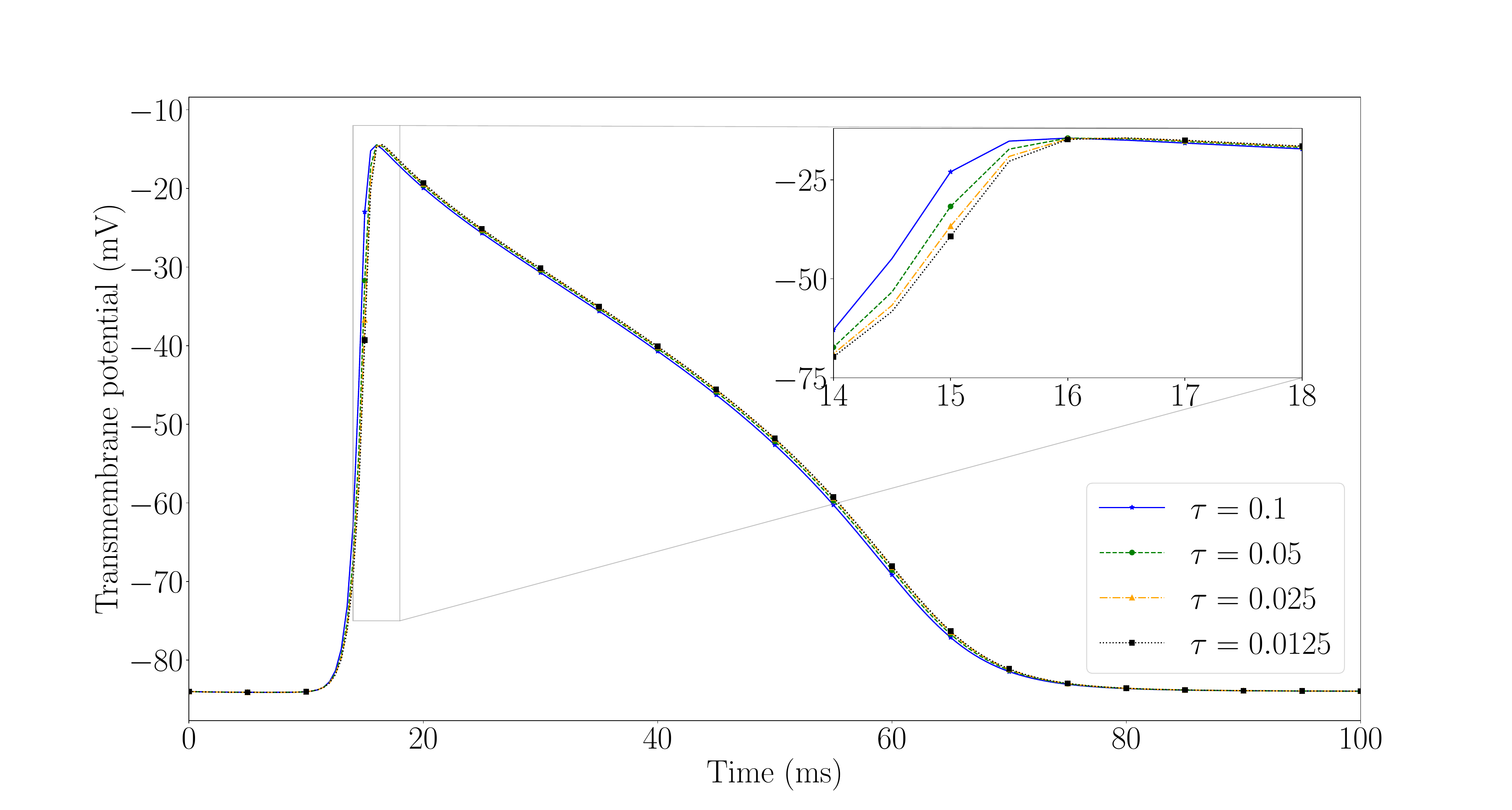}
  % \hfill
  \includegraphics[width=0.495\textwidth]{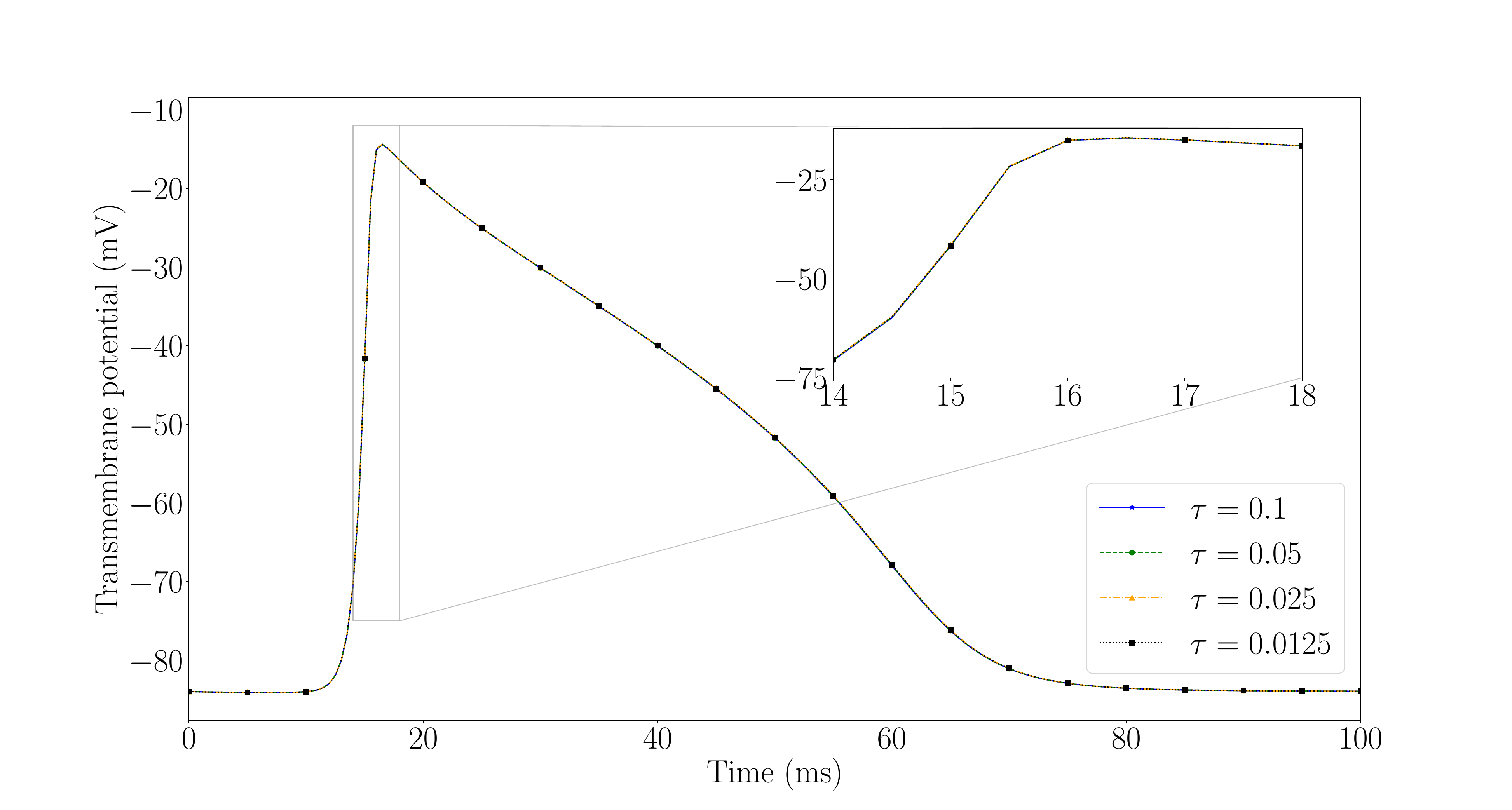}
% \end{subfigure}
\caption{Convergence for Partitioned strategy without SDC (on the left) and with SDC (on the right) using LR91 model for different time steps.}
\label{fig:lr_partitioned}
\end{figure}

It is observed that the partitioned strategy without SDC exhibits convergence towards the reference solution as the time step is refined, consistent with its first-order temporal accuracy. In contrast, the partitioned strategy with SDC is accurate even with coarser time steps. Furthermore, the number of inner iterations per time step is notably reduced under the same convergence criterion and time step size. On average, for $\tau=0.05$, the partitioned strategy without SDC requires 12 inner iterations, while only 5 inner iterations are required for the partitioned strategy with SDC. The solution for transmembrane potential over a duration of $100\text{ ms}$ at the location $(0.5,0.5)$ in the computational domain is depicted in Fig.~\ref{fig:lr_partitioned}.

\begin{figure}
\centering
\begin{subfigure}[t]{0.495\textwidth}
      \includegraphics[width=\linewidth]{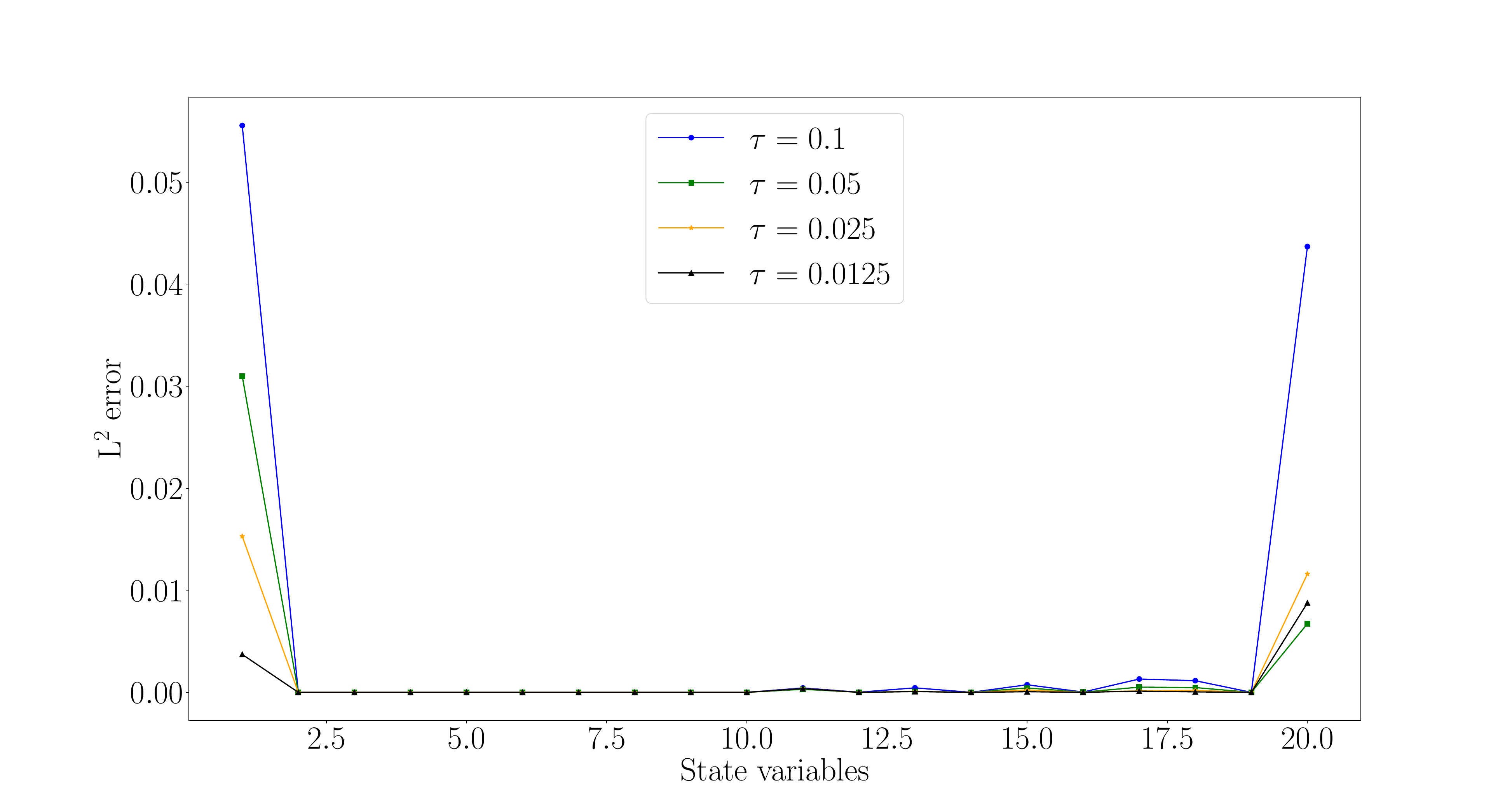}
      \caption{ The absolute L$^2$ error for the state variables at $60\text{ ms}$.}
    \end{subfigure}
\hfill
\begin{subfigure}[t]{0.495\textwidth}
     \includegraphics[width=\linewidth]{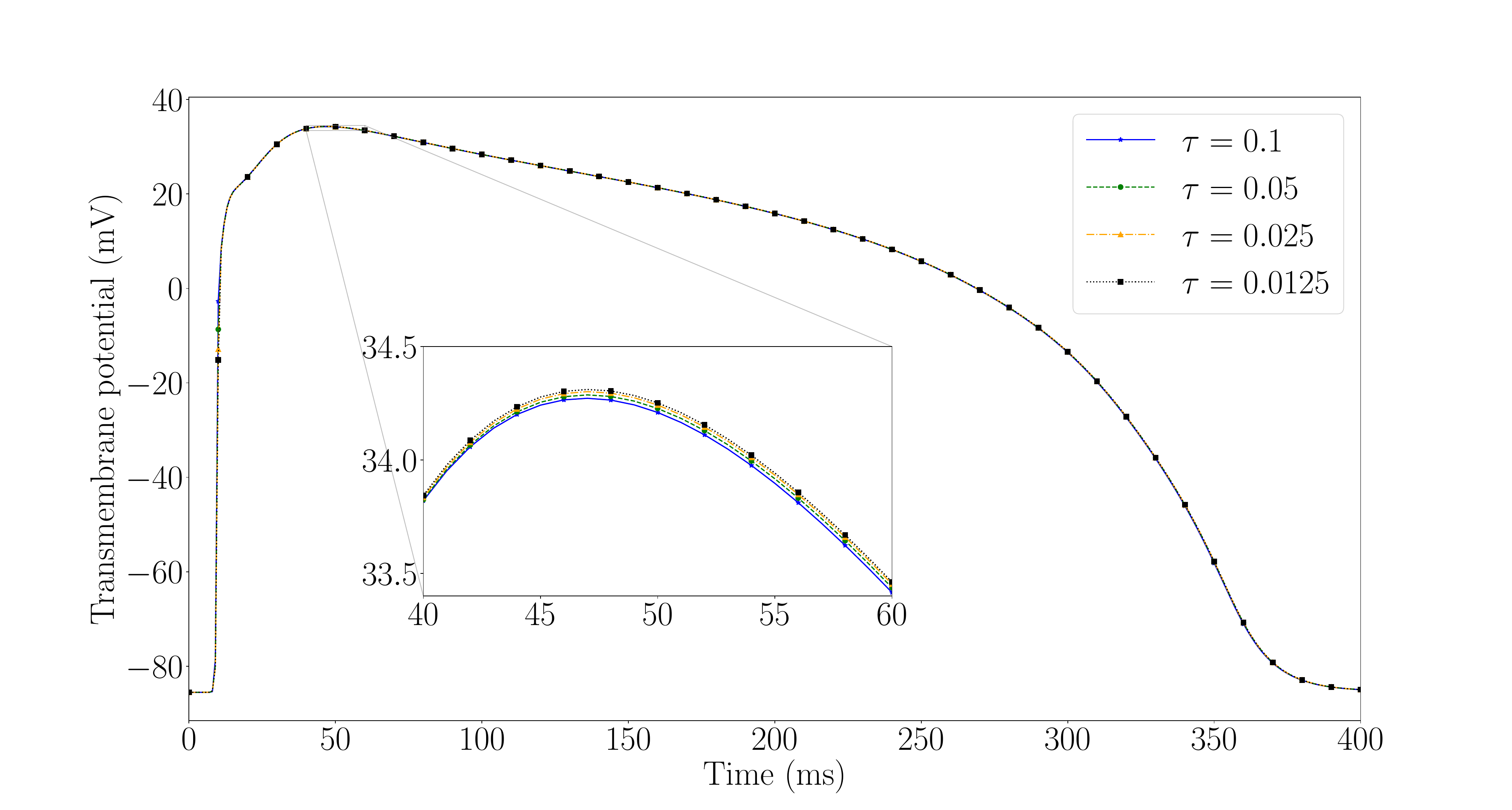}
      \caption{Transmembrane potential $v$ over $400\text{ms}$ at the midpoint of the domain.}
\end{subfigure}
\caption{Convergence for Partitioned strategy without SDC using TP06 model for different time steps.}
\label{fig:tp_partitioned}
\end{figure}

\begin{figure}
\centering
\begin{subfigure}[t]{0.495\textwidth}
     \includegraphics[width=\linewidth]{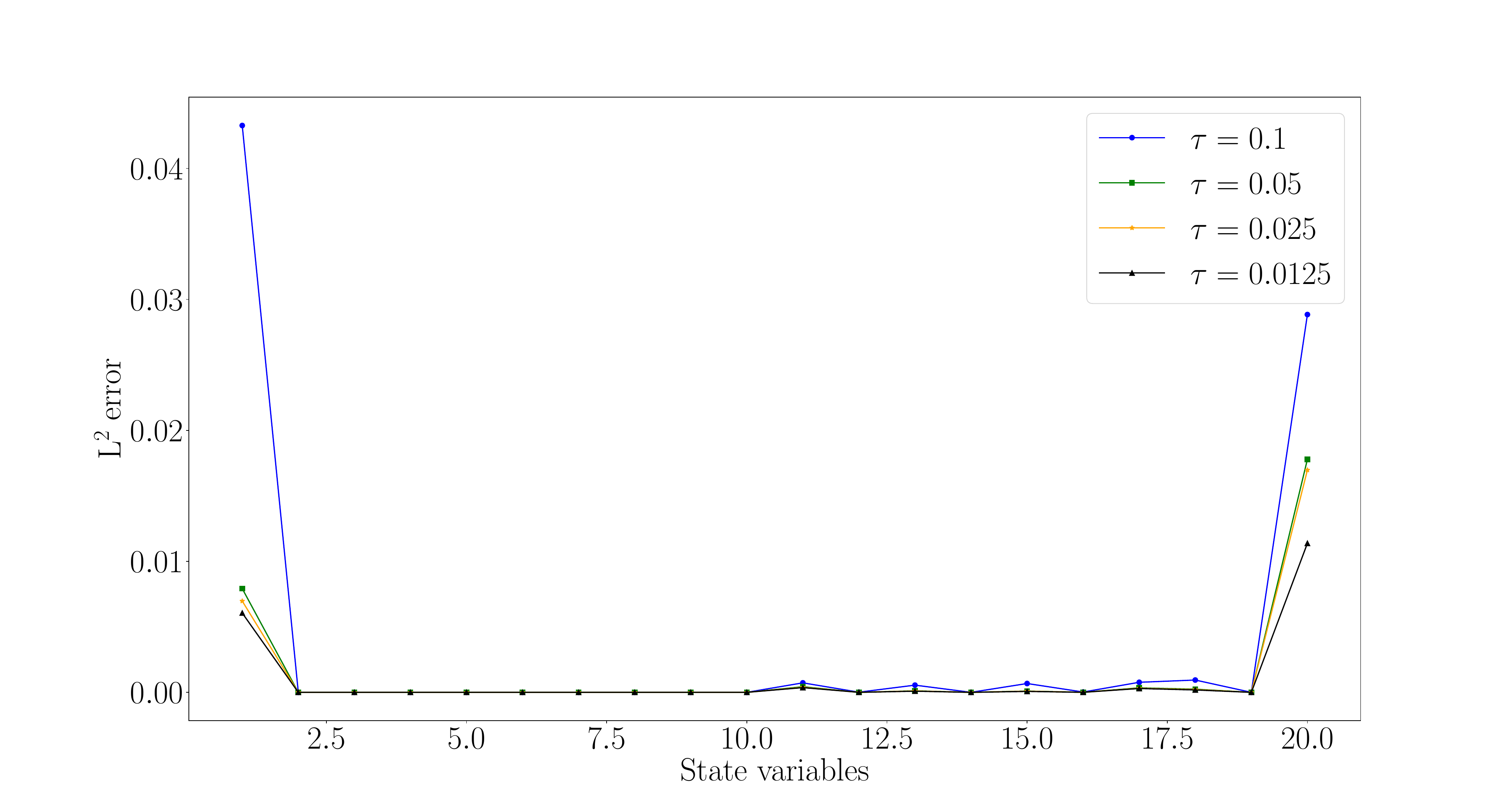}
      \caption{The absolute L$^2$ error for the state variables at $60\text{ ms}$.}
    \end{subfigure}
\hfill
\begin{subfigure}[t]{0.495\textwidth}
     \includegraphics[width=\linewidth]{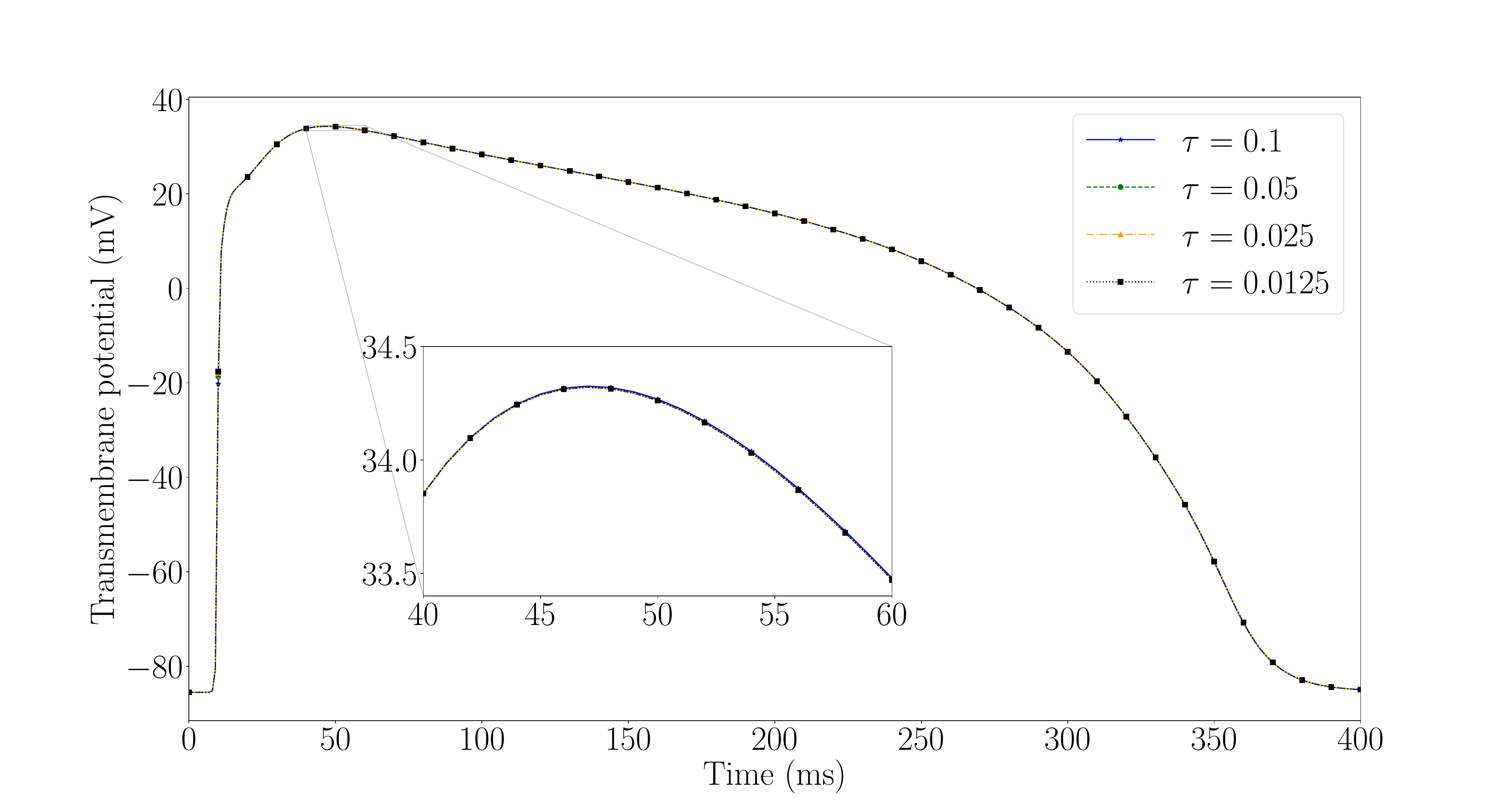}
      \caption{Transmembrane potential $v$ over $400\text{ms}$ at the midpoint of the domain.}
\end{subfigure}
\caption{Convergence for Partitioned strategy with SDC using TP06 model for different time steps.}
\label{fig:tp_partitioned_sdc}
\end{figure}

We performed a numerical convergence analysis of the TP06 model with respect to various spatial and temporal resolutions and observed results consistent with those obtained for the LR91 model.  Comparable results are also observed with the TP06 model. Fig.~\ref{fig:tp_partitioned} and Fig.~\ref{fig:tp_partitioned_sdc} show the absolute L$^2$ error at $60\text{ ms}$ for the 20 state variables (left), along with the evolution of the transmembrane potential $v$ over a time period of $400\text{ ms}$ in the center of the domain (right). Fig.~\ref{fig:tp_partitioned} corresponds to the partitioned strategy without SDC, while Fig.~\ref{fig:tp_partitioned_sdc} presents results with SDC, for time step sizes $\tau = 0.1$, $0.05$, $0.025$, and $0.0125$. In this human ventricular ionic model, as theoretically expected, we note that the standard partitioned strategy demonstrates only first-order temporal accuracy, whereas the SDC-enhanced partitioned strategy successfully achieves second-order accuracy, further validating its effectiveness in improving the convergence rate.


\subsection{Convergence analysis of the decoupled strategy}
\begin{figure}
\centering
\begin{subfigure}[t]{0.495\textwidth}
     \includegraphics[width=\linewidth]{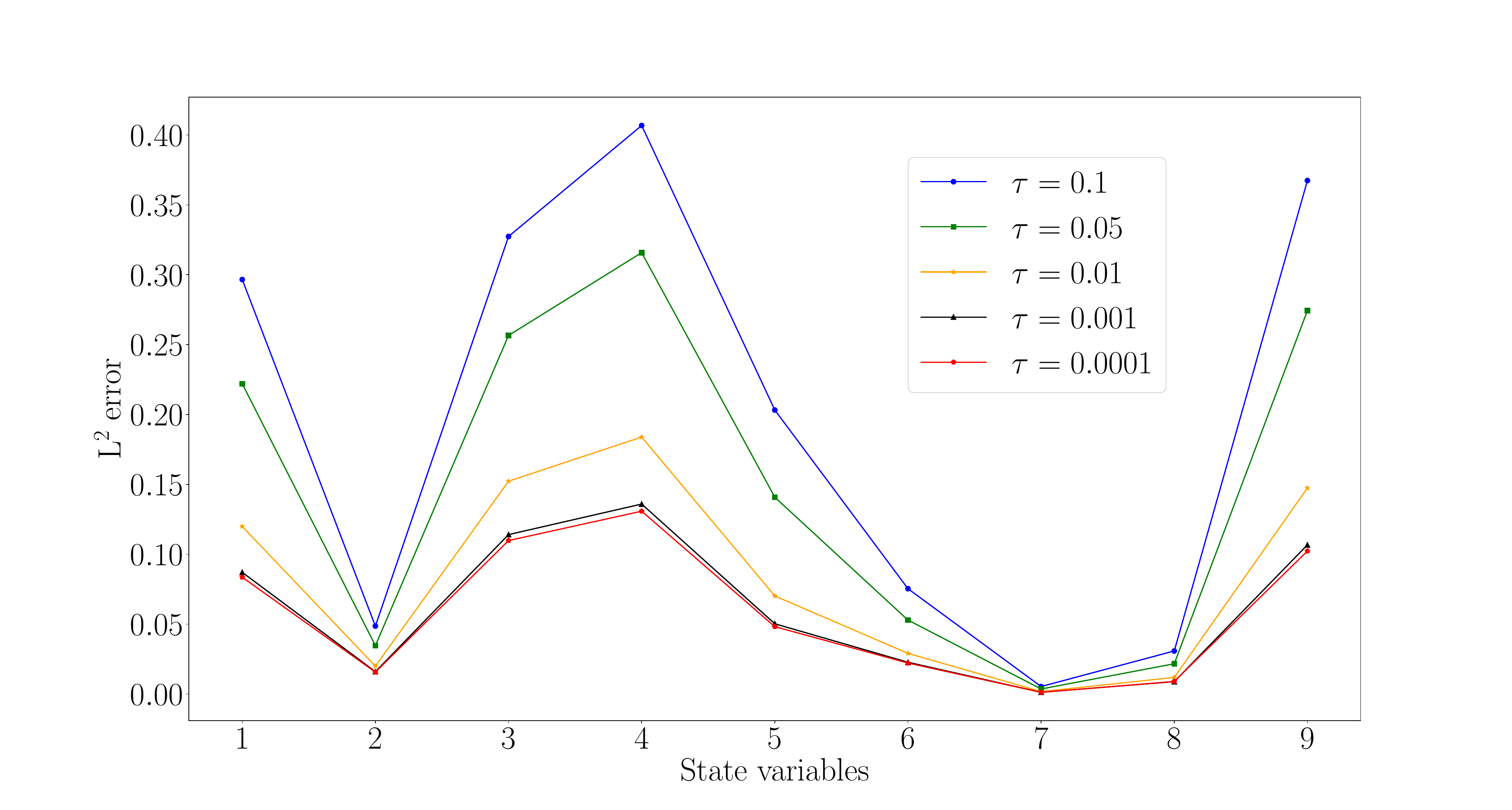}
      \caption{Relative L$^2$ for the state variables at $t=20\text{ ms}$.}
    \end{subfigure}
\hfill
\begin{subfigure}[t]{0.495\textwidth}
    \includegraphics[width=\linewidth]{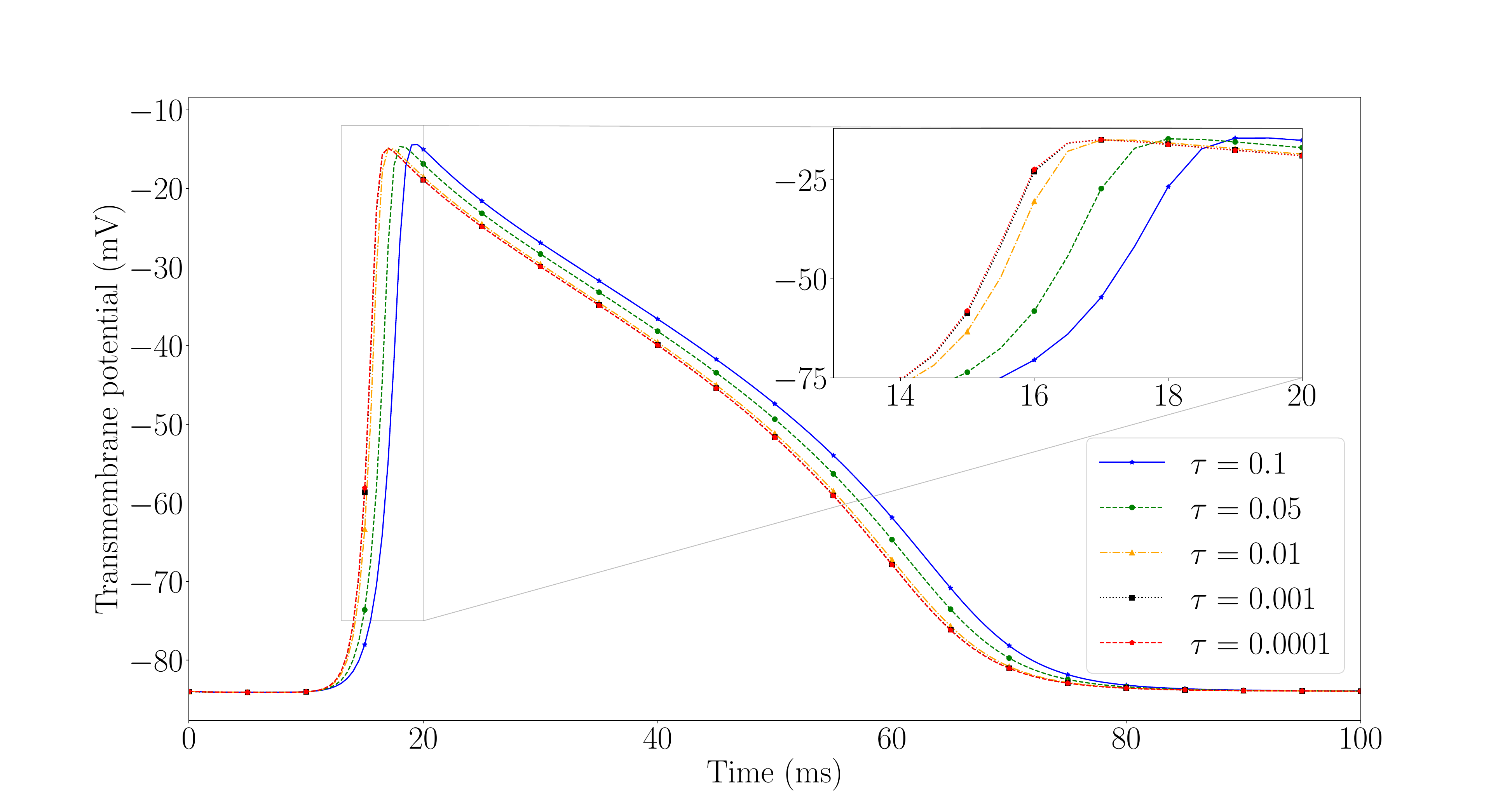}
      \caption{Transmembrane potential $v$ over time period $100\text{ ms}$ at the midpoint of the grid for various time steps.}
\end{subfigure}
\caption{Convergence in time for the decoupled strategy with IMEX scheme using LR91 model. }
\label{fig:lr_imex}
% \end{figure}
% \begin{figure}
% \centering
\begin{subfigure}[t]{0.495\textwidth}
     \includegraphics[width=\linewidth]{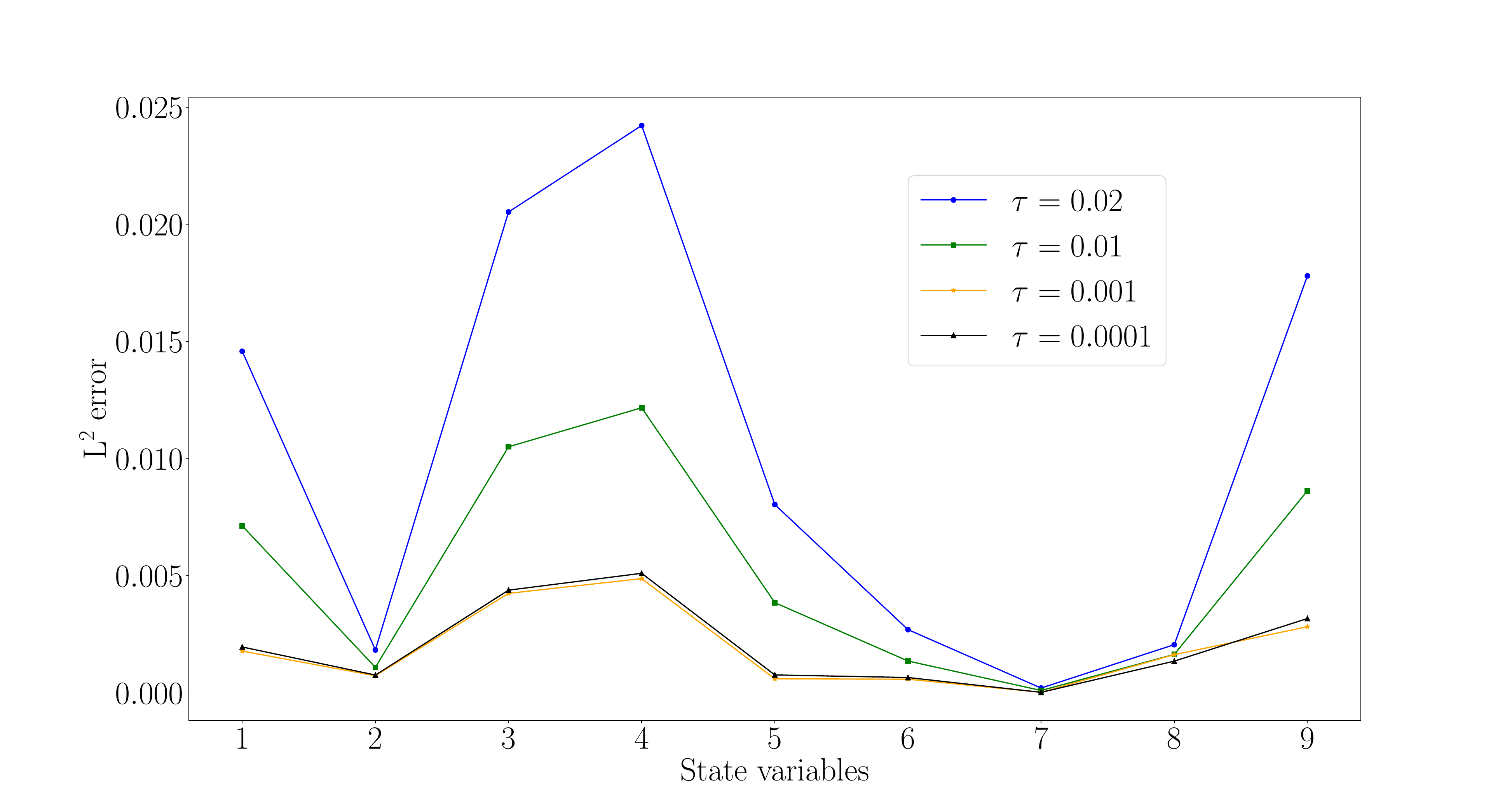}
      \caption{ Relative L$^2$ for the state variables at $t=20\text{ ms}$.}
    \end{subfigure}
\hfill
\begin{subfigure}[t]{0.495\textwidth}
   \includegraphics[width=\linewidth]{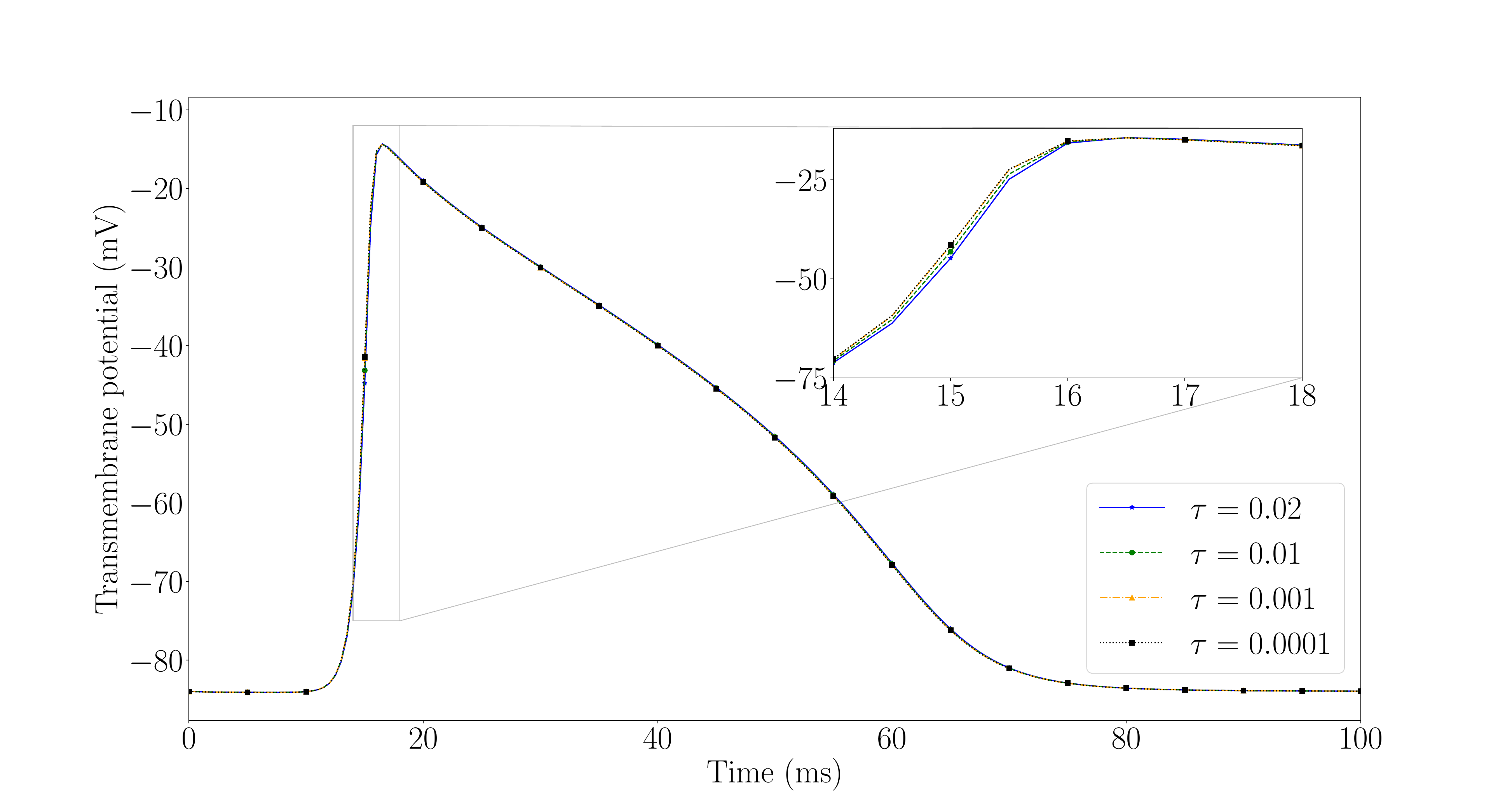}
      \caption{Transmembrane potential $v$ over time period $100\text{ ms}$ at the midpoint of the grid for various time steps.}
\end{subfigure}
\caption{Convergence in time for the decoupled strategy with CN-RK2 scheme using LR91 model.}
\label{fig:lr_cnrk}
\end{figure}

The convergence behaviour of various semi-implicit schemes in the decoupled strategy is evaluated for different time step sizes. The computational domain and grid resolution are identical to those used in the partitioned strategy analysis. Fig.~\ref{fig:lr_imex} shows the relative L$^2$ error of all the eight state variables at $t = 20\text{ ms}$ (left) and the transmembrane potential $v$ over a period of $100\text{ ms}$ (right) for the LR91 model using the IMEX scheme. For the IMEX scheme, time steps of $0.1$, $0.05$, $0.01$, $0.001$, and $0.0001$ are tested. For the CN-RK2 scheme, convergence is examined using $\tau = 0.02$, $0.01$, $0.001$, and $0.0001$, with larger time steps omitted due to stability constraints. The corresponding results obtained using the CN-RK2 scheme are presented in Fig.~\ref{fig:lr_cnrk}. Both IMEX and CN-RK2 demonstrate convergence at $\tau = 0.001$. , which is therefore adopted as the time step for all subsequent computations in this manuscript. Notably, the $L^2$ error associated with the CN-RK2 method is significantly smaller than that of the IMEX scheme at the same time step, despite both methods demonstrating solution convergence. Other schemes, including IMEX Gear, CN-AB, CN-FE, CN-IMEX, and SBDF, exhibit similar accuracy and convergence behaviour as the IMEX scheme.
\begin{figure}
\centering
\begin{subfigure}[t]{0.495\textwidth}
      \includegraphics[width=\linewidth]{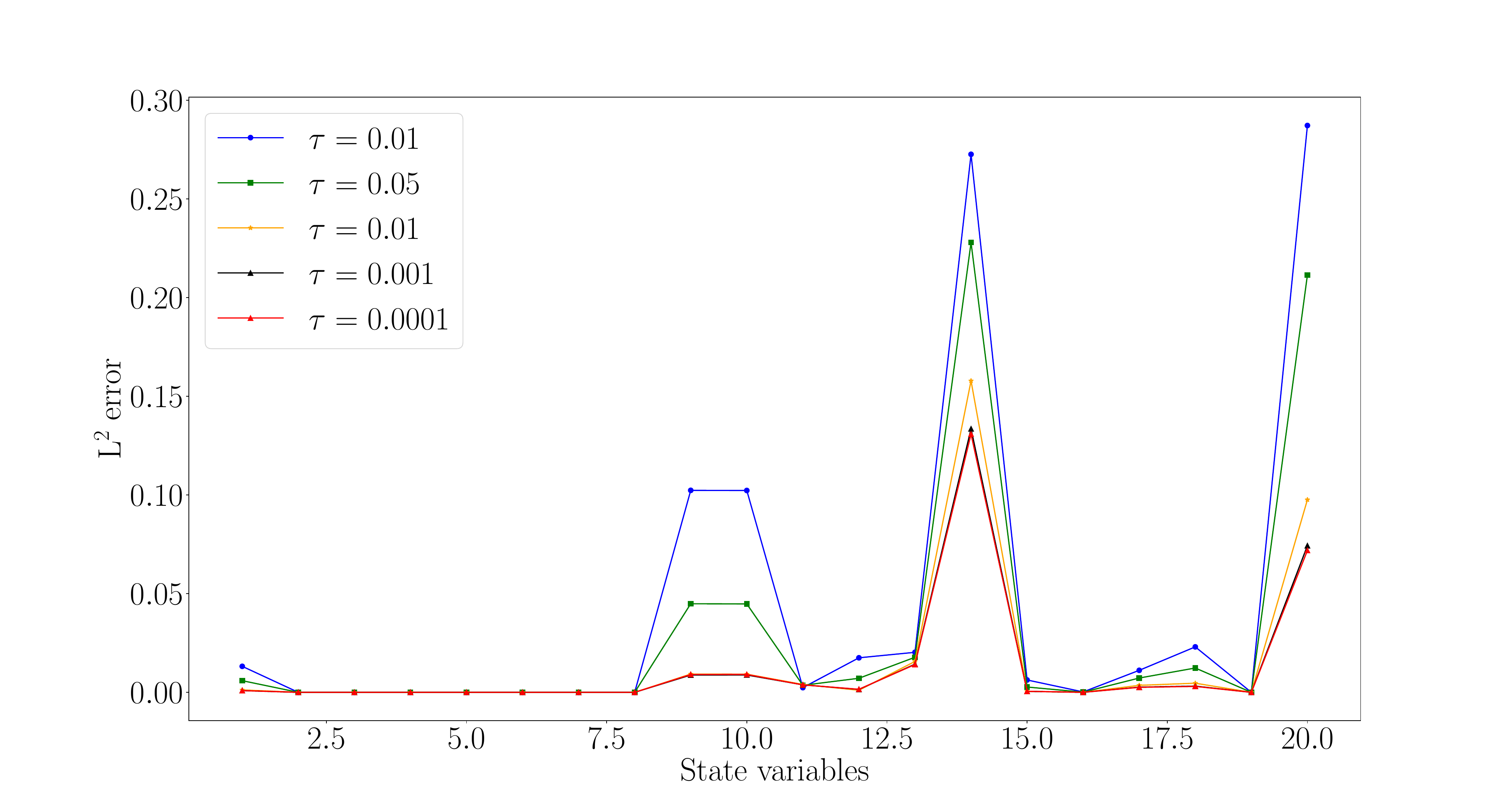}
      \caption{Relative L$^2$ for the state variables at $t=60\text{ ms}$.}
    \end{subfigure}
\hfill
\begin{subfigure}[t]{0.495\textwidth}
      \includegraphics[width=\linewidth]{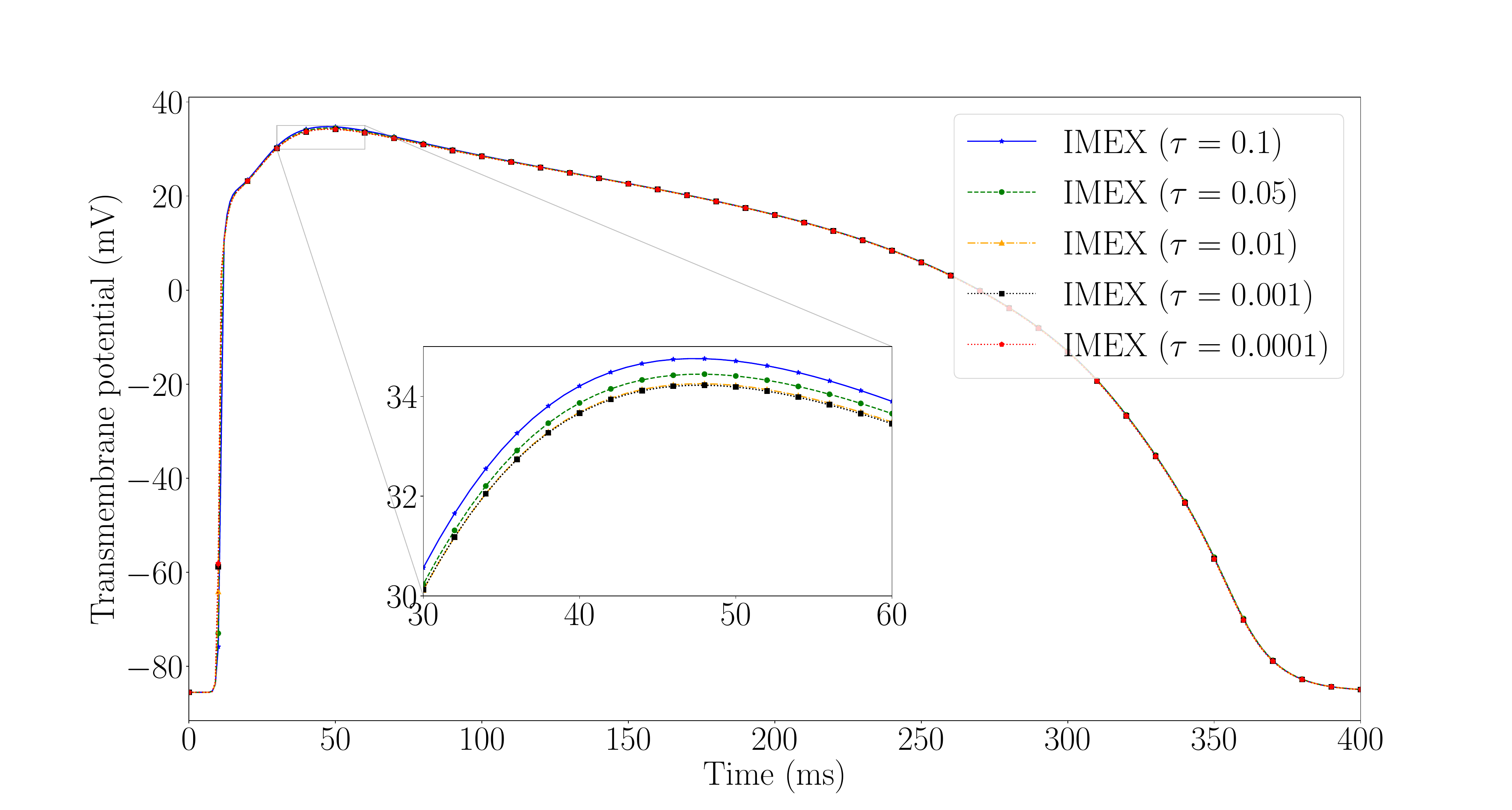}
      \caption{Transmembrane potential $v$ over time period $400\text{ ms}$ at the midpoint of the domain.}
\end{subfigure}
\caption{Convergence in time for the decoupled strategy with IMEX scheme using TP06 model.}
\label{fig:tp_imex}
\end{figure}
\begin{figure}
\centering
\begin{subfigure}[t]{0.495\textwidth}
      \includegraphics[width=\linewidth]{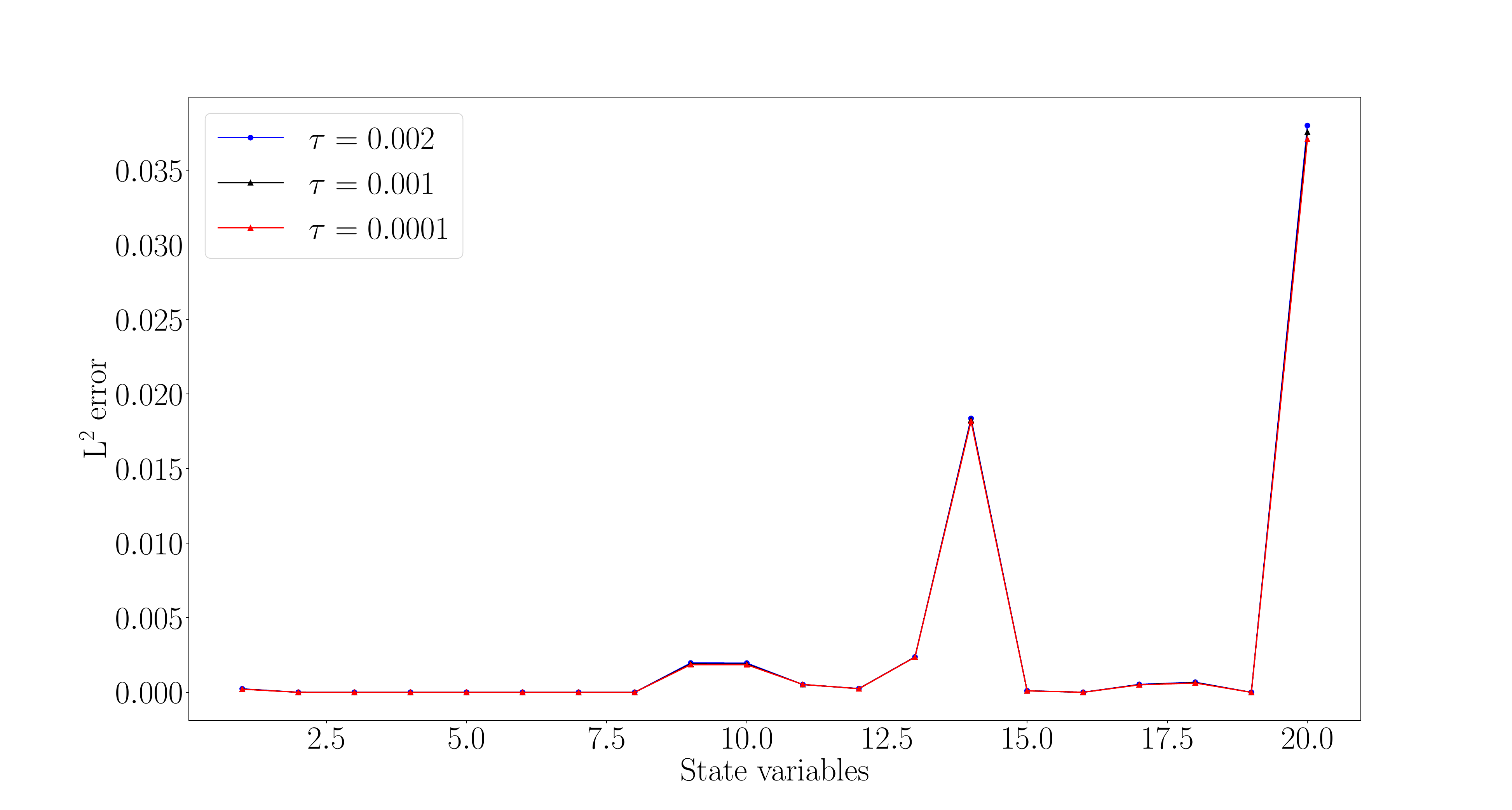}
      \caption{Relative L$^2$ for the state variables at $t=60\text{ ms}$.}
    \end{subfigure}
\hfill
\begin{subfigure}[t]{0.495\textwidth}
     \includegraphics[width=\linewidth]{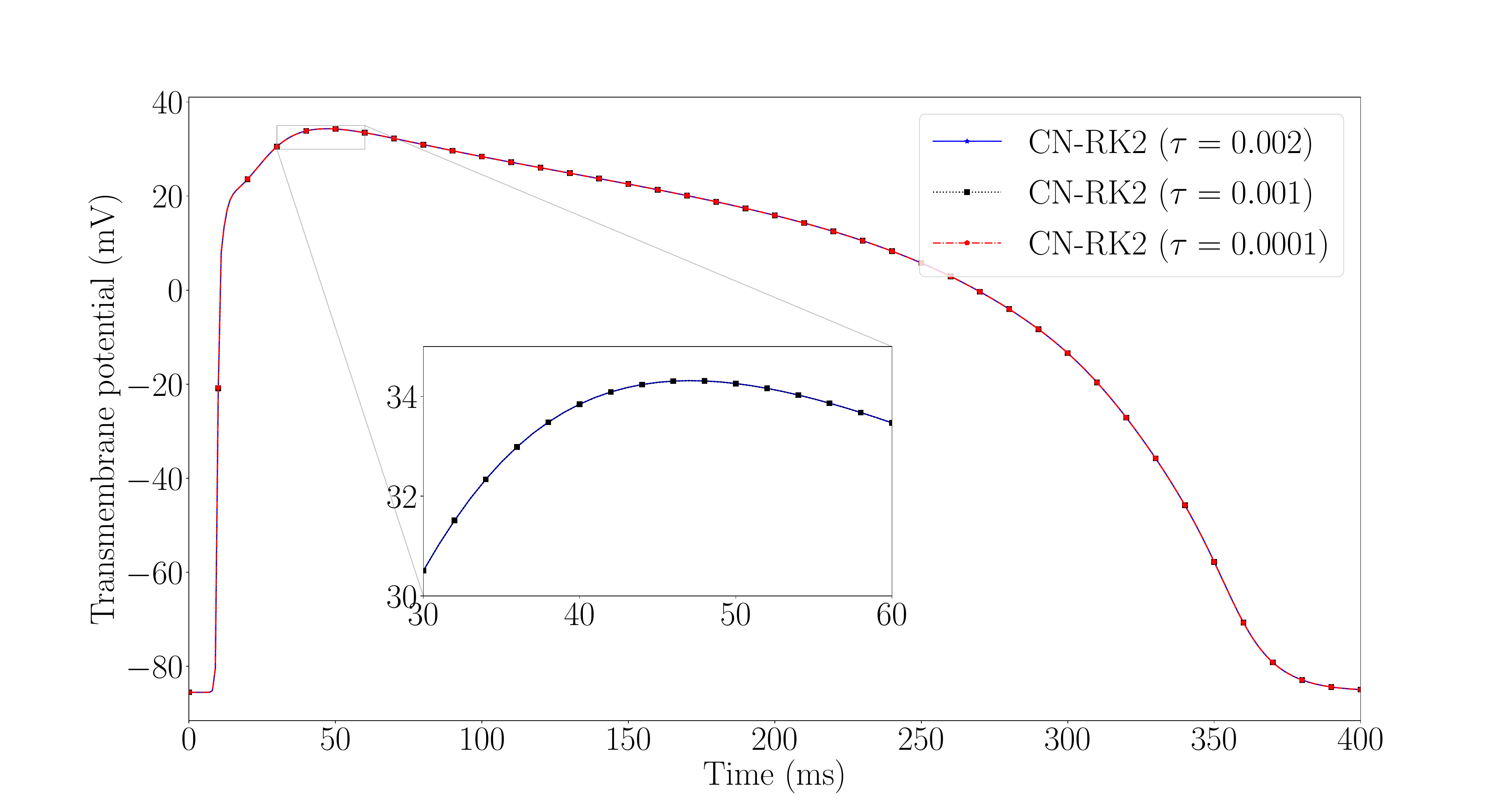}
      \caption{Transmembrane potential $v$ over time period $400\text{ ms}$ at the midpoint of the domain.}
\end{subfigure}
\caption{Convergence in time for the decoupled strategy with CN-RK2 scheme using TP06 model.}
\label{fig:tp_cnrk}
\end{figure}

In the case of the TP06 model, time steps of $0.1$, $0.05$, $0.01$, $0.001$, and $0.0001$ are tested for the IMEX scheme. For the CN-RK2 scheme, we could only consider time steps $\tau = 0.002$, $0.001$, and $0.0001$ due to stability issues of the solution. The relative L$^2$ error is computed at $t=60\text{ ms}$. Fig.~\ref{fig:tp_imex} and Fig.~\ref{fig:tp_cnrk} show the L$^2$ error for all the twenty state variables and the time evolution of $v$ for $400\text{ ms}$ using the IMEX and CN-RK2 schemes, respectively. Similar to the observations made for the LR91 model, both schemes exhibit convergence at $\tau = 0.001$ for the TP06 model as well, confirming the stability and accuracy of the solutions of the decoupled strategy across different time stepping schemes. A similar trend is also observed in the $L^2$ error: the CN-RK2 method yields a significantly smaller error compared to the IMEX scheme at the same time step, despite both methods achieving convergence.

\subsection{Comparison of different strategies}
\begin{figure}
 \begin{minipage}[]{0.4\textwidth}
    \centering
  \begin{tabular}{|p{2.1cm}|c|c|}
    \hline
    \multirow{2}{*}{Strategy} & \multicolumn{2}{c|}{$\text{L}^2$ error} \\
    \cline{2-3}
        &  $u$ & $v$ \\
    \hline
    % Coupled (Backward Euler ) &0.356378 & 0.711674\\
    Fully coupled & 0.01618 & 0.03385\\
    Partitioned without SDC &1.07751 &2.13239\\
    Partitioned with SDC &0.02527 &0.05177\\
    IMEX  & 2.22728 &4.48863\\
    IMEX Gear  & 2.22507 &4.48411\\
    CN-RK2  & 0.04587 &0.09216\\
    SBDF  & 2.24306 & 4.52046\\
    CN-AB  & 2.18353 & 4.40159\\
    CN-FE  & 2.15594 & 4.34586\\
    CN-IMEX  & 2.19309& 4.42086\\
    \hline
    \end{tabular}
      \captionof{table}{Absolute $\text{L}^2$ error for different strategies using LR91 model at $t=20\text{ ms}$.}
      \label{tab:l2error2}
    \end{minipage}
      \hfill
     \begin{minipage}[]{0.58\textwidth}
    \centering
     \includegraphics[trim={2cm 0cm 2cm 1cm},clip,width=\linewidth]{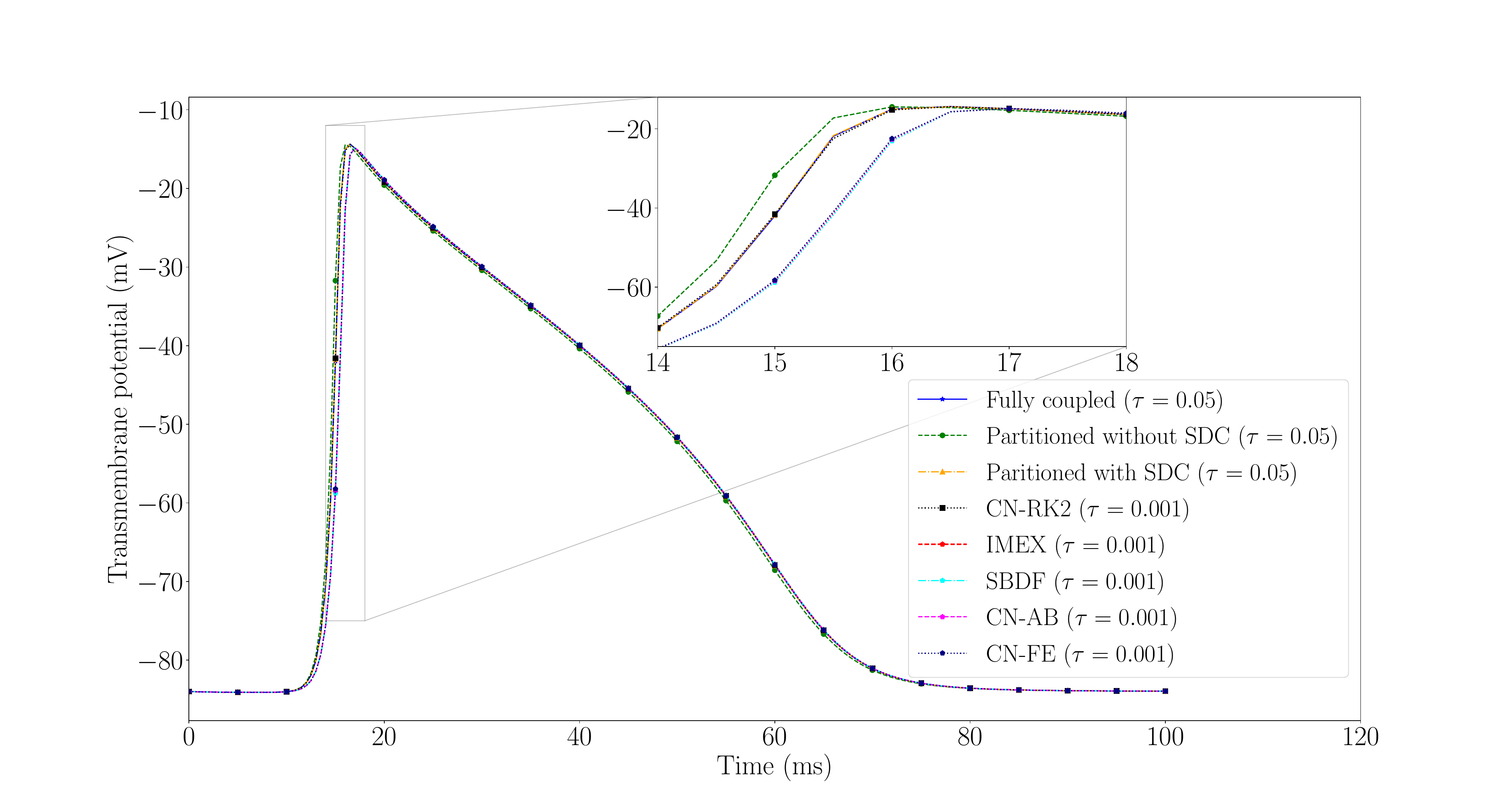}
    \captionof{figure}{Transmembrane potential ($v$) using LR91 model for different strategies for a time period of $100\text{ ms}$ at the midpoint of the computational domain.}
    \label{fig:lr}
    \end{minipage}
\end{figure}

In this subsection, the convergence and accuracy of the bidomain model solutions are compared across the coupled, partitioned, and various decoupled strategies for two physiologically realistic cell models. The coupled and partitioned strategies use a fixed time step of $\tau=0.05$. A time step size of $\tau=0.001$ is considered for the decoupled strategies. 
These time step sizes are selected based on the results of the previous convergence analysis for the three strategies. The L$^2$ error of the solution from the partitioned and decoupled strategies is computed by using the reference solution obtained using the fully coupled strategy with a time step size of $\tau = 0.001$. The absolute error for the LR91 model is evaluated at $20 \text{ms}$, corresponding to the onset of the repolarization phase, and measured at the midpoint of the computational domain. Table~\ref{tab:l2error2} presents the absolute errors in the transmembrane potential $v$ and the extracellular potential $u$ across all numerical strategies. It is evident that the fully coupled approach, the partitioned method with SDC, and the CN-RK2 decoupling strategy yield smaller $L^2$ errors, whereas the other methods result in significantly larger errors. These error metrics are essential for achieving highly accurate numerical solutions. Although such differences may appear negligible in single action potential duration (APD) simulations, where the errors are often averaged out over time in the healthy heart simulations, they become critically important in more complex scenarios, such as reentry simulations, which will be discussed in the following subsection.
% \begin{table}
%     \centering
%     \begin{tabular}{|l|l|l|}
%     \hline
%     Schemes    & $\text{L}^2$ error ($u$) &$\text{L}^2$ error ($v$) \\
%     \hline
%     % Coupled (Backward Euler ) &0.356378 & 0.711674\\
%     Fully coupled & 0.0161834 & 0.0338488\\
%     Partitioned without SDC &1.07751 &2.13239\\
%     Partitioned with SDC &0.0252667 &0.0517726\\
%     % IMEX ($\tau$=0.05) & 5.72461 &11.4682\\
%     % IMEX ($\tau$=0.01) & 3.07635 &6.18543\\
%     % CN-IMEX ($\tau$=0.01) & 2.78837& 5.61455\\
%     % CN-RK2 ( $\tau$=0.01) & 0.178415 & 0.355724\\
%     IMEX ($\tau$=0.001)  & 2.22728 &4.48863\\
%     IMEX Gear ($\tau$=0.001) & 2.22507 &4.48411\\
%     CN-RK2 ($\tau$=0.001) & 0.0458756 &0.0921597\\
%     SBDF ($\tau$=0.001) & 2.24306 & 4.52046\\
%     % SBDF3 ($\tau$=0.001) & 2.24297& 4.52031\\
%     CN-AB ($\tau$=0.001) & 2.18353 & 4.40159\\
%     CN-FE ($\tau$=0.001) & 2.15594 & 4.34586\\
%     CN-IMEX ($\tau$=0.001) & 2.19309& 4.42086\\
%     \hline
%     \end{tabular}
%     \caption{Absolute $\text{L}^2$ error for different strategies and various time-stepping schemes using LR91 model at $t=20\text{ ms}$.}
%     \label{tab:l2error2}
% \end{table}
% \begin{figure}
% \centering
%   \includegraphics[trim={2cm 0cm 2cm 1cm},clip,scale=0.2]{images/lr91/lr.pdf}
% \caption{Transmembrane potential ($v$) using LR91 model for different strategies for a time period of $100\text{ ms}$ at the midpoint of the computational domain.}
% \label{fig:lr}
% \end{figure}
   
\begin{figure}
 \begin{minipage}[]{0.41\textwidth}
    \centering
  \begin{tabular}{|p{2.1cm}|c|c|}
    \hline
    \multirow{2}{*}{Strategy} & \multicolumn{2}{c|}{$\text{L}^2$ error} \\
    \cline{2-3}
        &  $u$ & $v$ \\
    \hline
    Fully coupled &0.005786 & 0.011872\\
    Partitioned without SDC &0.030982& 0.006729\\
    Partitioned with SDC  &0.007919& 0.017790\\
    % IMEX (0.1) & 0.130943 &0.439306\\
    % IMEX ($\tau$=0.01) & 0.0899852 &0.153482\\
    IMEX   & 0.029546 &0.035071\\
    IMEX-Gear  & 0.029508 &0.035233\\  
    SBDF  & 0.029900&0.035377\\
    CN-AB &0.029930 & 0.034626\\
    CN-FE  & 0.029584 & 0.034332\\
    CN-IMEX & 0.029914 &0.034792\\
    CN-RK2 & 0.007484& 0.017738\\
    \hline
    \end{tabular}
      \captionof{table}{Absolute L$^2$ error for different strategies using TP06 model at $t=60\text{ ms}$.}
      \label{tab:tp06error}
    \end{minipage}
      \hfill
     \begin{minipage}[]{0.57\textwidth}
    \centering
     \includegraphics[trim={2cm 0cm 2cm 1cm},clip,width=\linewidth]{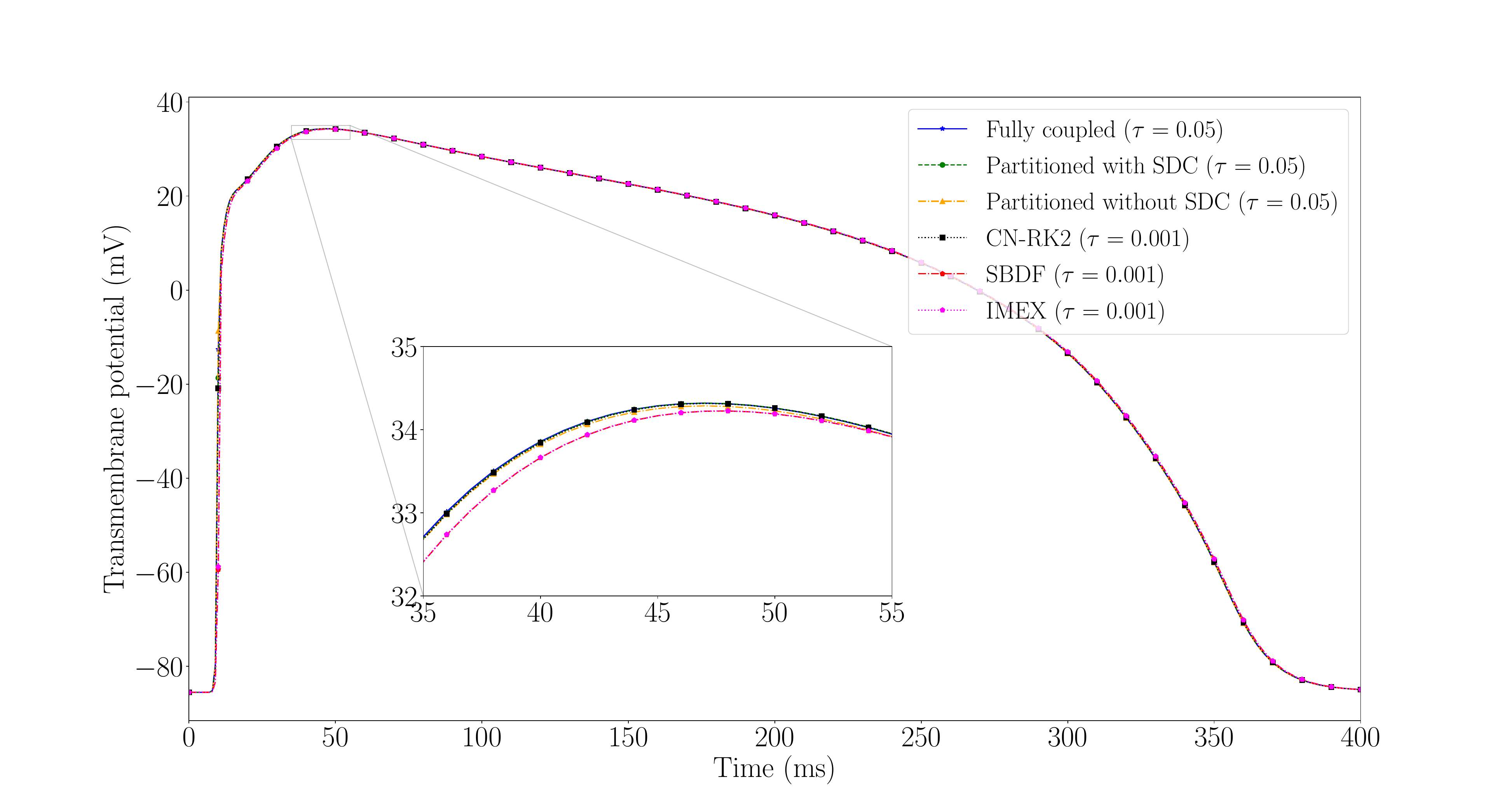}
    \captionof{figure}{Transmembrane potential ($v$) using TP06 model for different strategies for a time period of $400\text{ ms}$ at the point $(0.5,0.5)$, located at the midpoint of the computational domain.}
    \label{fig:tp}
    \end{minipage}
\end{figure}

 Similarly, Table~\ref{tab:tp06error} presents the corresponding errors for the TP06 model. For the TP06 model, the absolute error is computed at $60\text{ ms}$ after the depolarization that is measured at the midpoint of the computational domain. A similar behavior was observed for the TP06 model, consistent with the findings discussed above. To further support the above findings, Fig.~\ref{fig:lr} and Fig.~\ref{fig:tp} illustrate the temporal evolution of the transmembrane potential $v$ over a duration of $100 \text{ ms}$ for the LR91 model and $400\text{ ms}$ for the TP06 model. Across both ionic models, noticeable discrepancies are observed during the depolarization phase when using different numerical strategies. However, the solutions converge well in the resting state. The results presented in both the figures and the corresponding tables confirm that the partitioned strategy with SDC, as well as the CN-RK2 scheme, achieve good accuracy, closely replicating the behavior of the fully coupled reference solution.
% \begin{table}
%     \centering
%     \begin{tabular}{|l|l|l|}
%     \hline
%     Schemes    & L2 error ($u$) &L2 error ($v$) \\
%     \hline
%     Fully coupled &0.00578581 & 0.0118724\\
%     Partitioned without SDC &0.0309822& 0.00672941\\
%     Partitioned with SDC  &0.00791913& 0.0177903\\
%     % IMEX (0.1) & 0.130943 &0.439306\\
%     % IMEX ($\tau$=0.01) & 0.0899852 &0.153482\\
%     IMEX ($\tau$=0.001)  & 0.0295458 &0.0350712\\
%     IMEX Gear ($\tau$=0.001) & 0.0295075 &0.0352333\\  
%     SBDF ($\tau$=0.001) & 0.029897&0.0353773\\
%     % SBDF3 (dt=0.001) & \\not converging
%     CN-AB ($\tau$=0.001) &0.0299304 & 0.0346264\\
%     CN-FE ($\tau$=0.001) & 0.0295839 & 0.0343315\\
%     CN-IMEX ($\tau$=0.001)& 0.0299143 &0.0347921\\
%     CN-RK2 ($\tau$=0.001) & 0.0074843& 0.0177375\\
%     \hline
%     \end{tabular}
%     \caption{Absolute L$^2$ error for different strategies and various time-stepping schemes using TP06 model at $t=60\text{ ms}$.}
%     \label{tab:tp06error}
% \end{table}
% \begin{figure}
% \centering
%   \includegraphics[trim={2cm 0cm 2cm 1cm},clip,scale=0.2]{images/tp06/tp.pdf}
% \caption{Transmembrane potential ($v$) using TP06 model for different strategies for a time period of $400\text{ ms}$ at the point $(0.5,0.5)$, located at the midpoint of the computational domain.}
% \label{fig:tp}
% \end{figure}

\subsubsection*{Biomarkers}
We compare key action potential metrics, namely, the action potential duration at $90\%$ repolarization (APD90, in ms), the maximum upstroke velocity ($dv/dt$, in V/s), and the conduction velocity (CV, in cm/s), across the fully coupled, partitioned, and decoupled strategies. For this comparison, the fully coupled strategy employs a fixed time step of $0.05$, while partitioned strategies with and without SDC use adaptive time stepping. And, for the decoupled strategy, we compare the CN-RK2 and SBDF schemes with a fixed time step of $0.001$.
These characteristics for the potential with the LR91 model at the midpoint of the domain are presented in Table~\ref{tab:lr_bio}. Similarly, Table~\ref{tab:tp_bio} presents the action potential characteristics for the TP06 model. The results in both tables indicate that the action potential characteristics obtained using the partitioned strategy with SDC and the CN-RK2 scheme closely match those of the fully coupled strategy, whereas small deviations are observed with the partitioned strategy without SDC and the SBDF scheme.
\begin{figure}
 \begin{minipage}[]{0.49\textwidth}
    \centering
  \begin{tabular}{|p{2.1cm}|p{1.15cm}|p{1cm}|p{1.1cm}|}
    \hline
       Strategy & APD90 (ms) & $dv /dt$ (V/s) & CV (cm/s)\\
       \hline
        Fully coupled &47.720 &36.297 & 27.53 \\
        Partitioned without SDC&47.649 &38.321 & 28.04 \\ 
        Partitioned with SDC&47.692 &36.547 & 27.52\\ %
        CN-RK2 &47.727 &35.656 & 27.52\\ %
        SBDF&47.390 & 34.782 & 26.55\\ 
        \hline
    \end{tabular}
      \captionof{table}{Biomarkers for LR91 model}
      \label{tab:lr_bio}
    \end{minipage}
      \hfill
     \begin{minipage}[]{0.49\textwidth}
    \centering
   \centering
    \begin{tabular}{|p{2.1cm}|p{1.15cm}|p{1cm}|p{1.1cm}|}
    \hline
       Strategy & APD90 (ms) & $dv /dt$ (V/s) & CV (cm/s)\\
       \hline
       Fully coupled &353.509 &44.011 &44.26\\
        Partitioned without SDC&353.492 &44.081 & 43.86  \\ 
        Partitioned with SDC&353.548 &44.296 &44.01 \\ %
        CN-RK2&353.560 &44.209 & 44.56 \\ %
        SBDF&353.763 &43.431 &45.45 \\ 
        \hline
    \end{tabular}
    \captionof{table}{Biomarkers for TP06 model}
    \label{tab:tp_bio}
 \end{minipage}
\end{figure}
% \begin{table}
%     \centering
%     \begin{tabular}{|l|c|c|c|}
%     \hline
%        Strategy  & APD90 (ms) & $dv /dt$ (V/s) & CV (cm/s)\\
%        \hline
%        Fully Coupled &47.720 &36.297 & 27.53 \\
%         Partitioned without SDC&47.649 &38.321 & 28.04 \\ 
%         Partitioned with SDC&47.692 &36.547 & 27.52\\ %
%         CN-RK2 scheme&47.727 &35.656 & 27.52\\ %
%         SBDF scheme&47.390 & 34.782 & 26.55\\ 
%         \hline
%     \end{tabular}
%     \caption{Biomarkers for LR91 model}
%     \label{tab:lr_bio}
% \end{table}
% \begin{table}
%     \centering
%     \begin{tabular}{|l|c|c|c|}
%     \hline
%        Strategy  & APD90 (ms) & $dv /dt$ (V/s) & CV (cm/s)\\
%        \hline
%        Fully Coupled &353.509 &44.011 &44.26\\
%         Partitioned without SDC&353.492 &44.081 & 43.86  \\ 
%         Partitioned with SDC&353.548 &44.296 &44.01 \\ %
%         CN-RK2 scheme&353.560 &44.209 & 44.56 \\ %
%         SBDF scheme&353.763 &43.431 &45.45 \\ 
%         \hline
%     \end{tabular}
%     \caption{Biomarkers for TP06 model}
%     \label{tab:tp_bio}
% \end{table}
% \begin{table}
\begin{wraptable}{R}{5cm}
    \centering
    \begin{tabular}{|p{2.1cm}|r|}
    \hline
     Strategy   & Time (s) \\
         \hline
    Fully coupled  & 65738.651\\
    Partitioned with SDC  &5200.755\\
    CN-RK2 & 20381.856    \\
    \hline
    \end{tabular}
    \caption{Total Time for solving with LR91 model for $t=100\text{ ms}$.}
    \label{tab:lrtime}
\end{wraptable}
% \end{table}
\subsubsection*{Computational Cost}
We compare the total computational times required to solve the bidomain model using fully coupled, partitioned, and decoupled numerical strategies. For this analysis, simulations based on the LR91 model are performed over a duration of $100\text{ ms}$ as shown in Fig.~\ref{fig:lr}—which captures a full APD. The corresponding sequential computational times for all strategies are reported in Table~\ref{tab:lrtime}. Specifically, we compare the performance of the fully coupled scheme and the partitioned strategy with SDC, both using a fixed time step of $\tau = 0.05$, against the Strang splitting approach (CN-RK2) with a smaller time step of $\tau = 0.001$. We employed the GMRES solver with the MRAS preconditioner across all strategies to ensure a fair comparison of sequential runtimes. In the fully coupled strategy, the solver was used to solve the linear systems in Newton iteration for the coupled PDE-ODE system. In contrast, in the partitioned and decoupled strategies, it was applied only to the discretised elliptic equation \eqref{eq:elliptic-ODE}. The results demonstrate that the partitioned strategy with SDC achieves substantial computational savings, being approximately 12.64 times faster than the fully coupled strategy and 3.92 times faster than the CN-RK2 scheme, while maintaining a better level of accuracy.

\subsection{Reentry}
\begin{figure}[]
    \centering
    \includegraphics[trim={2cm 0cm 2cm 1cm},clip,scale=0.2]{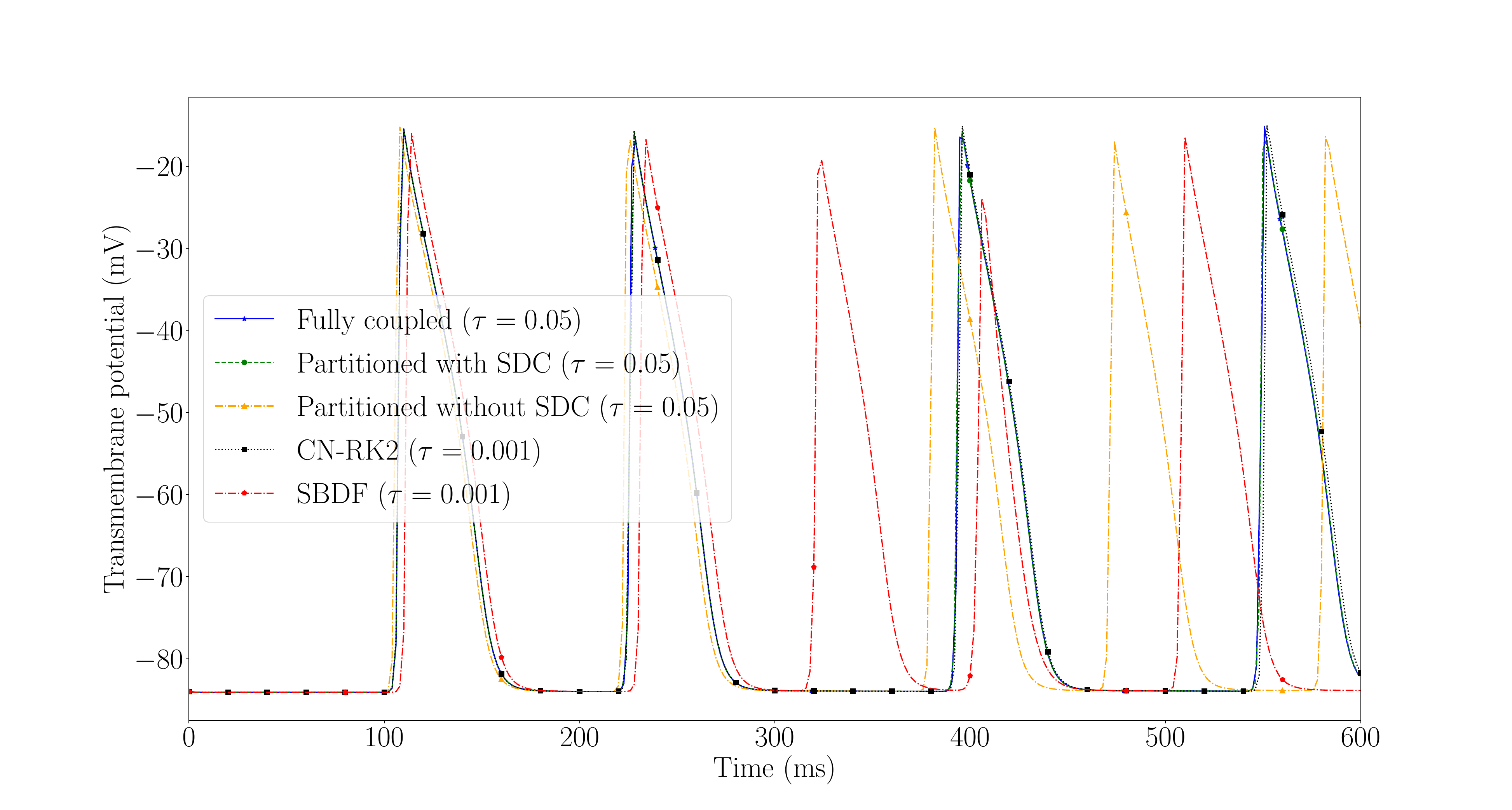}
    \caption{Transmembrane potential for LR91 model over time period $600\,\text{ms}$ at the midpoint of the domain.}
    \label{fig:lr91}
\end{figure}

\begin{figure}
    \centering
\scalebox{0.9}{
    \includegraphics[trim={2cm 2cm 4.8cm 4cm},clip,width=0.175\linewidth]{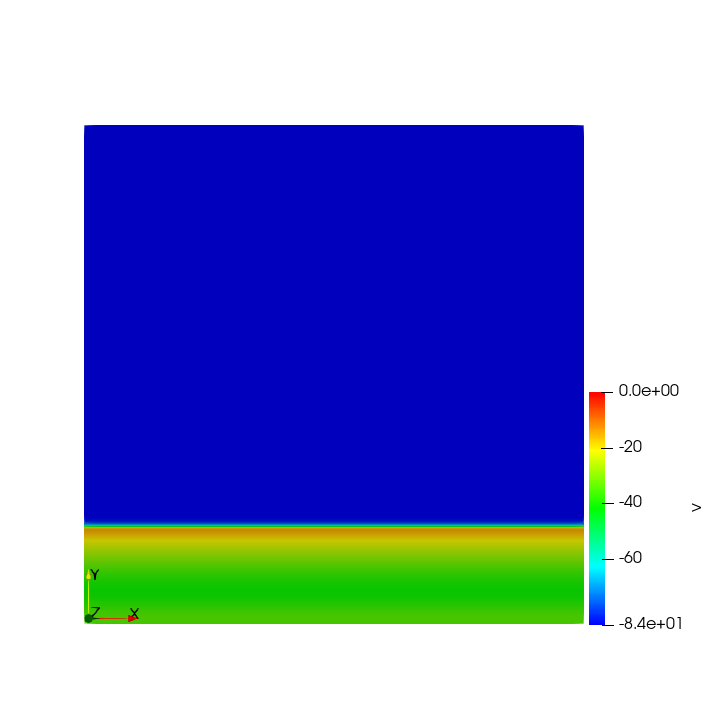}
    \includegraphics[trim={2cm 2cm 4.8cm 4cm},clip,width=0.175\linewidth]{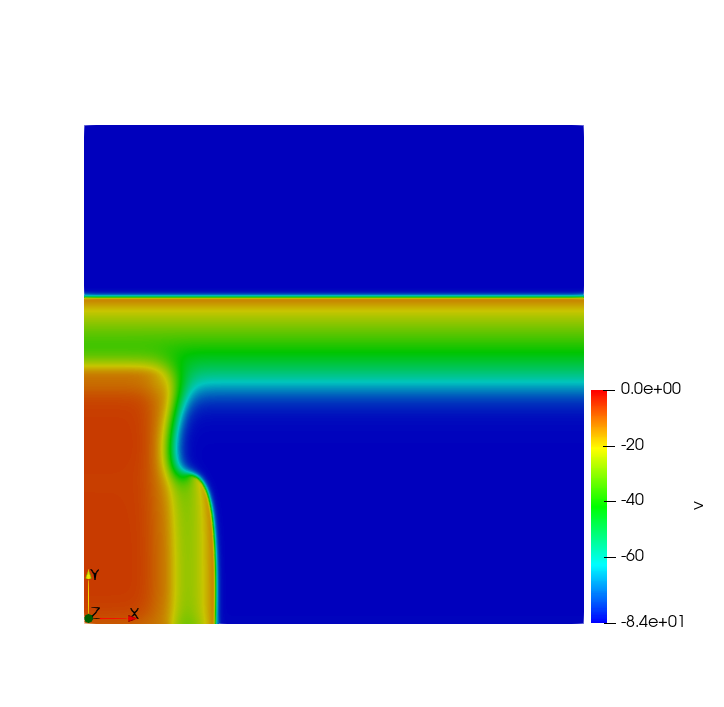}
    \includegraphics[trim={2cm 2cm 4.8cm 4cm},clip,width=0.175\linewidth]{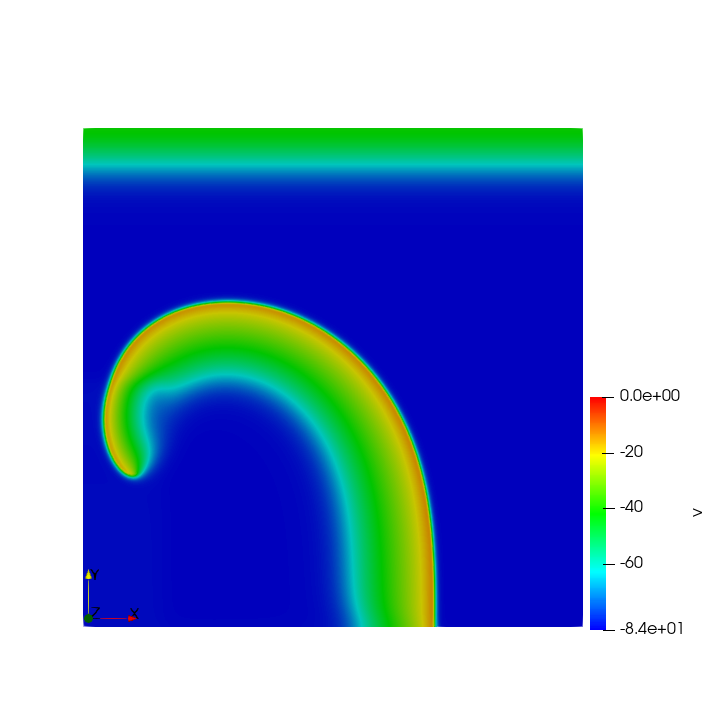}
    \includegraphics[trim={2cm 2cm 4.8cm 4cm},clip,width=0.175\linewidth]{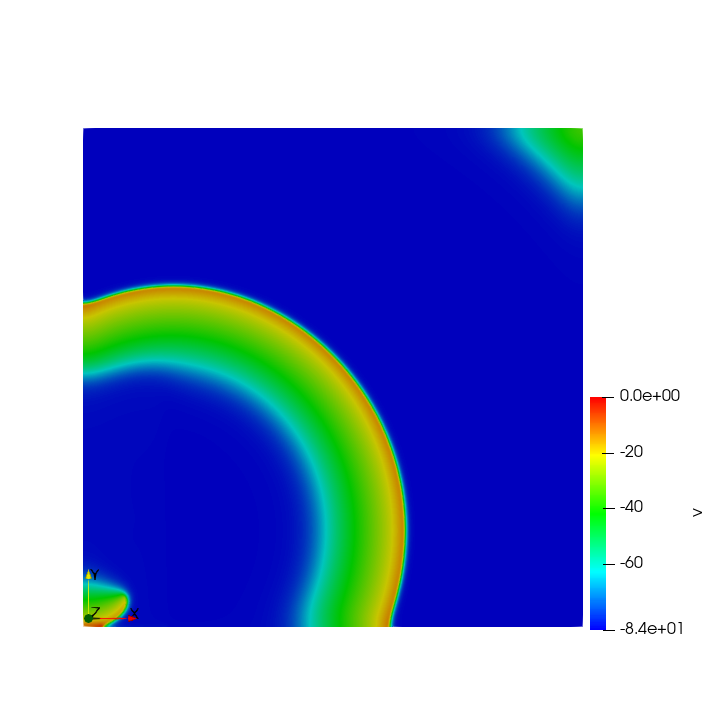}
    \includegraphics[trim={2cm 2cm 4.8cm 4cm},clip,width=0.175\linewidth]{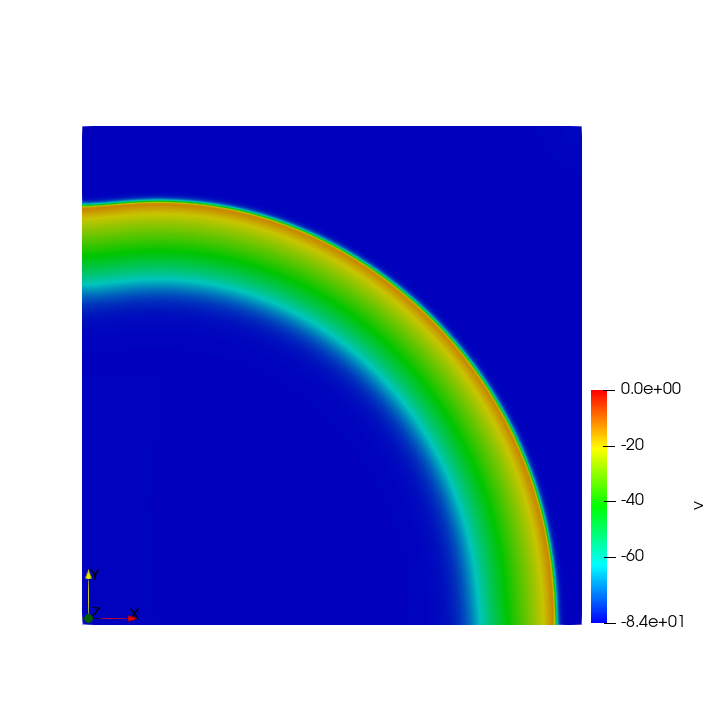}
    \includegraphics[trim={0.5cm 0cm 0.5cm 0cm},clip,height=1.09in]{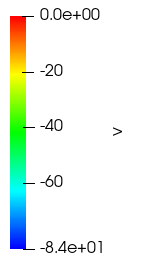}}
\scalebox{0.9}{
     \includegraphics[trim={2cm 2cm 4.8cm 4cm},clip,width=0.175\linewidth]{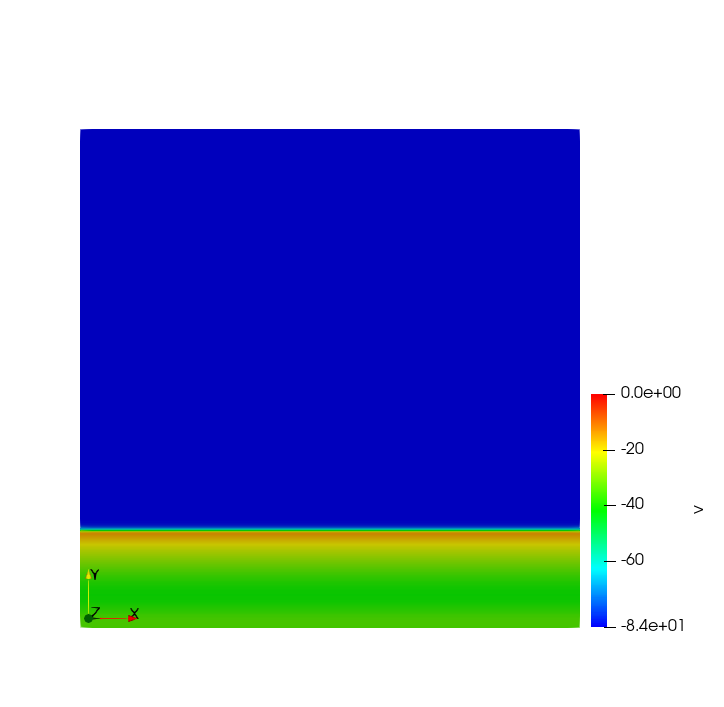}
    \includegraphics[trim={2cm 2cm 4.8cm 4cm},clip,width=0.175\linewidth]{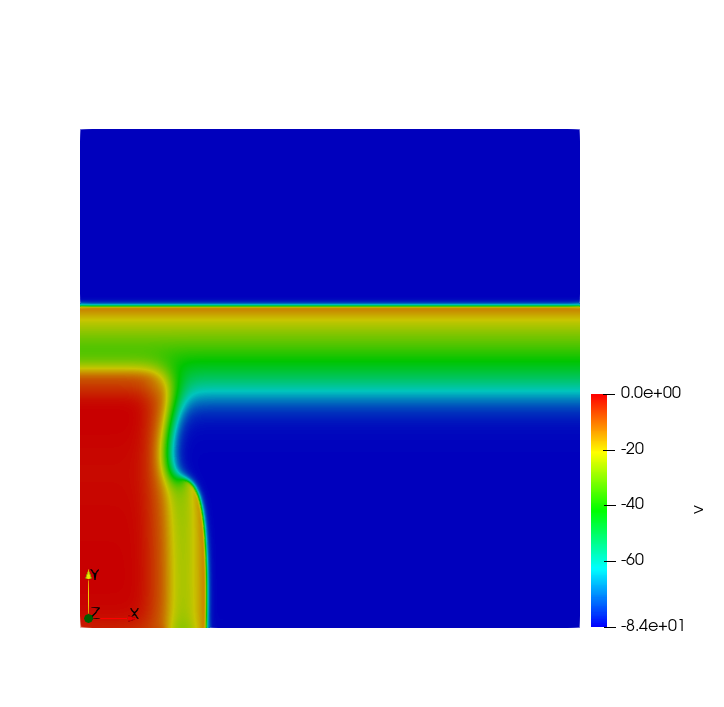}
    \includegraphics[trim={2cm 2cm 4.8cm 4cm},clip,width=0.175\linewidth]{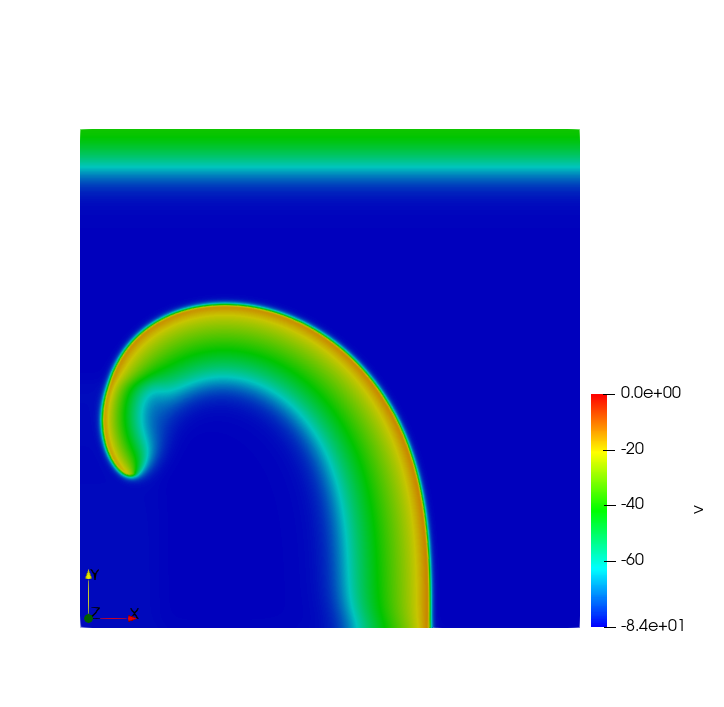}
    \includegraphics[trim={2cm 2cm 4.8cm 4cm},clip,width=0.175\linewidth]{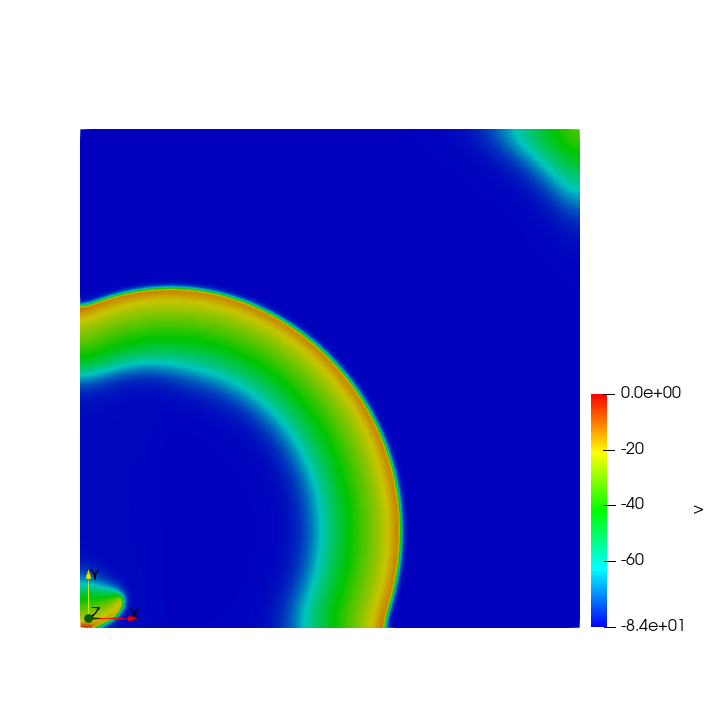}
    \includegraphics[trim={2cm 2cm 4.8cm 4cm},clip,width=0.175\linewidth]{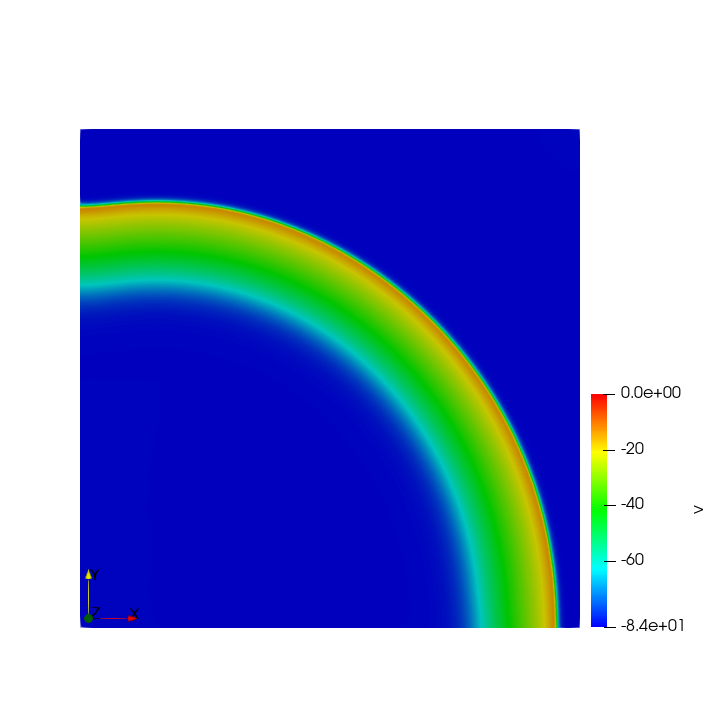}
    \includegraphics[trim={0.5cm 0cm 0.5cm 0cm},clip,height=1.09in]{images/lr91/re-entry/v_m.png}}
\scalebox{0.9}{
    \includegraphics[trim={2cm 2cm 4.8cm 4cm},clip,width=0.175\linewidth]{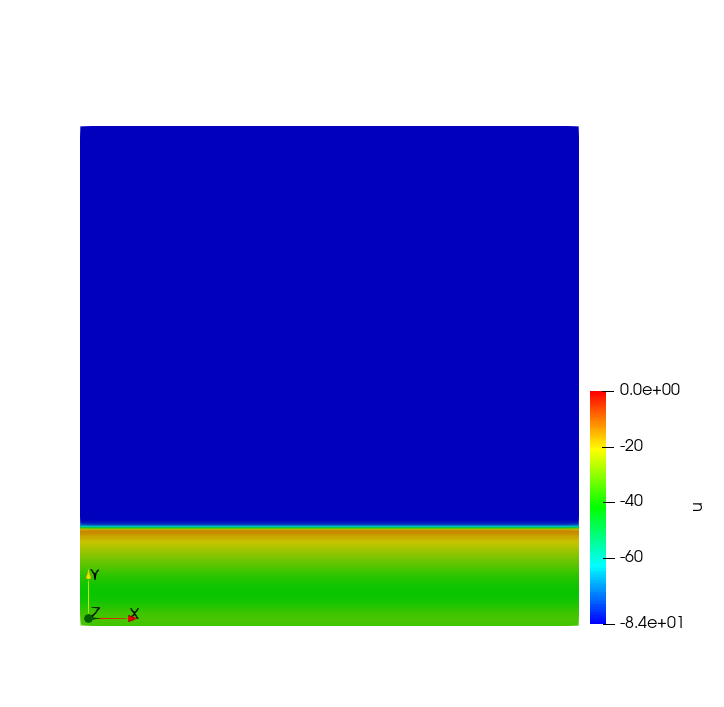}
    \includegraphics[trim={2cm 2cm 4.8cm 4cm},clip,width=0.175\linewidth]{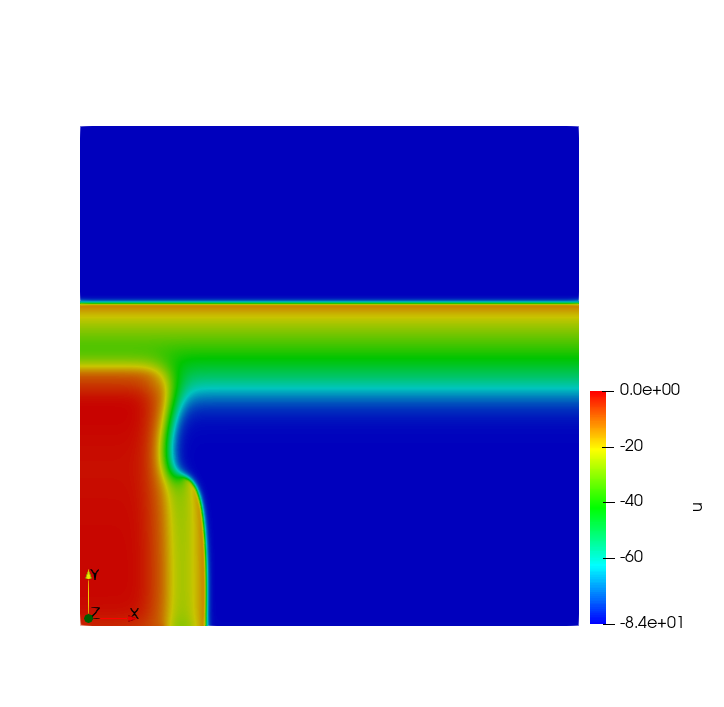}
    \includegraphics[trim={2cm 2cm 4.8cm 4cm},clip,width=0.175\linewidth]{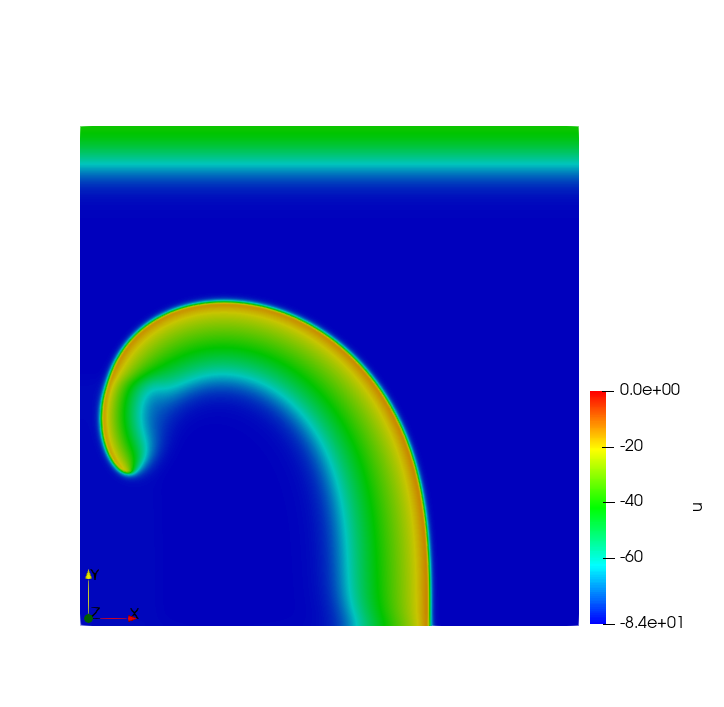}
    \includegraphics[trim={2cm 2cm 4.8cm 4cm},clip,width=0.175\linewidth]{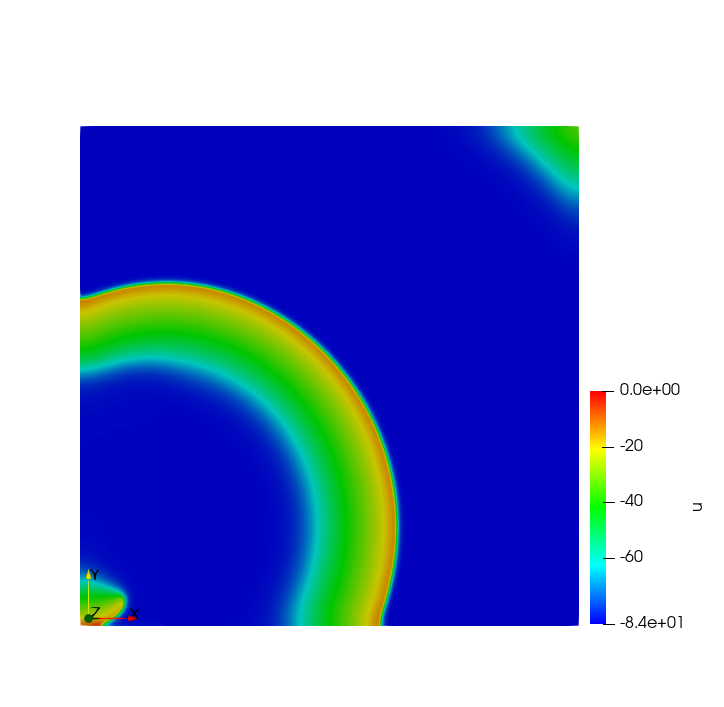}
    \includegraphics[trim={2cm 2cm 4.8cm 4cm},clip,width=0.175\linewidth]{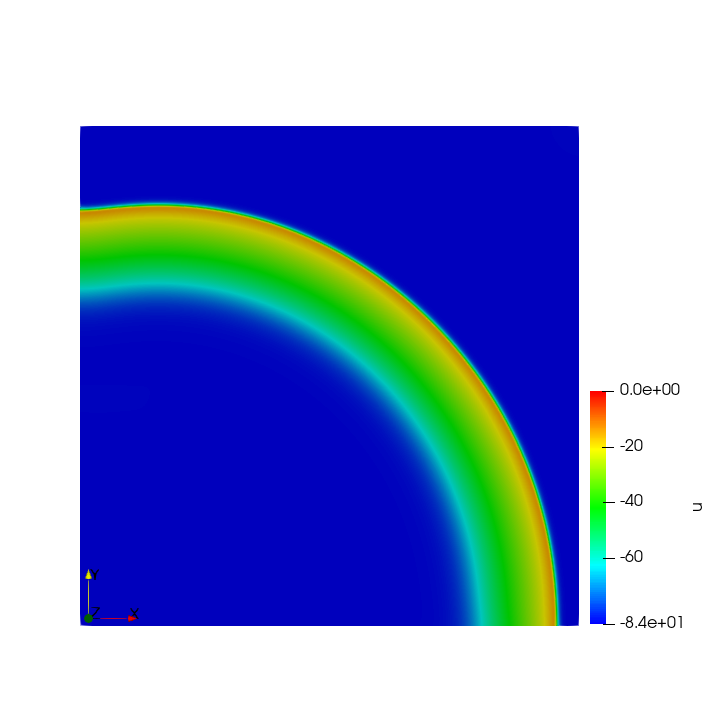}
    \includegraphics[trim={0.5cm 0cm 0.5cm 0cm},clip,height=1.09in]{images/lr91/re-entry/v_m.png}}
\scalebox{0.9}{
     \includegraphics[trim={2cm 2cm 4.8cm 4cm},clip,width=0.175\linewidth]{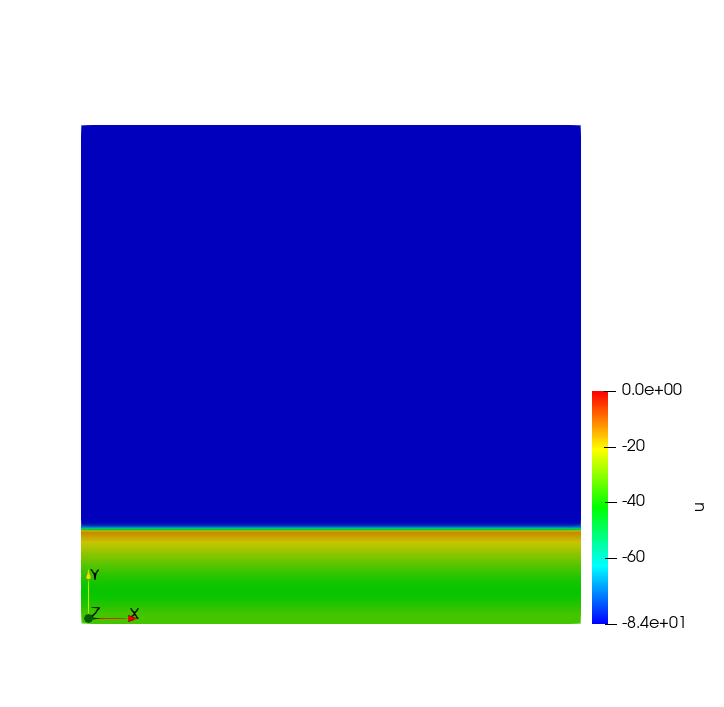}
    \includegraphics[trim={2cm 2cm 4.8cm 4cm},clip,width=0.175\linewidth]{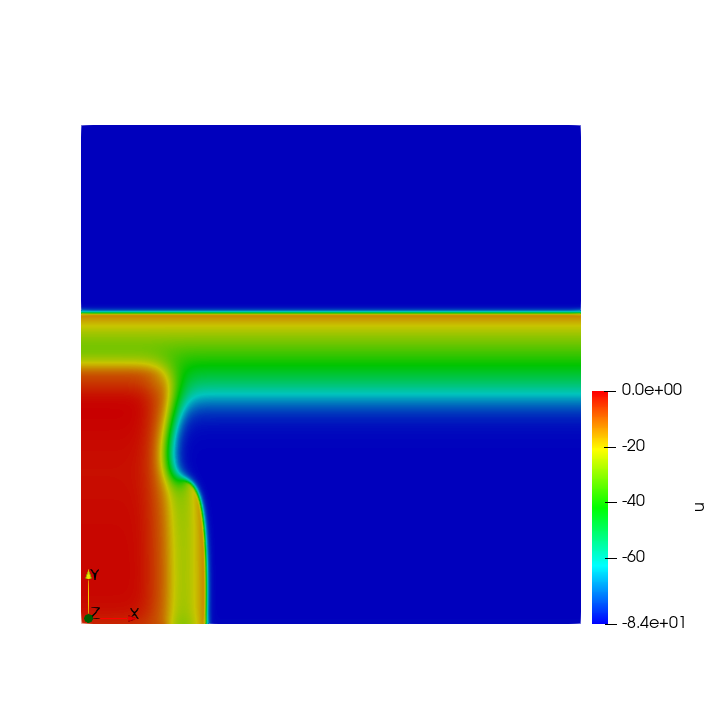}
    \includegraphics[trim={2cm 2cm 4.8cm 4cm},clip,width=0.175\linewidth]{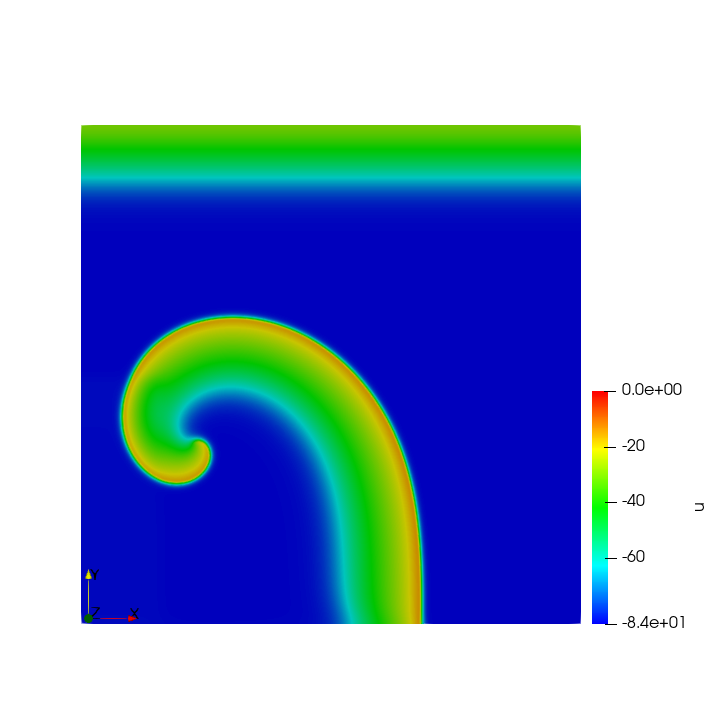}
    \includegraphics[trim={2cm 2cm 4.8cm 4cm},clip,width=0.175\linewidth]{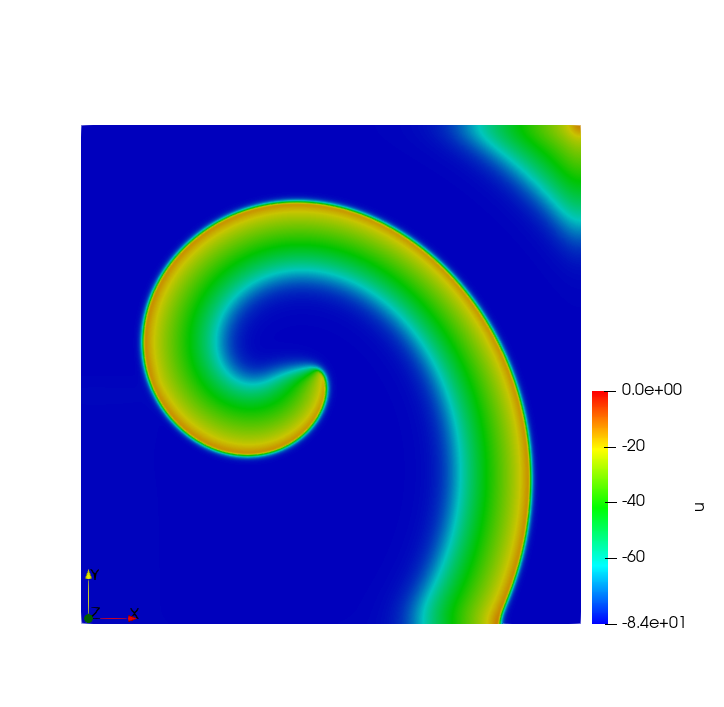}
    \includegraphics[trim={2cm 2cm 4.8cm 4cm},clip,width=0.175\linewidth]{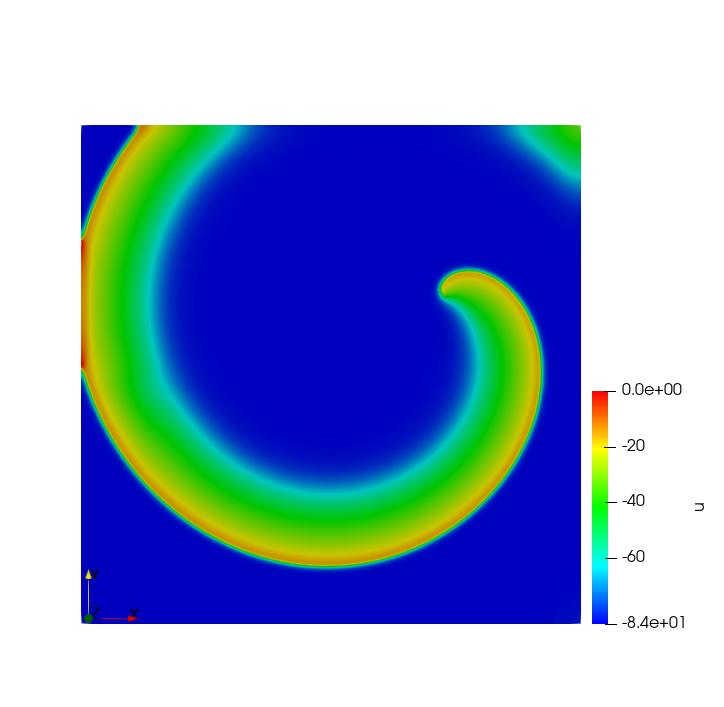}
    \includegraphics[trim={0.5cm 0cm 0.5cm 0cm},clip,height=1.09in]{images/lr91/re-entry/v_m.png}}
    \caption{Transmembrane potential $v$ for LR91 model at $t=40$, $140$, $240$, $400$, and $600$ ms respectively. \textit{First row:} Fully coupled strategy, \textit{Second row:} Partitioned strategy with SDC, \textit{Third row:} CN-RK2 scheme with $\tau=0.001$, and \textit{Fourth row:} SBDF scheme with $\tau=0.001$.}
    \label{fig:lr_2D}
\end{figure}

A domain $\Omega=[0,6]^2$ with grid resolution $600\times 600$ is chosen for the reentry analysis of the LR91 model. An initial stimulus $\I_{stim}=100\mu \text{A}/\text{cm}^2$ is applied for $2 \text { ms}$ at $t=0$ along the bottom edge of the computational domain. The second stimulus is applied at the tail of the excitation wavefront, specifically at $t=120 \text{ ms}$, in the region $[0,1]\times[0,3]$ to generate spiral wavefronts. Fig.~\ref{fig:lr91} shows the transmembrane potential $v$ over a simulation period of $600\,\text{ms}$ at the point $(3,3)$. The fully coupled strategy is computed with a fixed time step of $0.05$, the partitioned strategy uses an adaptive time-stepping approach with a maximum time step of $0.05$, and the decoupled strategy employs a fixed time step of $0.001$ which are required to get the accurate numerical solutions as like in single APD simulations. To further illustrate the spatiotemporal evolution of the solution, Fig.~\ref{fig:lr_2D} presents a sequence of snapshots of $v$ at $t = 40$, $140$, $240$, $400$, and $600 \text{ms}$, obtained using the fully coupled strategy, the partitioned strategy with SDC, the CN-RK2 scheme, and the SBDF scheme. It is evident from both Fig.~\ref{fig:lr91} and Fig.~\ref{fig:lr_2D} that the partitioned strategy with SDC and the CN-RK2 scheme produce results that closely match those obtained from the fully coupled strategy, consistent with the observations from the convergence tests. It can be observed that the solutions obtained using the partitioned method without SDC, as well as those computed with the SBDF scheme, do not align closely with the results from the fully coupled strategy. It is important to note, however, that these discrepancies diminish when the time step is further reduced, leading to improved agreement with the fully coupled reference solution. This behavior has also been consistently observed in multiple APD simulations conducted in the context of cardiac ventricular arrhythmia modeling, where numerical accuracy is particularly sensitive to the choice of time integration scheme and time step size.

\begin{figure}[]
    \centering
    \includegraphics[trim={2cm 0cm 2cm 1cm},clip,scale=0.2]{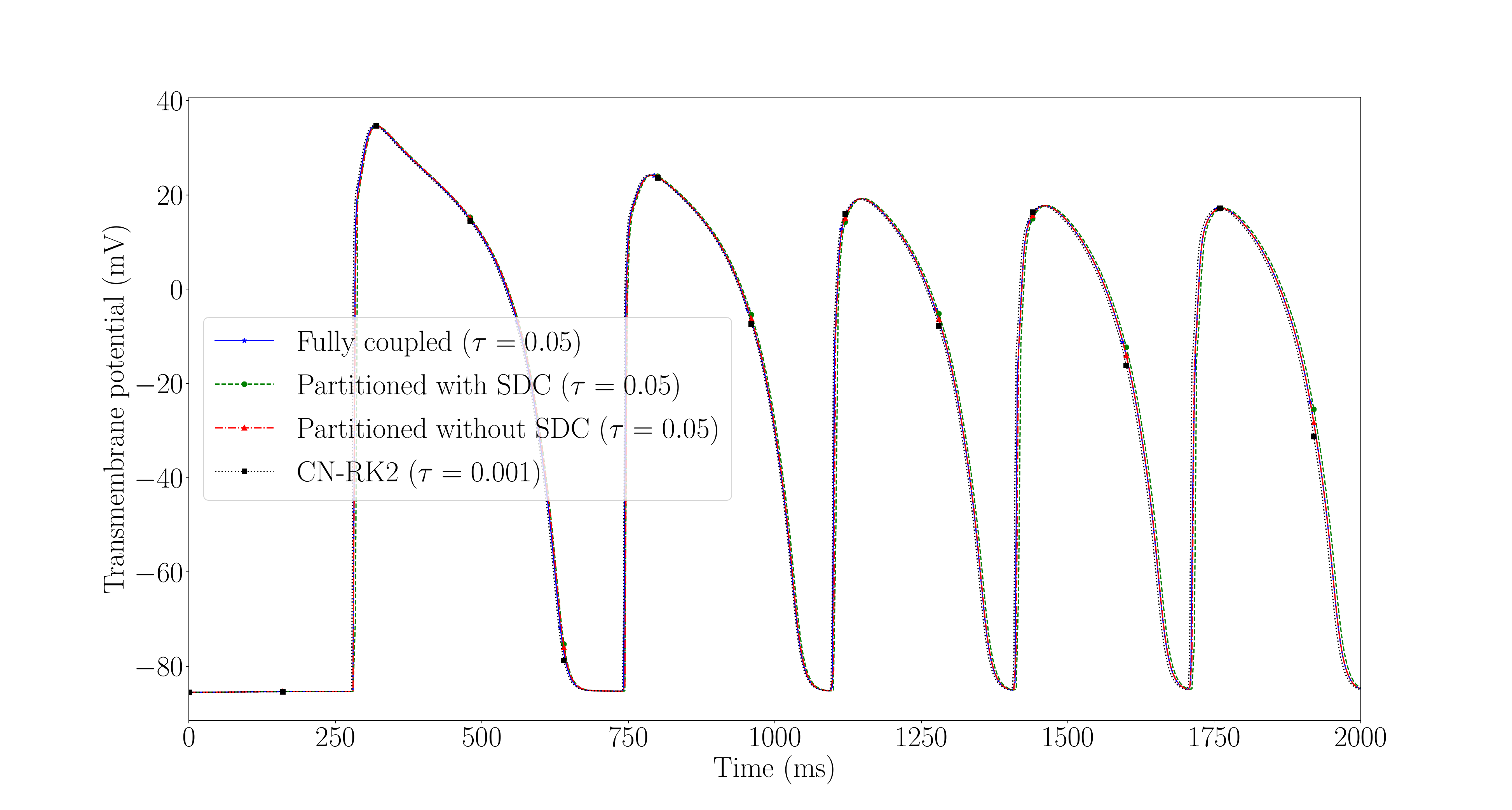}
    \caption{Transmembrane potential $v$ for TP06 model over a period of time $2000\text{ ms}$ at the midpoint of the domain.}
    \label{fig:tp06}
\end{figure}

\begin{figure}
    \centering
\scalebox{0.9}{
    \includegraphics[trim={2cm 2.5cm 6cm 4cm},clip,width=0.18\linewidth]{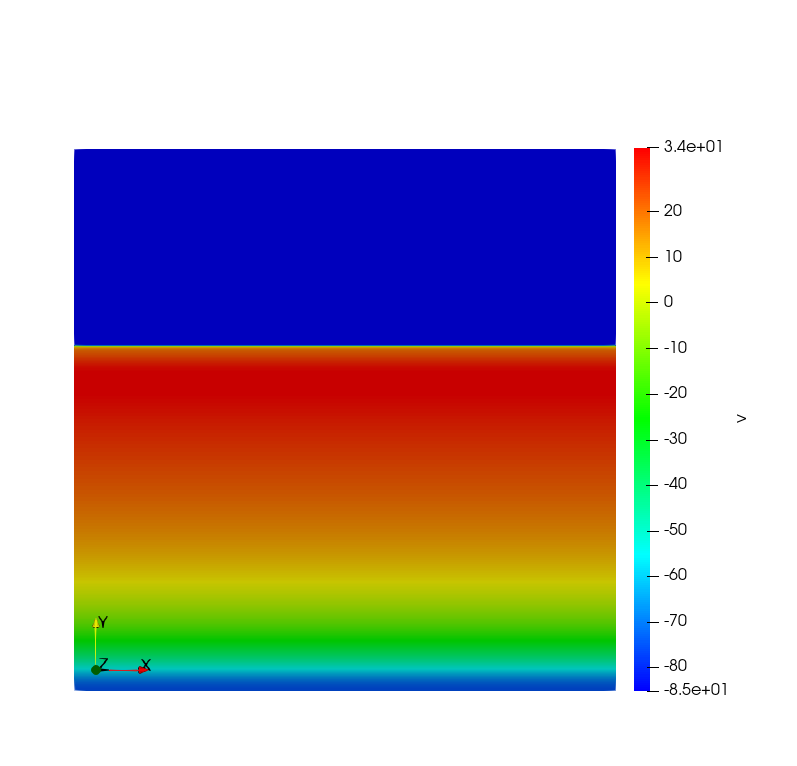}
    \includegraphics[trim={2cm 2.5cm 6cm 4cm},clip,width=0.18\linewidth]{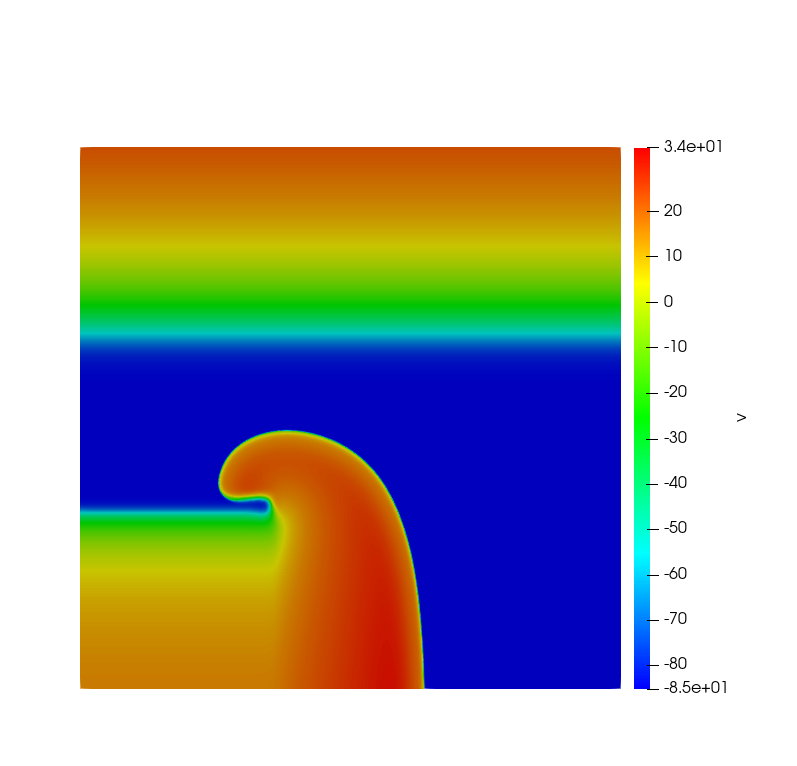}
    \includegraphics[trim={2cm 2.5cm 6cm 4cm},clip,width=0.18\linewidth]{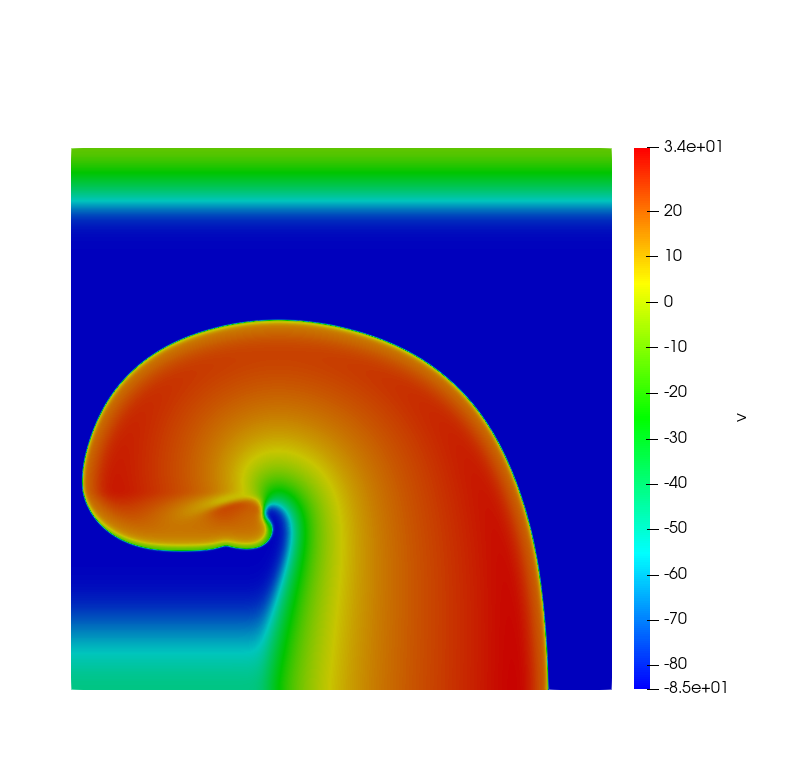}
    \includegraphics[trim={2cm 2.5cm 6cm 4cm},clip,width=0.18\linewidth]{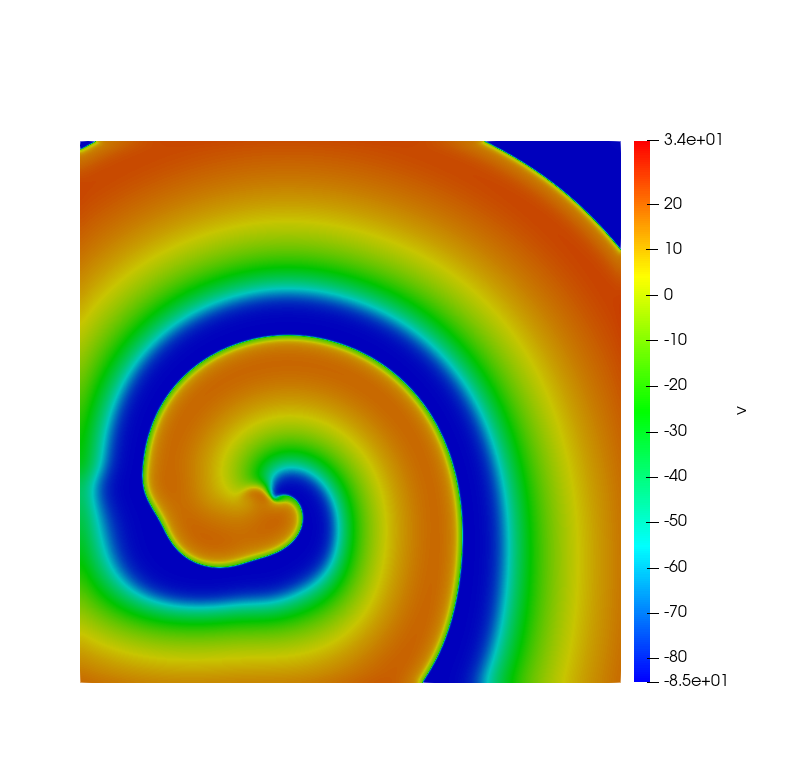}
    \includegraphics[trim={2cm 2.5cm 6cm 4cm},clip,width=0.18\linewidth]{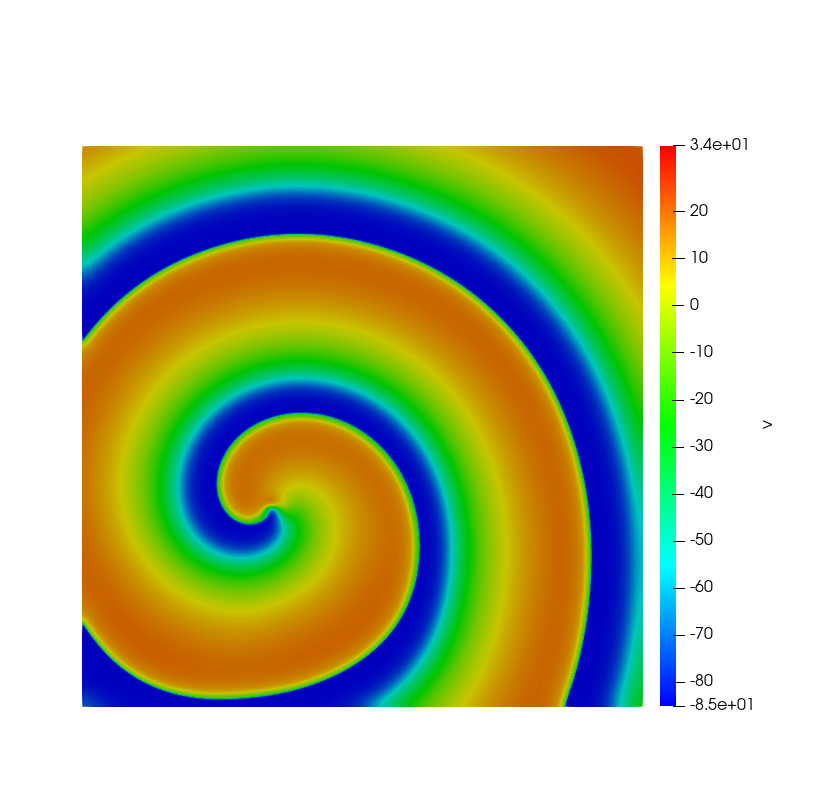}
    \includegraphics[trim={0.5cm 0cm 0cm 0cm},clip,height=1.03in]{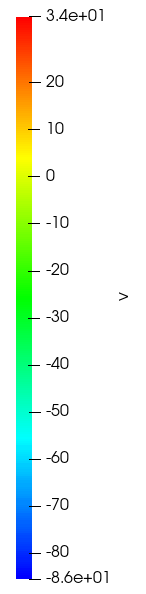}}
\scalebox{0.9}{
     \includegraphics[trim={2cm 2.5cm 6cm 4cm},clip,width=0.18\linewidth]{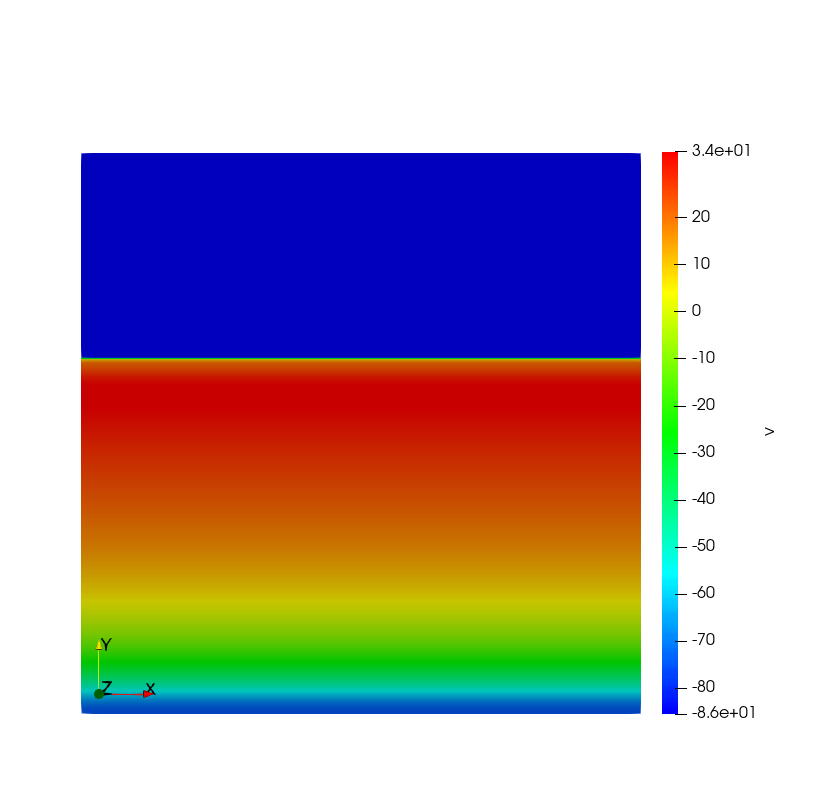}
    \includegraphics[trim={2cm 2.5cm 6cm 4cm},clip,width=0.18\linewidth]{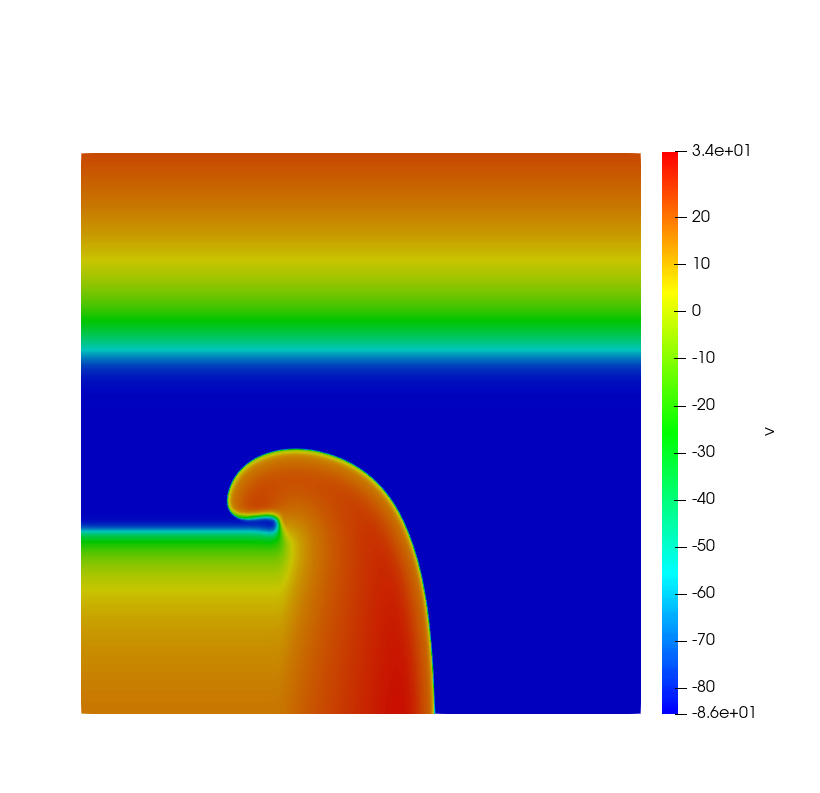}
    \includegraphics[trim={2cm 2.5cm 6cm 4cm},clip,width=0.18\linewidth]{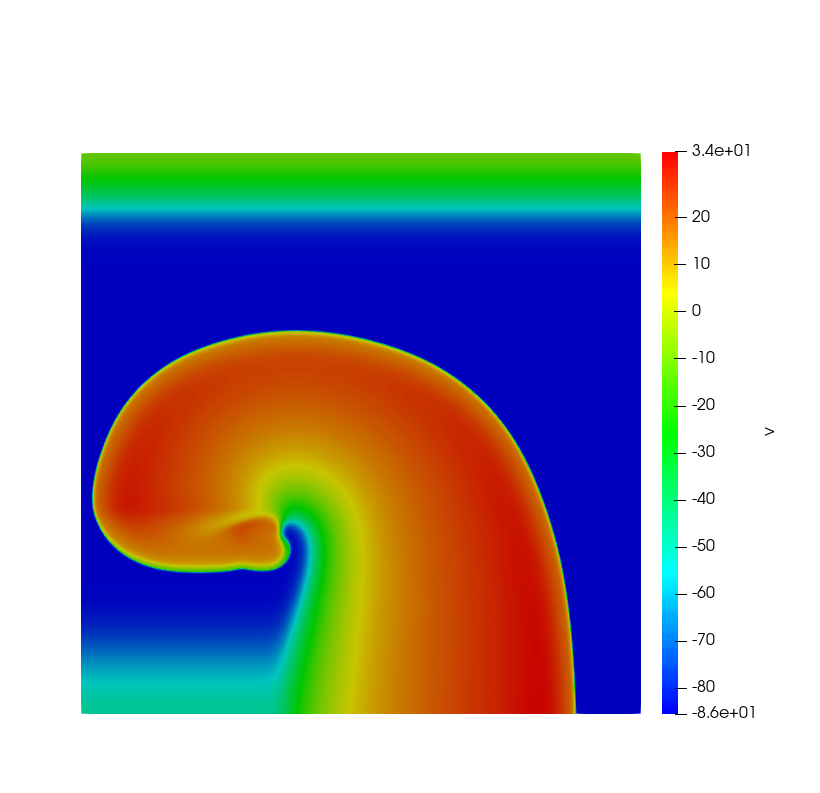}
    \includegraphics[trim={2cm 2.5cm 6cm 4cm},clip,width=0.18\linewidth]{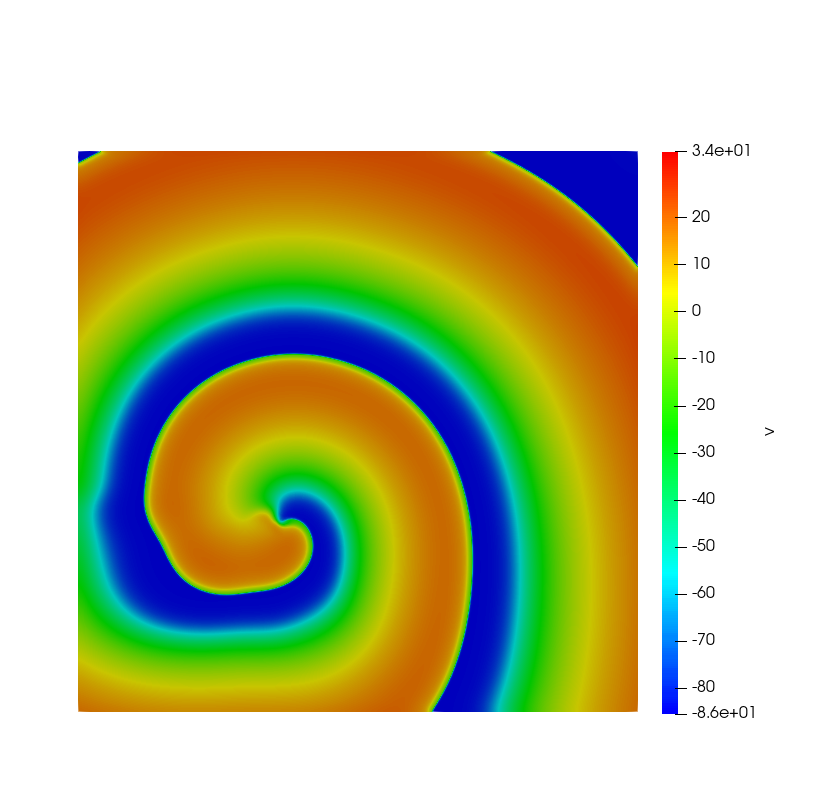}
    \includegraphics[trim={2cm 2.5cm 6cm 4cm},clip,width=0.18\linewidth]{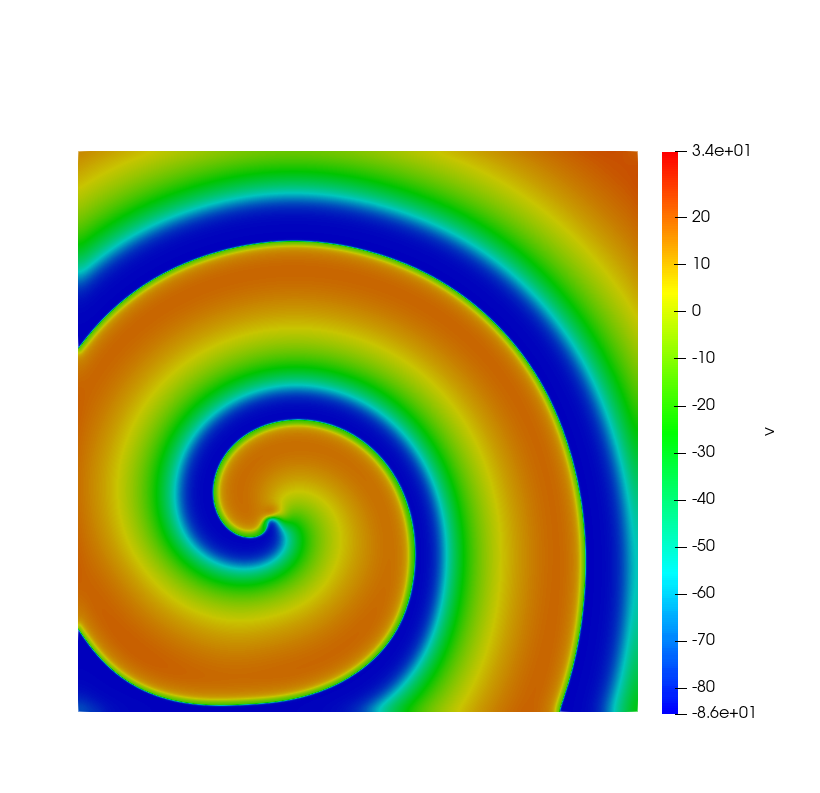}    
    \includegraphics[trim={0.5cm 0cm 0cm 0cm},clip,height=1.03in]{images/tp06/v_m.png}}
\scalebox{0.9}{
    \includegraphics[trim={2cm 2.5cm 6cm 4cm},clip,width=0.18\linewidth]{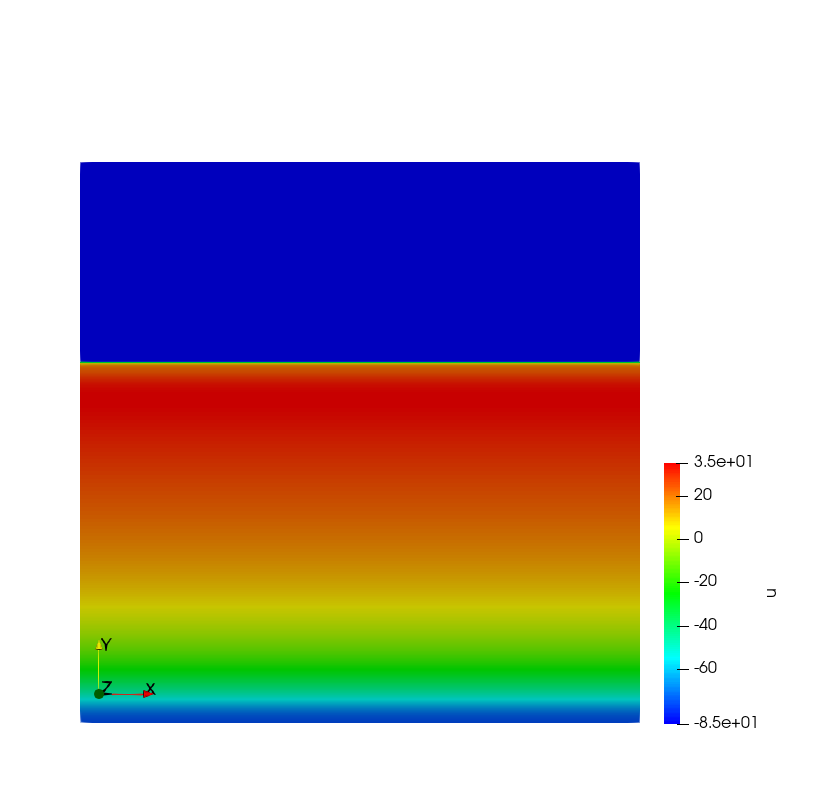}
    \includegraphics[trim={2cm 2.5cm 6cm 4cm},clip,width=0.18\linewidth]{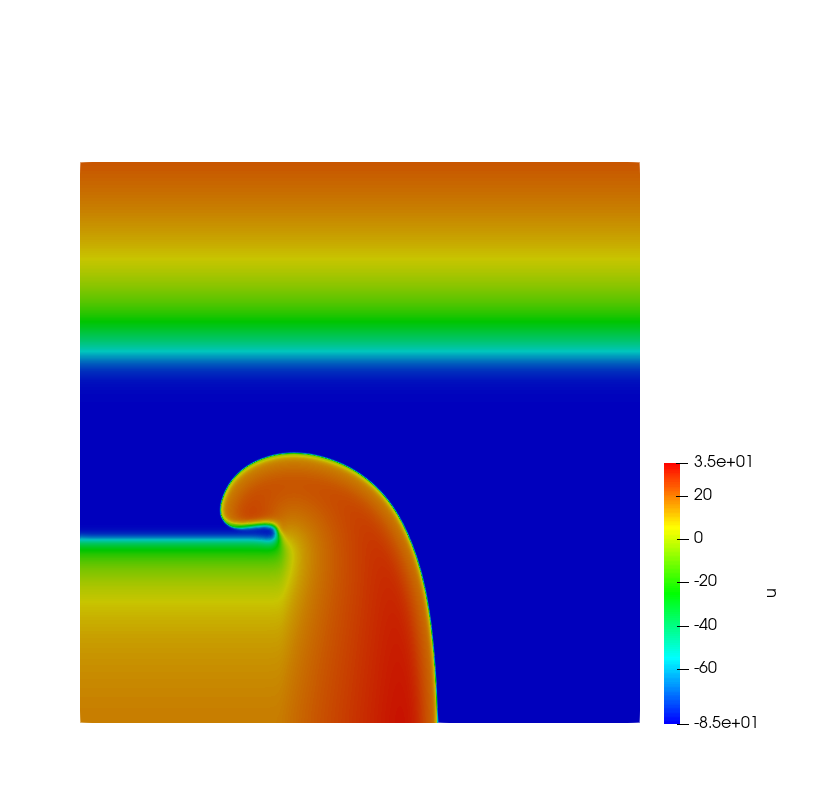}
    \includegraphics[trim={2cm 2.5cm 6cm 4cm},clip,width=0.18\linewidth]{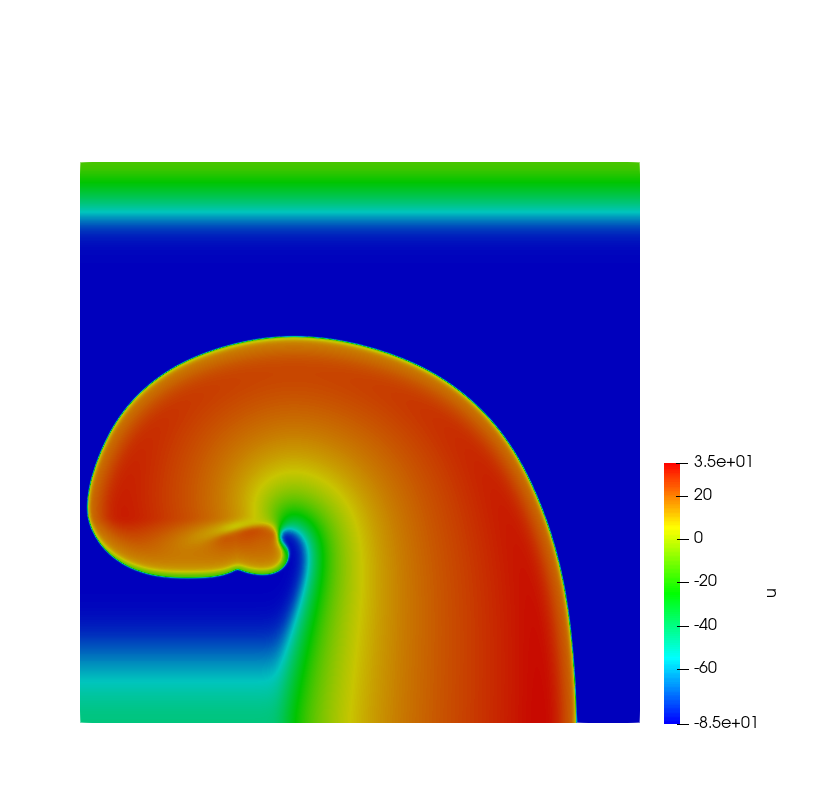}
   \includegraphics[trim={2cm 2.8cm 6cm 4cm},clip,width=0.18\linewidth]{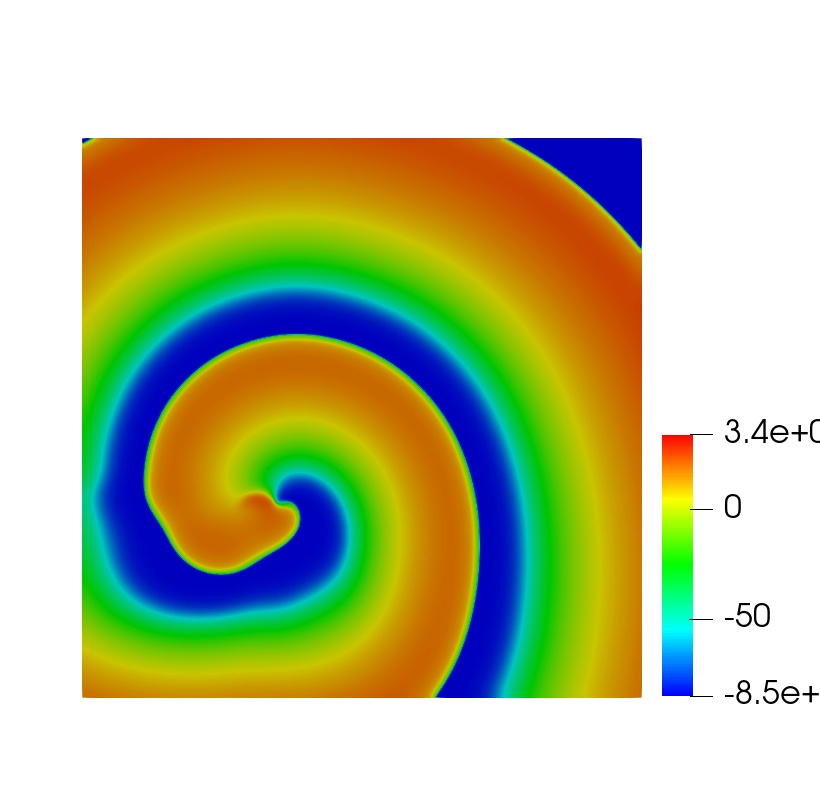}
   \includegraphics[trim={2cm 2.8cm 6cm 4cm},clip,width=0.18\linewidth]{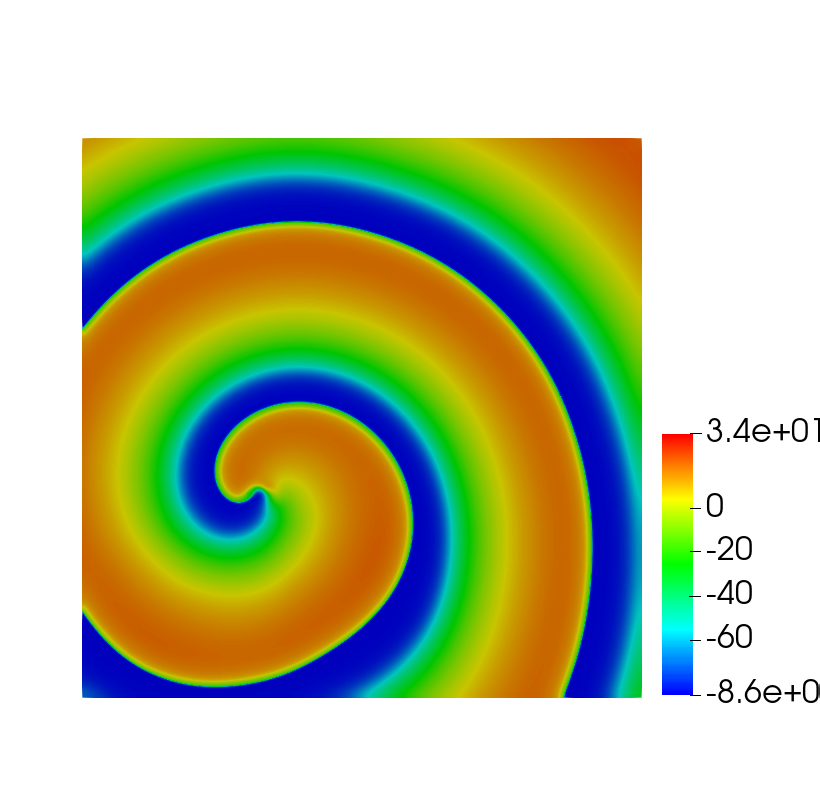}
    \includegraphics[trim={0.5cm 0cm 0cm 0cm},clip,height=1.03in]{images/tp06/v_m.png}}
    \caption{Transmembrane potential $v$ for TP06 model at $t=360$, $710$, $850$, $1510$, and $2000$ ms respectively. \textit{First row:} Fully coupled strategy, \textit{Second row:} Partitioned strategy with SDC, and \textit{Third row:} CN-RK2 scheme with $\tau=0.001$.}
    \label{fig:tp_2D}
\end{figure}

We report the computational times for the three numerical strategies that produced results in close agreement with the fully coupled solution. Due to the substantial computational cost of these simulations, all runs were performed on a single compute node equipped with 28 cores. The fully coupled strategy required approximately 140 hours to complete the simulation over a duration of $600 \text{ ms}$. In contrast, the partitioned strategy with SDC significantly reduced the computational time, completing the simulation in just 19 hours. The decoupled strategy, implemented via the CN-RK2 scheme, completed the same task in approximately 42 hours. This comparison clearly demonstrates the considerable computational savings achieved by the partitioned SDC method, which offers a favorable balance between efficiency and accuracy.
\newline
 % \begin{minipage}[]{0.49\textwidth}
 %     \centering
 %    \includegraphics[trim={2cm 0cm 2cm 1cm},clip,width=\linewidth]{images/lr91/lr_reentry.pdf}
 %   \captionof{figure}{Transmembrane potential for LR91 model over time period $600\,\text{ms}$ at the midpoint of the domain.}
 %    \label{fig:lr91}
 %    \end{minipage}
 %      \hfill
 %     \begin{minipage}[]{0.49\textwidth}
 %  \centering
 %    \includegraphics[trim={2cm 0cm 2cm 1cm},clip,width=\linewidth]{images/tp06/tp_reentry.pdf}
 %    \captionof{figure}{Transmembrane potential $v$ for TP06 model over a period of time $2000\text{ ms}$ at the midpoint of the domain.}
 %    \label{fig:tp06}
 % \end{minipage}
For the TP06 cell model, a computational domain of $\Omega=[0,25]^2$ with $1000\times1000$ cells is considered. An initial stimulus $\I_{stim}=100\,\mu \text{A}/\text{cm}^2$ is applied for $1\, \text {ms}$ at $t=0$ along the bottom edge of the grid. The second stimulus is applied at $t=55\, \text{ms}$ in the region $[0,9]\times[0,9]$ to generate spiral wavefronts.  The transmembrane potential $v$ at the point $(12.5,12.5)$ over a period of time $2000\,\text{ms}$ is presented in Fig.~\ref{fig:tp06}.  Fig.\ref{fig:tp_2D} provides a sequence of snapshots of $v$ at $t = 360$, $710$, $850$, $1510$, and $2000,\text{ms}$, computed using the fully coupled strategy, the partitioned strategy with SDC, and the CN-RK2 scheme.  
The time-stepping configurations employed here are consistent with those used for the LR91 model: the fully coupled strategy uses a fixed time step of $\tau = 0.05$; the partitioned strategy adopts an adaptive time-stepping scheme with a maximum step size of $\tau = 0.05$; and the CN-RK2 scheme is implemented with a fixed time step of $\tau = 0.001$. The corresponding simulation times on two compute nodes (each with 28 cores) were approximately 1930 hours for the fully coupled approach, 230 hours for the partitioned strategy with SDC, and 397 hours for the CN-RK2 scheme. These results highlight the significant reduction in computational time achieved by the proposed partitioned strategy with SDC while maintaining high accuracy.
%965 for fully coupled using 5 nodes
\subsection{Results on 3D slab}
\begin{figure}
    \centering
     \includegraphics[trim={2cm 0cm 4.8cm 10cm},clip,width=0.175\linewidth]{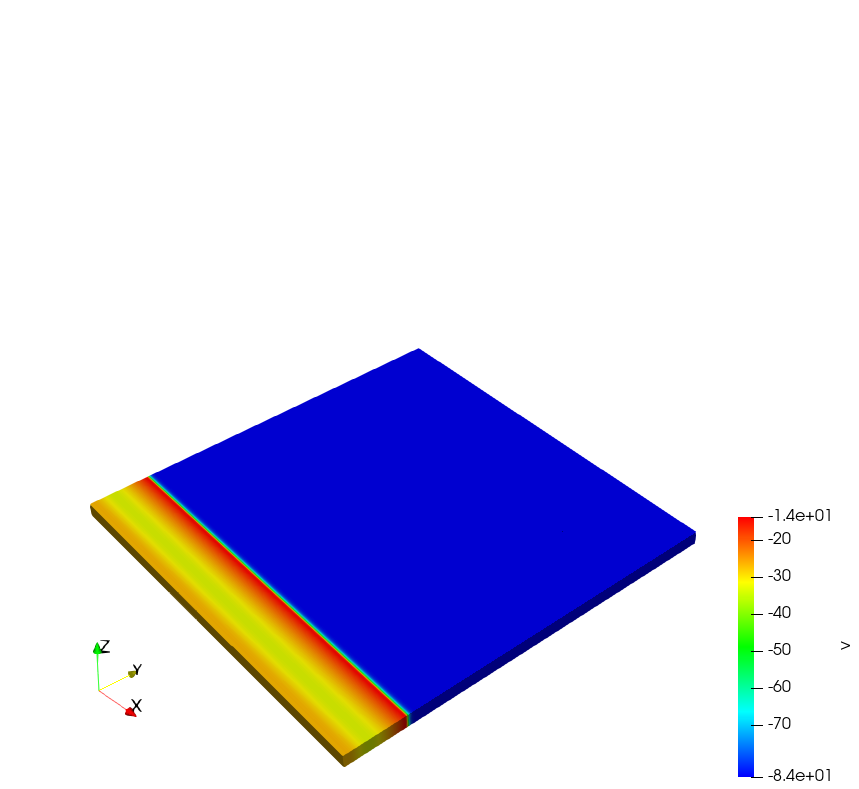}
    \includegraphics[trim={2cm 0cm 4.8cm 10cm},clip,width=0.175\linewidth]{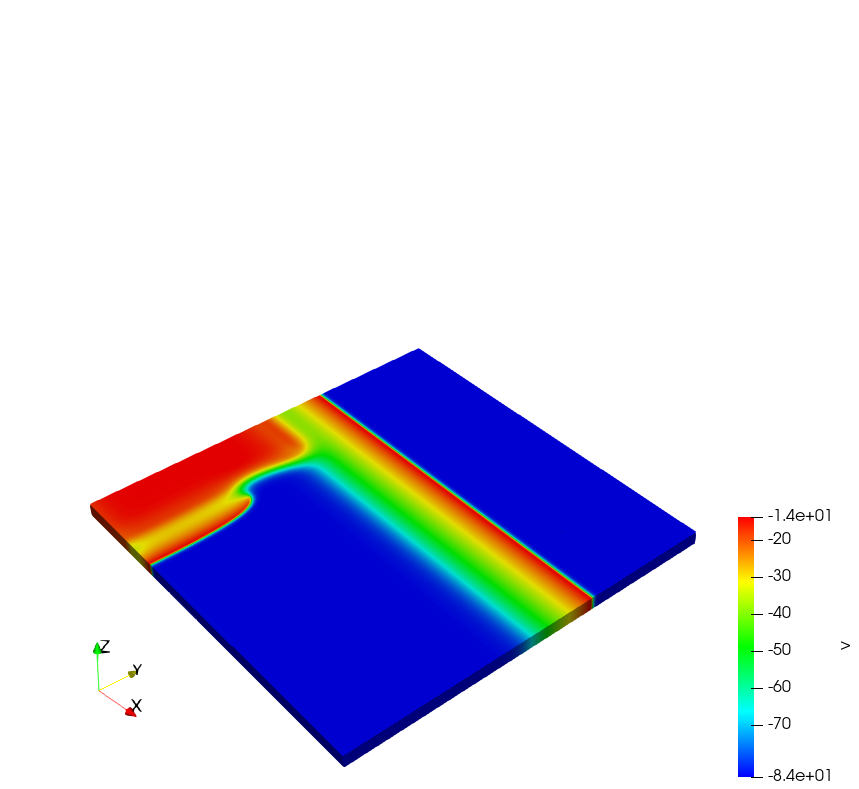}
    \includegraphics[trim={2cm 0cm 4.8cm 10cm},clip,width=0.175\linewidth]{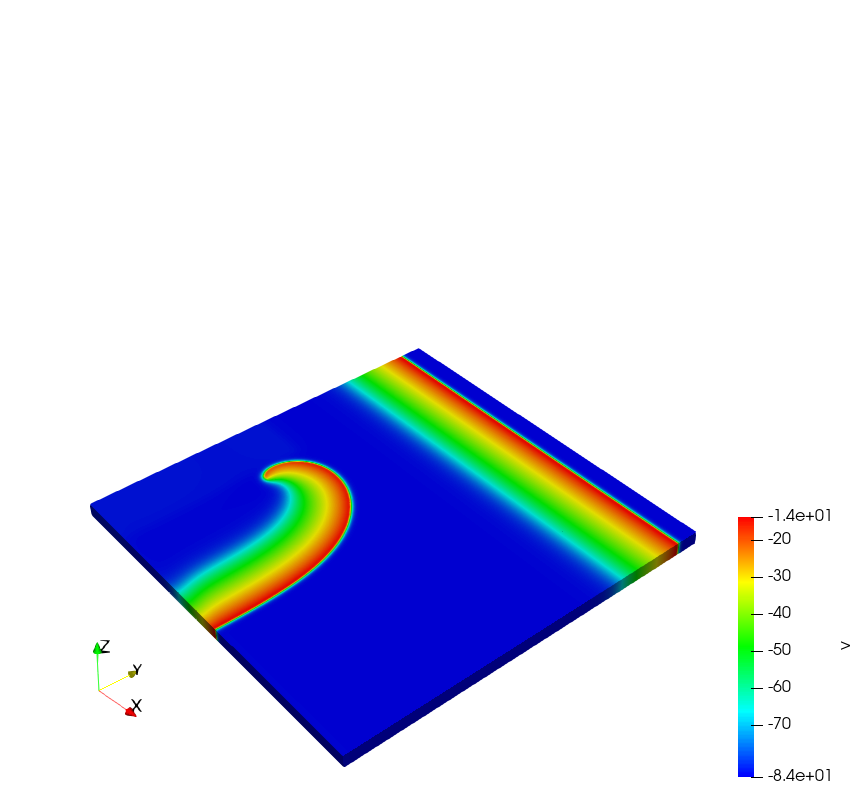}
    \includegraphics[trim={2cm 0cm 4.8cm 10cm},clip,width=0.175\linewidth]{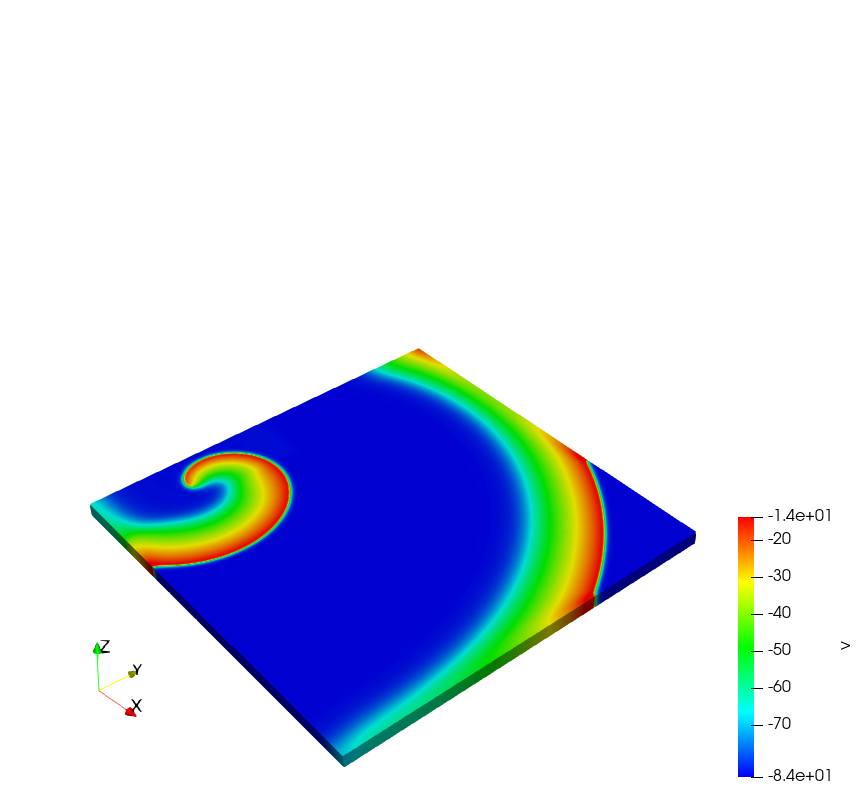}
     \includegraphics[trim={2cm 0cm 4.8cm 10cm},clip,width=0.175\linewidth]{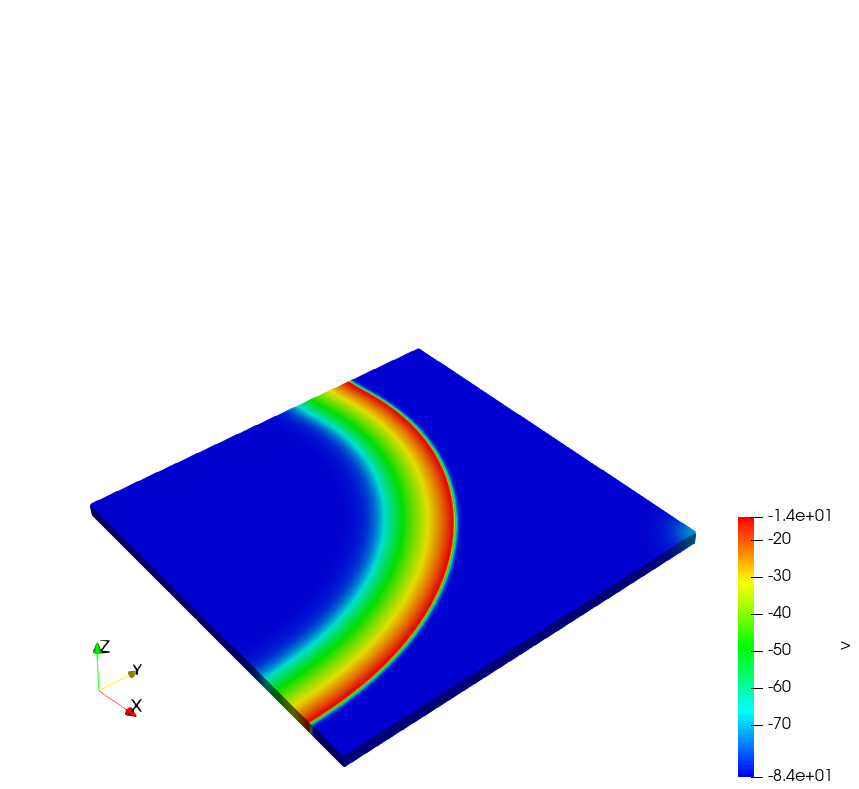}
    \includegraphics[trim={0.5cm 0cm 0.5cm 0cm},clip,height=0.7in]{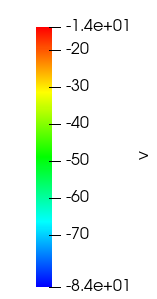}
    
    \includegraphics[trim={2cm 0cm 4.8cm 10cm},clip,width=0.175\linewidth]{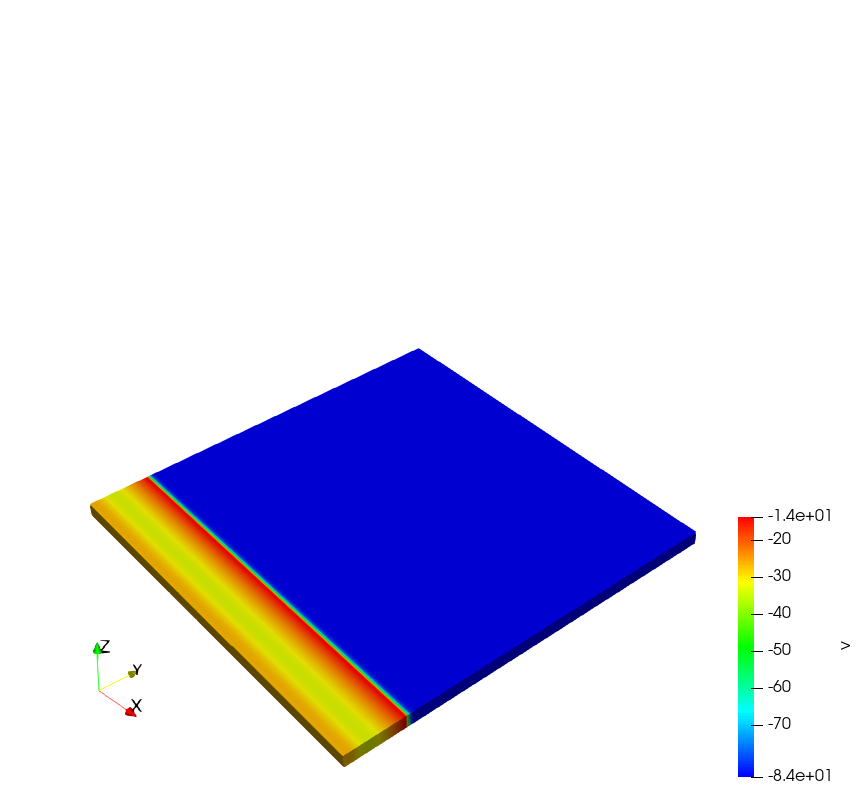}
    \includegraphics[trim={2cm 0cm 4.8cm 10cm},clip,width=0.175\linewidth]{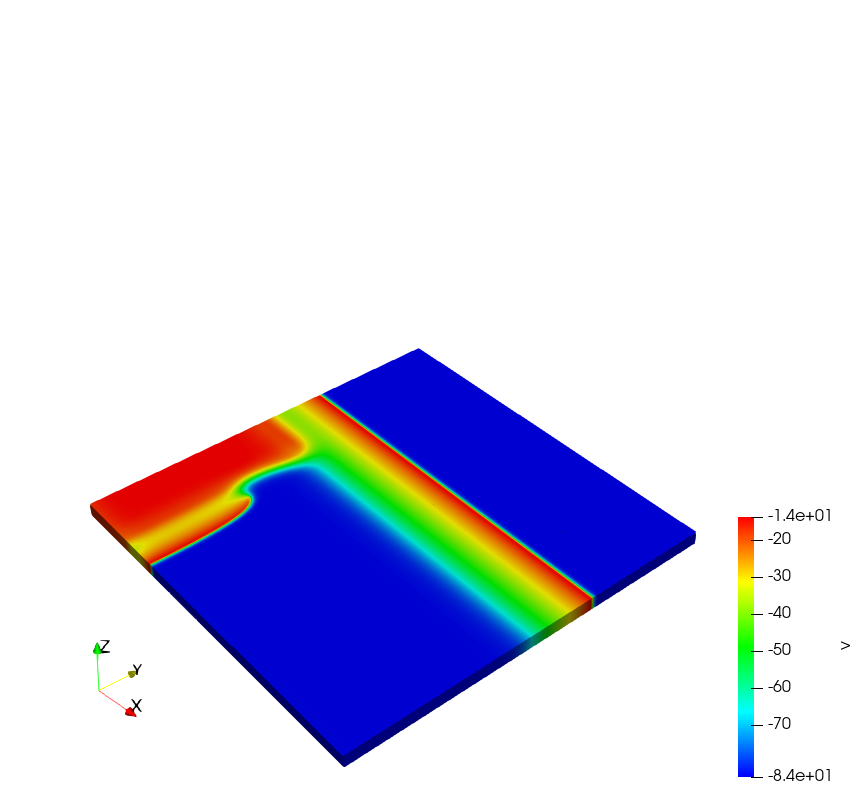}
    \includegraphics[trim={2cm 0cm 4.8cm 10cm},clip,width=0.175\linewidth]{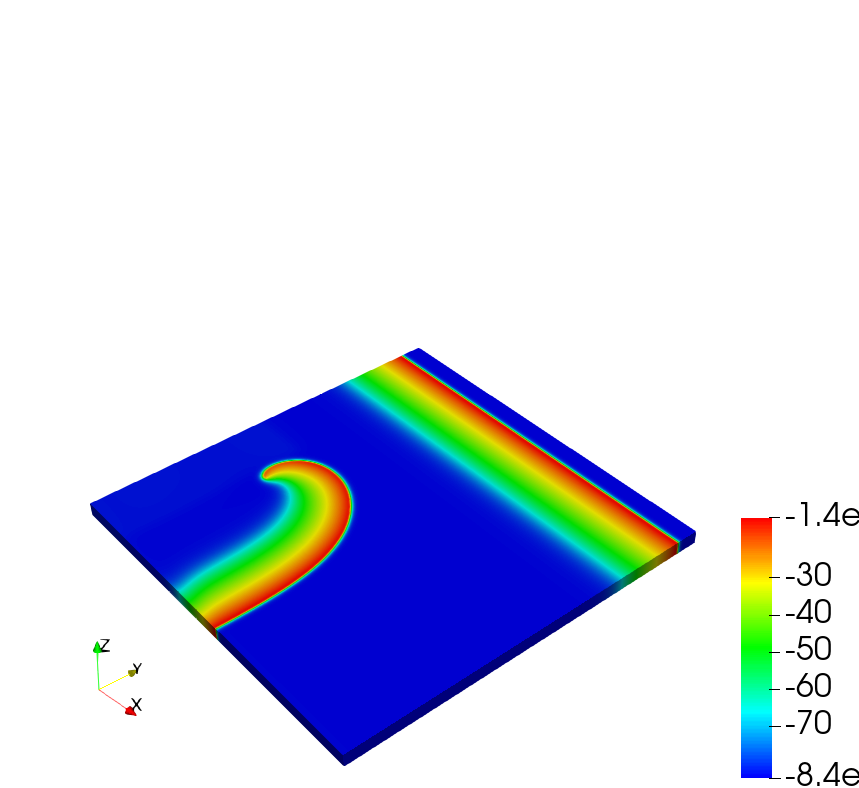}
    \includegraphics[trim={2cm 0cm 4.8cm 10cm},clip,width=0.175\linewidth]{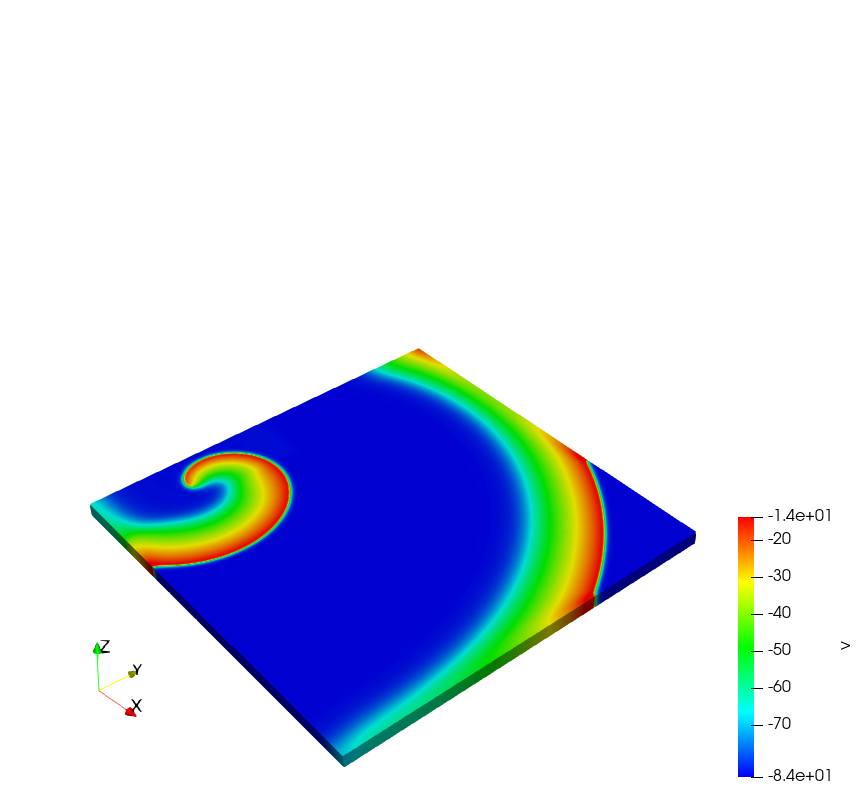}
    \includegraphics[trim={2cm 0cm 4.8cm 10cm},clip,width=0.175\linewidth]{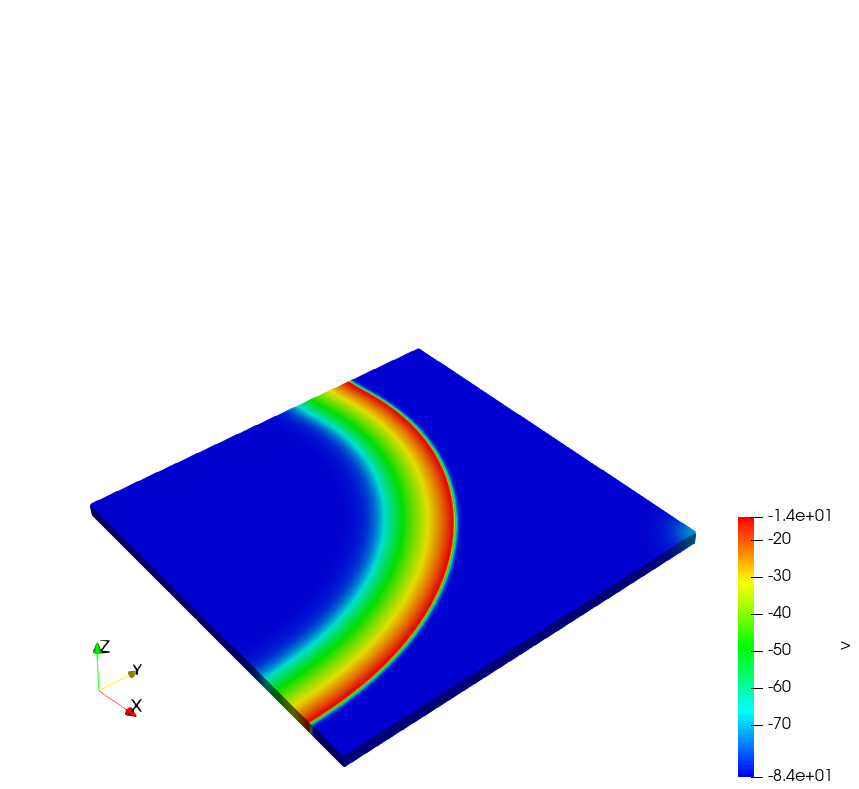}
    \includegraphics[trim={0cm 0cm 0cm 0cm},clip,height=0.7in]{images/3d_cubic/v_scale.png}
    \caption{Transmembrane potential $v$ for LR91 model at $t=32$, $146$, $206$, $340$, and $600$ ms respectively on a 3D slab domain. \textit{First row:} Partitioned strategy with SDC using $\tau=0.05$, \textit{Second row:} CN-RK2 scheme using $\tau=0.001$.}
    \label{fig:lr_3D}
\end{figure}

The numerical results for the LR91 model on a three-dimensional slab domain of size $[0,6]\times[0,6]\times[0,0.2]\,\text{cm}^3$ are presented in this subsection to evaluate the performance of the partitioned strategy with SDC and the decoupled approach. The computational domain is discretized using a structured grid with $600\times60\times20$ cells. An initial stimulus and a second stimulus are applied at the same locations and time as in the two-dimensional case to generate and sustain reentrant electrical activity.
Figure~\ref{fig:lr_3D} shows the evolution of the transmembrane potential at time $t = 32$, $146$, $206$, $340$, and $600\,\text{ms}$ on the three-dimensional slab geometry. The partitioned strategy with SDC and the CN-RK2 scheme produce comparable spatial patterns, accurately capturing the propagation of transmembrane potential. While the two approaches yield similar results, the partitioned strategy completes the simulation significantly faster, requiring approximately 432 hours compared to 647 hours for the CN-RK2 scheme, highlighting its computational efficiency.
\subsection{Results on Realistic Geometry}
The numerical results for the partitioned strategy with SDC and the decoupled strategy for the LR91 model on a realistic three-dimensional geometry are presented in this subsection. A full rabbit ventricle geometry, which is generated based on histological images \cite{plank09:_generation}, is considered for the computations. The computational domain measures $2.78\times2.92\times2.58\,\text{cm}^3$, comprising 3,073,529 tetrahedra and 656,421 nodal points. For the construction of parallel grids, the ALUGrid library \cite{Dedner_ALUGrid_CVS04} was employed, which in turn uses METIS \cite{Metis_SISC98} graph partitioner for the decomposition of the grid.

\begin{figure}
    \centering
    \includegraphics[trim={1cm 0 3.5cm 4cm},clip,width=0.175\linewidth]{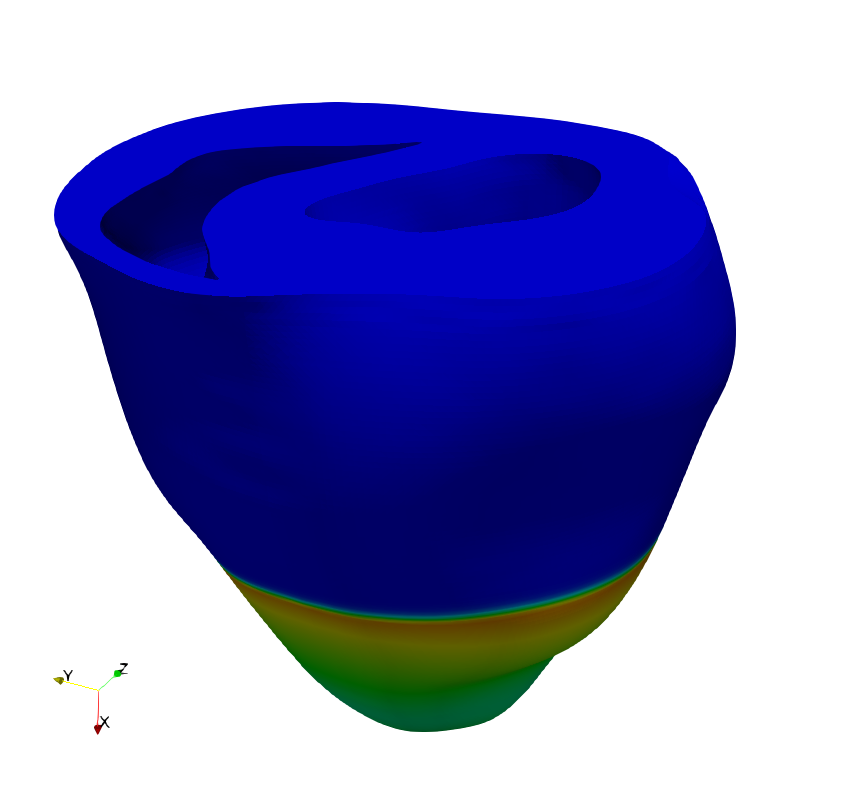}
    \includegraphics[trim={1cm 0 3.5cm 4cm},clip,width=0.175\linewidth]{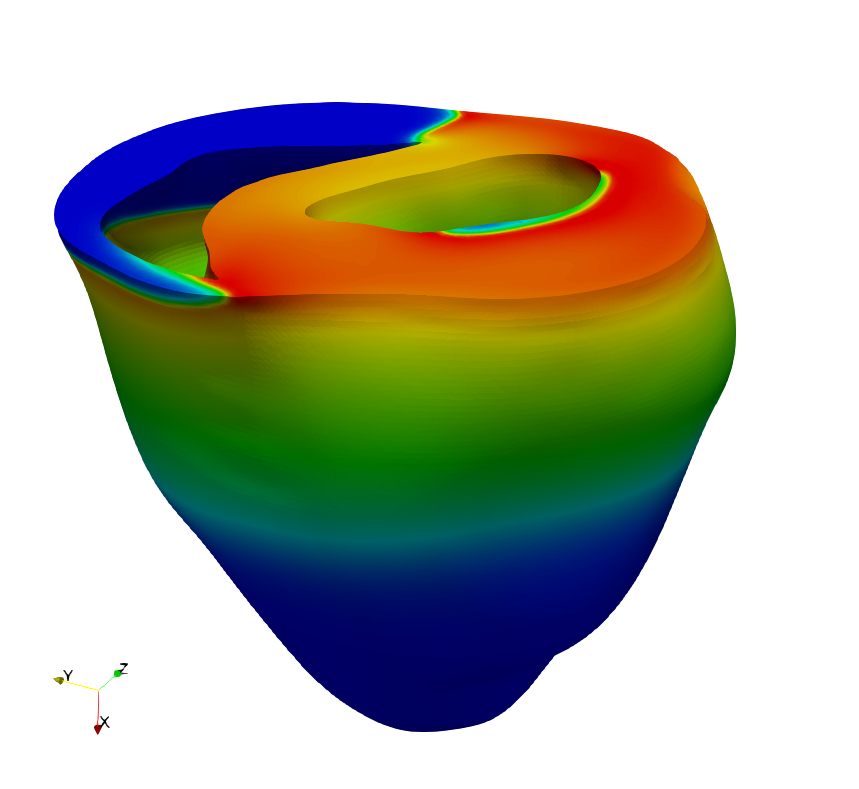}
    \includegraphics[trim={1cm 0 3.5cm 4cm},clip,width=0.175\linewidth]{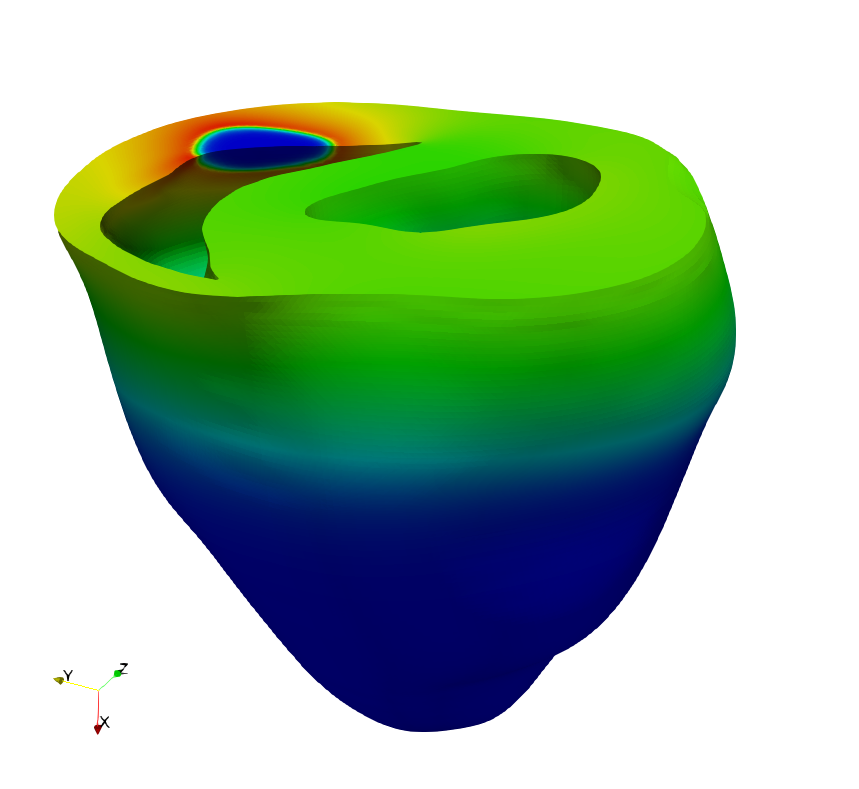}
    \includegraphics[trim={1cm 0 3.5cm 4cm},clip,width=0.175\linewidth]{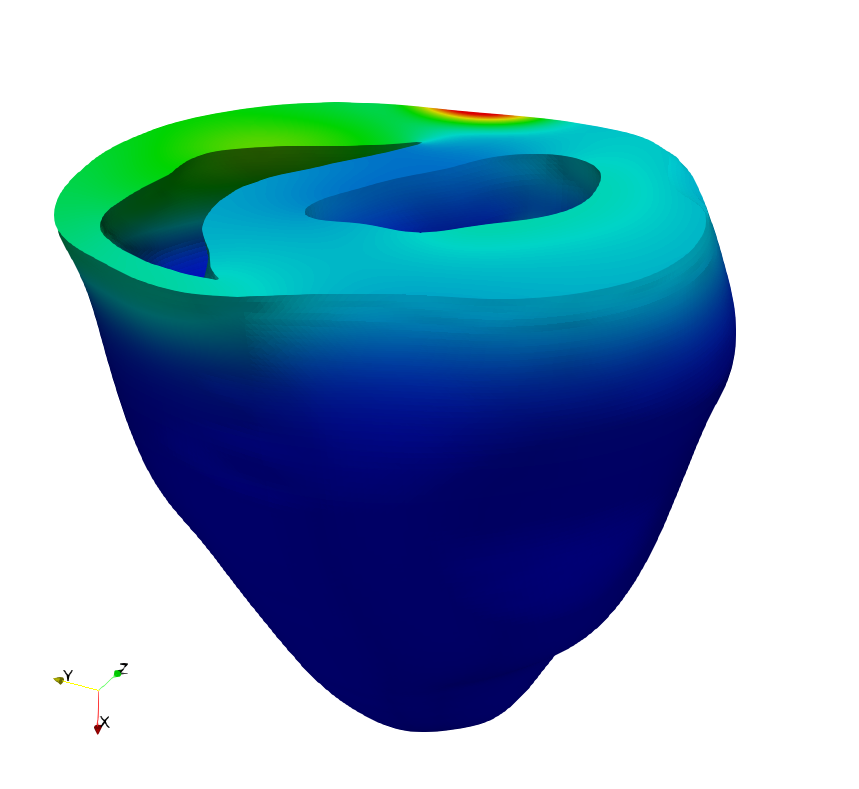}
    \includegraphics[trim={1cm 0 3.5cm 4cm},clip,width=0.175\linewidth]{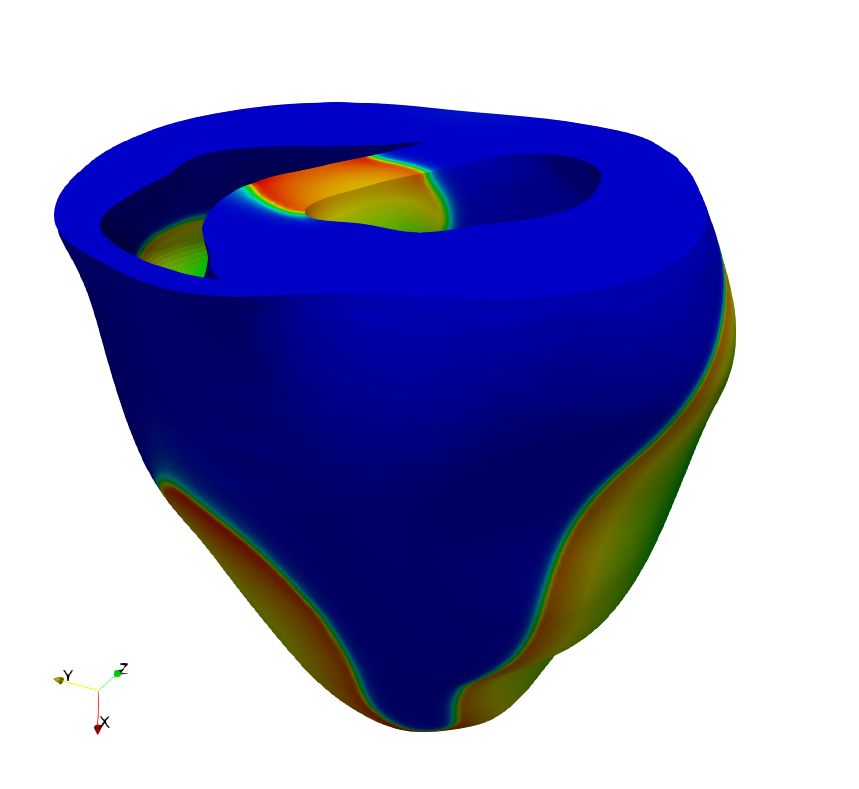}
   \includegraphics[trim={0cm 0cm 0cm 0cm},clip,height=0.9in]{images/lr91/re-entry/v_m.png}

    \includegraphics[trim={1cm 0 3.5cm 4cm},clip,width=0.175\linewidth]{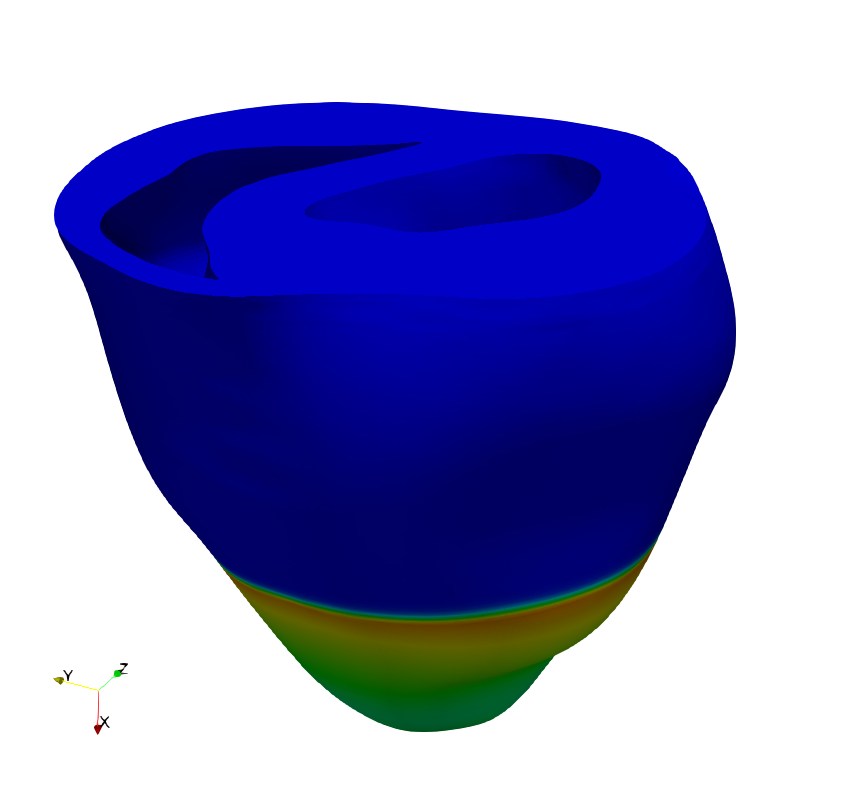}
    \includegraphics[trim={1cm 0 3.5cm 4cm},clip,width=0.175\linewidth]{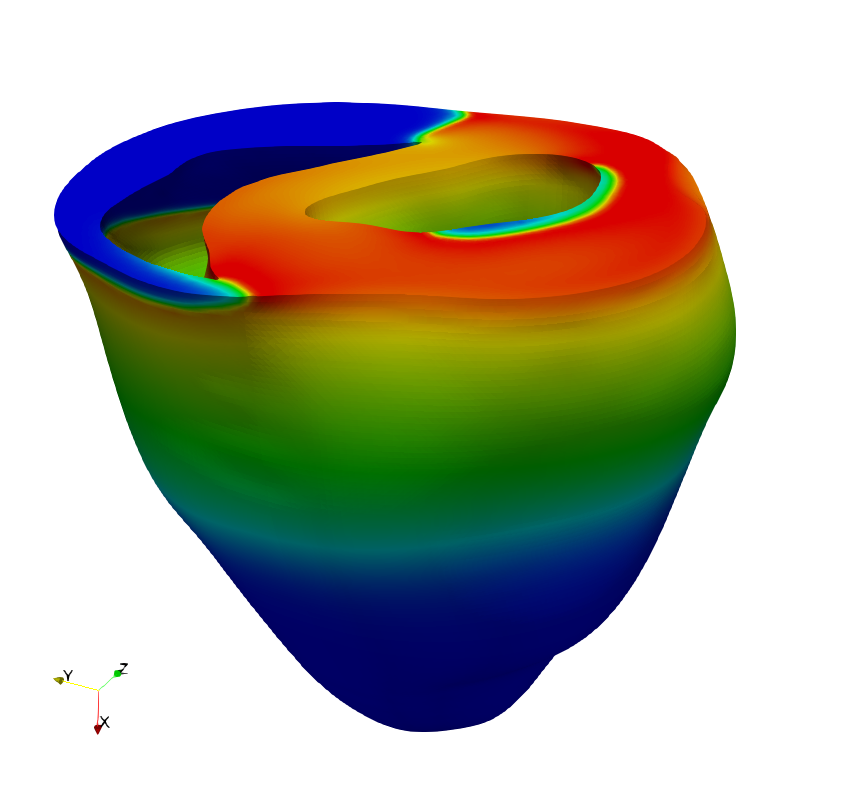}
    \includegraphics[trim={1cm 0 3.5cm 4cm},clip,width=0.175\linewidth]{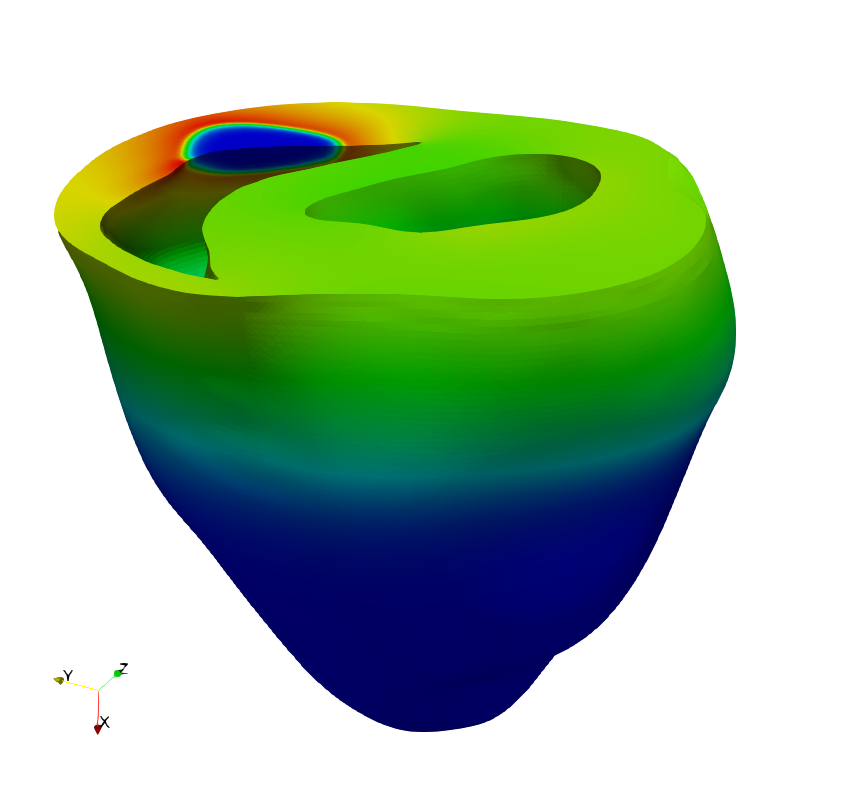}
     \includegraphics[trim={1cm 0 3.5cm 4cm},clip,width=0.175\linewidth]{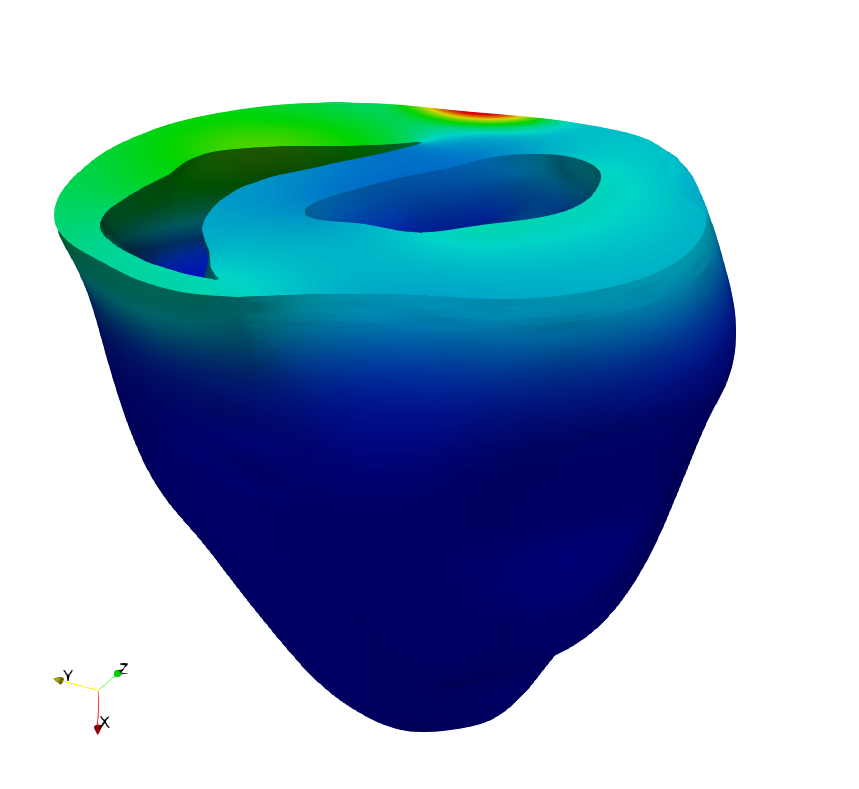}
     \includegraphics[trim={1cm 0 3.5cm 4cm},clip,width=0.175\linewidth]{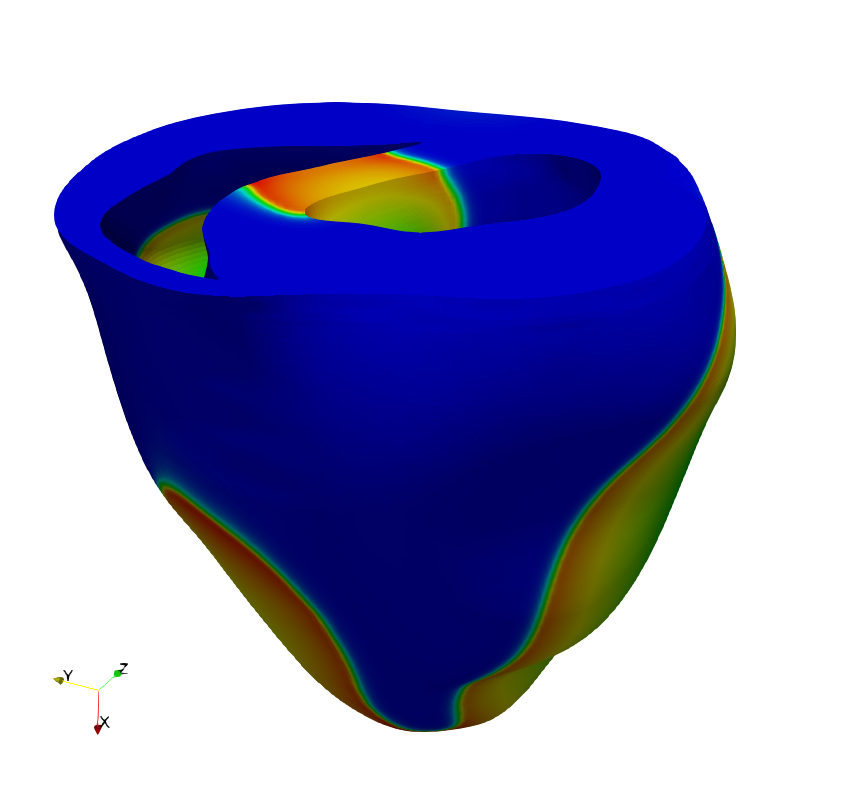}
    \includegraphics[trim={0cm 0cm 0cm 0cm},clip,height=0.9in]{images/lr91/re-entry/v_m.png}
    \caption{Transmembrane potential $v$ for LR91 cell model at $t=50$, $106$, $120$, $140$, and $200\text{ ms}$ on a realistic rabbit ventricular geometry. \textit{First row}: Partitioned strategy with SDC with time step $\tau=0.05$. \textit{Second row}: CN-RK2 scheme with time step $\tau=0.001$.}
    \label{fig:lr_geo}
\end{figure}

A stimulus of $100~\mu A/\text{cm}^3$ is applied initially at $t=0\text{ ms}$ at the bottom of the geometry for a duration of $2\text{ ms}$, similar to the case for the 2D domain. The second stimulus is applied at $t=135\text{ ms}$ for $2\text{ ms}$ in a small region of radius $0.8\text{ cm}$ from the point $(0.0281,0.5744,1.1987)$. Figure~\ref{fig:lr_geo} presents a sequence of snapshots at time $t=50$,  $106$, $120$, $140$ and $200\,\text{ms}$ for the partitioned strategy with SDC using adaptive time-stepping approach with a maximum time step of $0.05$, and CN-RK2 scheme using fixed time step size $\tau=0.001$. The two strategies yield comparable numerical results. However, regarding computational efficiency, the partitioned strategy completed the simulation in approximately 71 hours, while the CN-RK2 scheme required nearly 287 hours to achieve the same.

\subsection{Results for Bidomain-Bath}
In this subsection, we consider the bidomain-bath configuration, where the electrical activity of the cardiac tissue is embedded within a surrounding nonconductive bath. The integrated computational domain is $[-1,7] \times [-1,7]\subset \mathbb{R}^2$, which is discretized using a uniform grid of $800 \times 800$ cells. The embedded tissue domain is $[0,6]\times[0,6]$ with a uniform grid of $600\times600$ cells.
The partitioned strategy uses an adaptive time-stepping approach with a maximum time step of $0.05$, and the decoupled strategy employs a fixed time step of $0.001$, consistent with that used in the continuum (tissue) domain. 

\begin{figure}
    \centering
    \includegraphics[trim={2cm 0cm 2cm 1cm},clip,scale=0.2]{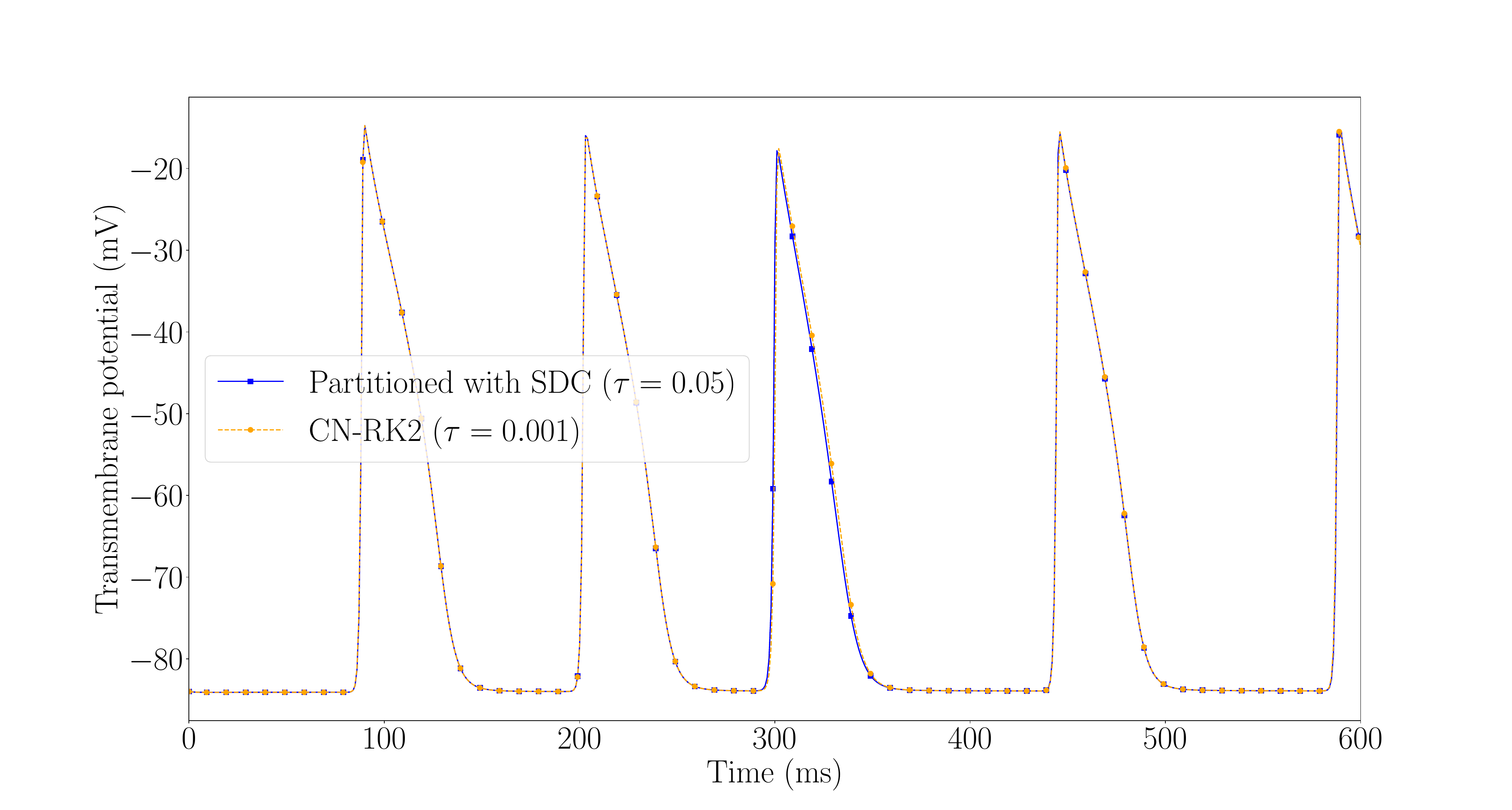}
    \caption{Transmembrane potential over a period of $600\,\text{ms}$ at the midpoint of the tissue domain.}
    \label{fig:bath_2d}
\end{figure}

Figure \ref{fig:bath_2d} presents the transmembrane potential $v$ over a simulation period of $600\,\text{ms}$ at the point $(3,3)$ on the tissue domain. The results from the partitioned strategy with SDC and the CN-RK2 scheme are nearly identical. The computational simulations were performed on two compute nodes, each with 28 cores. Regarding computational cost, the partitioned strategy with SDC completed the simulation in approximately 34 hours, whereas the CN-RK2 scheme required about 65 hours. This substantial reduction in computation time demonstrates the efficiency of the proposed partitioned strategy with SDC in the bidomain-bath setting, offering a favorable balance between accuracy and computational cost.

\section{Conclusion}
In this work, we have introduced a novel partitioned bidomain model strategy incorporating a spectral deferred correction (SDC) approach. This strategy improves accuracy while preserving a flexible and independent treatment of the parabolic-ODE and elliptic components, allowing for easier integration and modification of individual sub-solvers. The unique aspect of this scheme is its ability to enable high-order time integration without sacrificing the computational advantages of operator splitting. Comprehensive numerical experiments were performed using the  LR91 and TP06 realistic ionic models across various test cases, including an idealized two-dimensional domain, realistic three-dimensional ventricular geometry, and an integrated bidomain-bath configuration. These tests demonstrated that the proposed partitioned strategy with SDC achieves accuracy comparable to the fully coupled strategy while incurring significantly lower computational costs. In particular, the CN-RK2 scheme and the SDC-enhanced algorithm exhibit superior accuracy among the decoupled and partitioned strategies, aligning closely with the reference fully coupled solution in transient behaviour and error metrics.
Importantly, the proposed partitioned strategy with SDC significantly reduced computational time using CTSM, achieving more than seven times faster than the fully coupled strategy while maintaining high accuracy. The strategy was also effectively used to simulate phenomena like spiral wave reentry, accurately capturing the spatiotemporal evolution of transmembrane potential. These results demonstrate the effectiveness of the proposed partitioned scheme in achieving a balance between accuracy and efficiency, highlighting its potential as a robust and scalable framework for large-scale cardiac simulations.

\bibliographystyle{siamplain}
 \bibliography{ms}
\end{document}